\documentclass[reqno,11pt]{amsart}
\usepackage{amsthm}
\usepackage{amsmath,amsfonts,amssymb,bm}
\usepackage{xr}
\usepackage{mathrsfs}
\usepackage{mathtools}
\usepackage{hyperref}
\usepackage{mathabx}
\usepackage[initials,nobysame]{amsrefs}
\newtheorem{theorem}{Theorem}[section]
\newtheorem{lemma}[theorem]{Lemma}
\newtheorem{proposition}[theorem]{Proposition}
\newtheorem{corollary}[theorem]{Corollary}
\newtheorem{remark}[theorem]{Remark}
\usepackage[top=2.0cm,bottom=2.0cm,left=3cm,right=3cm]{geometry}
\numberwithin{equation}{section}
\mathtoolsset{showonlyrefs}

\newcommand{\eps}{\varepsilon}
\newcommand{\p}{\partial}
\newcommand{\R}{\mathbb{R}} 
\def\L{\mathcal{L}}

\begin{document}
\begin{sloppypar}
\title[Compressible Navier-Stokes approximation for the Boltzmann Equation]{ Compressible Navier-Stokes system with slip boundary from Boltzmann equations with reflection boundary: derivations and justifications}
\date{}
\author[N. Jiang and Y.-L. Wu]{Ning Jiang \and Yulong Wu}

\address[N. Jiang]{School of Mathematics and Statistics, Wuhan University, Wuhan 430072, China
}

\email{
njiang@whu.edu.cn}
\address[Y.-L. Wu]{
  School of Mathematics and Statistics,
    Wuhan University,  Wuhan 430072, China
}
  \email{
  yulong\string_wu@whu.edu.cn
}

\begin{abstract}

This is the first in a series of papers connecting the boundary conditions for the compressible Navier-Stokes system from the Boltzmann equations with the Maxwell reflection boundary. The slip boundary conditions are formally derived from the Boltzmann equation with both specular and almost specular reflection boundary conditions. That is, the accommodation coefficient $\alpha_\eps=O(\eps^\beta)$ with $\beta>0$ or $\alpha_\eps =0$. Here, the small number $\eps>0$ denotes the Knudsen number. The systematic formal analysis is based on the Chapman-Enskog expansion and the analysis of the Knudsen layer. In particular, for the first time, we employ the appropriate ansatz for the general $\beta>0$. This completes the program started in \cite{aoki2017slip}. In the second part, the compressible Navier-Stokes-Fourier approximation for the Boltzmann equation with specular reflection in general bounded domains is rigorously justified. The uniform regularity for the compressible Navier-Stokes system with the derived boundary conditions is investigated. For the remainder equation,  the $L^2\mbox{-}L^6\mbox{-}L^\infty$ framework is employed to obtain uniform estimates in $\eps$. \\


  \end{abstract}
    \maketitle


\section{Introduction}
\subsection{Boltzmann equations and fluid limits}
In this paper, the rescaled Boltzmann equation in a three dimensional smooth bounded domain $\Omega$ is considered:
\begin{equation}
\begin{cases}\label{be}
\partial_tF_\varepsilon+v\cdot\nabla_xF_\varepsilon=\frac 1\varepsilon Q(F_\varepsilon,F_\varepsilon)  &\text{on} \quad \mathbb{R}_+\times\Omega\times\mathbb{R}^3, \\
F_\varepsilon(0,x,v)=F^{in}_\varepsilon(x,v)\geq 0     &\textrm{on}\quad \Omega\times\mathbb{R}^3,\\
\end{cases}
\end{equation}
where $F_\varepsilon=F_\varepsilon(t,x,v) \geq0$ denotes the density distribution function of the particle gas at time $t \geq0$, position $x\in \Omega$ and velocity $v\in \mathbb{R}^3$. The positive parameter $\varepsilon$ is the Knudsen number defined as the ratio of the molecular mean free path length to a representative physical length scale. Furthermore, $Q(F_1,F_2)$ denotes the Boltzmann collision operator and is defined as
\begin{equation}\label{coll}
Q(F_1,F_2)=\int_{\mathbb{R}^3 \times \mathbb{S}^2}b(|v_1-v|,\omega )[F_1(v_1^{\prime})F_2(v^{\prime})-F_1(v_1)F_2(v)]dv_1 d\omega,
\end{equation}
where $v,v_1 $ and $v',v'_1 $ denotes the velocities of the particles before and after the collision respectively, which satisfy the conservation of momentum and energy for partial pairs after a collision, namely,
$$\left\{
\begin{aligned}
&v'=v+[(v_1 -v)\cdot w]w,\quad v_1'=v_1-[(v_1 -v)\cdot w]w,\\
&|v_1 |^2+|v|^2=|v'_1 |^2+|v'|^2.
\end{aligned}\right.
$$
for any given $(\omega,v_1,v)\in \mathbb{S}^2\times\mathbb{R}^3\times\mathbb{R}^3$. The non-negative function $b(|v_1 - v|, \omega)$ represents the collision kernel, which is determined by the intermolecular potential. For hard-sphere molecules, $b(|v_1 - v|, \omega) = |(v_1 - v) \cdot \omega|$.
In the past three decades, significant progress has been made in justifying the hydrodynamic limits from the Boltzmann equations. The Knudsen number (denoted as $\eps>0$ in the current paper), which is the ratio of the mean free path to the macroscopic length scale, characterizes the fluidity of dilute gases. The smaller the Knudsen number, the more the gas behaves like fluids. Under different physical scalings, different fluid equations, such as compressible (or incompressible) Navier-Stokes (or Euler) equations, can be systematically derived from the scaled Boltzmann equations. For the formal derivations of the fluid equations (without boundary conditions), see \cite{bardos1991fluid}. The most typical methods of derivations are moment methods, and expansion methods. In particular, the expansion methods include the Hilbert and Chapman-Enskog expansion.

There have been many notable advances in the convergence of these fluid limits. Based on the limiting fluid equations, the convergence results can be classified into the following four groups:

(1) \underline{Incompressible Navier-Stokes limits}. One of the most successful progress was the Bardos-Golse-Levermore program starting from the late 80s, aimed at justifying the limit to Leray's solutions of the incompressible Navier-Stokes system from the DiPerna-Lions renormalized solutions of Boltzmann equations. This program was initialized by Bardos-Golse-Levermore \cites{bardos1993fluid,bardos1991fluid} and completed by Golse and Saint-Raymond \cites{golse2004navier,golse2009incompressible}. The key issues are compactness (from entropy inequality) and weak convergence methods. For the case of domains with boundary, from the Boltzmann equations with Maxwell reflection boundary, see also \cites{Masmoudi-SRM-CPAM2003, jiang2017boundary}. From the incoming boundary condition, see \cite{Jiang-ZhangX-JDE2019, Jiang-ZhangX-SIMA2019}. Another types of results are in the framework of classical solutions, which are based on the nonlinear energy method, the semi-group method or hypercoercivity (the latter two further rely on the spectral analysis of the linearized Boltzmann operators); see \cites{bardos1991classical,briant2015boltzmann,briant2019boltzmann,gallagher2020convergence,guo2006boltzmann,jiang2018incompressible}.

(2) \underline{Incompressible Euler limits}. To justify the incompressible Euler limits might be the hardest one among the fluid limits from the Boltzmann equations. The key reason is that our understanding of the incompressible Euler equations is far from satisfactory. The only exiting global-in-time solutions to the incompressible Euler equations are the so-called dissipative solutions introduced by Lions \cite{Lions1996}. However, the dissipative solutions are even weaker than the weak solutions in the sense of distributions. So far, there are only some partial results on the incompressible Euler limits global in time, which belong to Saint-Raymond \cites{saint2003convergence, SRM-2009}. She proved that the fluctuations of the renormalized solutions converge to the dissipative solutions of the incompressible Euler equations. There is some recent progress in the analytical data in half-spaces \cite{Jang-Kim-APDE, Cao-Jang-Kim-JDE,kim2024validity}.

(3) \underline{Compressible Euler limits}. The justification of the compressible Euler limits started from late 70's and early 80's by the Japanese school in the framework of analytical solutions \cite{Nishida-CMP1978}, and by Caflisch \cite{caflisch1980fluid}, who first employed the Hilbert expansion method. There are also some progress in this direction, including the simpler acoustic limits (which are the linearization of the compressible Euler equations around constant states), see \cite{guo2010acoustic}. Recently, there are some progress on the Euler limits in domains with boundaries, in particular half-space, \cite{guo2021hilbert, Jiang-Luo-Tang-2024-TAMS}. The main feature of these works is the coupled kinetic-fluid (Knudsen-Prandtl) boundary layers. We remark that these works are all in the linear regimes. In this sense, the linear boundary layers appeared in compressible Euler and acoustic limits are the same, see \cite{Jang-Jiang-DCDS, Guo-Jang-Jiang-KRM, jiang2021compressiblea, Jiang-Wu-2025-KRM}. Currently,  we are far from satisfactory understanding of the nonlinear boundary layers, both in Knudsen and Prandtl layers. 

(4) \underline{Compressible Navier-Stokes-Fourier approximations}. Among the fluid dynamics of the Boltzmann equations as the Knudsen number $\eps\rightarrow 0$, compressible Navier-Stokes-Fourier (briefly, CNS) system is different with other three fluid systems. More specifically, starting from the scaled Boltzmann equation, it can only obtain the CNS depending on the Knudsen number $\eps$. More specifically, the diffusion terms of the CNS are of order $O(\eps)$. In this sense, CNS is not a limiting system of the Boltzmann equation as the Knudsen number $\varepsilon$ tends to zero. It is in fact an {\em asymptotics} derived from the Chapman-Enskog expansion. For more detailed formal derivations, see \cite{bardos1991fluid, SRM-LNM}. More specifically, Bardos-Golse-Levermore derived formally the asymptotics of the CNS from the scaled Boltzmann equation, see \eqref{chap}-\eqref{G}.

There have been some rigorous justifications on the compressible Navier-Stokes approximation, such as \cites{kawashima1979fluid,lachowicz1992solutions,liu2014compressible,duan2021compressible}. In particular, Kawashima, Matsumura, and Nishida \cite{kawashima1979fluid}
demonstrated that the Boltzmann solution converges asymptotically to the CNS solution in the long time for Cauchy problems in the entire space, provided that the initial data for both problems are sufficiently close to constant equilibrium states in smooth Sobolev spaces. In spirit of a previous work \cite{lachowicz1987initial}, Lachowicz \cite{lachowicz1992solutions} initiated a study of the CNS approximation including the initial layer to the Boltzmann equation in spatially periodic domains. In \cite{liu2014compressible}, Liu, Yang, and Zhao justified the CNS approximation in the framework of Chapman-Enskog expansion for the whole space. They assumed the well-prepared initial data, i.e. the difference between the initial data of CNS and Boltzmann is of order $O(\eps^2)$. Later on, in the same spirit, Duan and Liu \cite{duan2021compressible} extended the result of \cite{liu2014compressible} to the bounded domain case for a given homogeneous Dirichlet boundary conditions, $u=0$ and $\theta=0$. In other words, the boundary conditions for the CNS was not derived from the Boltzmann equations, but was directly endowed. 

\subsection{Derivation of the boundary conditions of CNS} 
Our project aims to derive and justify the boundary conditions for the compressible Navier-Stokes-Fourier system from the Boltzmann equations endowed with boundary conditions, Maxwell reflection, or incoming boundary conditions. Previous work in this direction includes Sone's books \cite{sone2002kinetic, sone2007molecular}, and the papers \cites{coron1989derivation,aoki2017slip}. The slip boundary conditions were derived in \cites{coron1989derivation,aoki2017slip}. Specifically, Coron \cite{coron1989derivation} considered the Couette flow between two flat plates endowed with simple Maxwell accommodation conditions. In \cites{aoki2017slip}, the authors considered the Maxwell reflection boundary condition with the accommodation coefficient $\alpha_\eps \in (0,1],\alpha_\eps = O(1)$. One of the key contributions of \cite{aoki2017slip} is that for the Maxwell reflection boundary, the homogeneous Dirichlet boundary conditions $u=0$ and $\theta=0$, are not correct. There are higher-order $O(\eps)$ correctors for the boundary conditions of the CNS. The main technique to determine the boundary conditions for the fluid equations is from the analysis of the Knudsen boundary layer equations. 

Many problems are left to be discussed after \cites{coron1989derivation,aoki2017slip}. First, there have been no rigorous convergence proofs for the CNS approximation with a boundary involving the Knudsen boundary layer. In particular, the justification for the systematic derivation of the boundary conditions for CNS from the Boltzmann equations has not been achieved so far. Second, even at the formal level, from the Boltzmann equations with the Maxwell reflection boundary, only the complete diffusion, that is, the accommodation coefficient $\alpha_\eps=O(1)$ was treated in \cites{coron1989derivation,aoki2017slip}. For the other cases, for example, the almost specular reflection case $\alpha_\eps= O(\eps^\beta)$, $\beta>0$, the corresponding boundary conditions for the CNS are not available in the literature, to the best of our knowledge. Note that when $\beta\rightarrow 0$, it should be consistent with the complete diffusion boundary, while when $\beta\rightarrow \infty$, it should be consistent with the specular reflection boundary. Furthermore, for the Boltzmann equations with other boundary conditions, for example, the incoming boundary, how to derive the boundary conditions for the CNS is completely open. There are only some numerical simulations on this problem; see \cite{sone2007molecular}. Third, even in the cases \cites{coron1989derivation,aoki2017slip} where the CNS boundary conditions were derived, the slip coefficients appearing in the boundary conditions are determined by solving the corresponding Knudsen boundary layer equations. Whether these coefficients make the boundary energy of the CNS dissipative is a highly non-trivial problem. This requires the more quantitative properties of the solutions to the Knudsen layer equations. 

Our project is to address the above problems on the derivations of the boundary conditions of the CNS from the Boltzmann equation. In this paper, which is the first paper in a series, we have two goals, both about the Maxwell reflection boundary of the Boltzmann equation. First, we treat the specular $\alpha_\eps=0$ or the almost specular $\alpha=O(\eps^\beta)$, $\beta>0$ reflection boundary conditions. We systematically derive the slip boundary conditions of the CNS that depend on the Knudsen number. Our analysis follows the procedure for the analysis of the Knudsen layer developed by the Japanese school, as referenced in \cites{sone2002kinetic,sone2007molecular,aoki2017slip}. This approach was initially pioneered by Sone and his team. One of the key difficulties is that, if the accommodation coefficient is $\alpha=O(\eps^\beta)$, and $\beta> 0$ is general, the expansion ansatz was not clear in the previous works. In this paper, a novel appropriate ansatz is used to treat these cases, when $\beta$ is irrational or rational numbers. 

In this paper, it is discovered that when $\alpha_\eps = 0$ or $\alpha_\eps = O(\eps^\beta)$ for $\beta>1$, the boundary conditions of the CNS can be determined directly by the boundary condition of the Boltzmann equation. In these cases, the Knudsen layer equation is not needed. On the other hand, when $0<\beta\leq 1$ the boundary conditions of the CNS variables must be determined by the solvability of the Knudsen layer equation; see Section \ref{Sec-boundary} for more details. For these almost specular reflection boundaries, we focus on the systematic formal analysis in the current paper, and leave the rigorous analysis to the next paper in this series. 

The second goal of this paper is to rigorously justify the CNS approximation from the Boltzmann equation with the specular reflection boundary condition. In this case, the boundary conditions for CNS is the complete slip boundary condition:
 \begin{equation*}
     u\cdot n = 0 ,\quad  ([\nabla_x u + \nabla_x u^{\mathrm{T}})\cdot n]^{\tan} = 0 , \quad  \nabla_x \theta \cdot n =0.
 \end{equation*}
 Here, $n= n (x)$ denotes the normal outward unit vector and $\mathrm{v}^{\tan} = \mathrm{v}  -( \mathrm{v} \cdot n) n$ denotes the tangential part of the vector field $\mathrm{v}$. We make use of the techniques developed in \cites{masmoudi2012uniform, wang2016uniform, duan2021compressible} to obtain uniform estimates in the Knudsen number for the CNS and prove the CNS approximation for the Boltzmann equation in general bounded domains.

The boundary condition to \eqref{be} considered in this paper is the following general moving Maxwell reflection boundary condition:
\begin{equation}
  \begin{aligned}\label{general_Maxwell_reflection}
    \gamma_- F_\eps(t,x,v) = &(1 - \alpha_\eps) F_\eps\left(t,x,v-2[(v-u_w)\cdot n]n\right) \\
     &+ \alpha_\eps \frac{1}{2\pi\theta_w^{2}} e^{-\frac{|v-u_w|^2}{2\theta_w}}\int_{(v-u_w)\cdot n>0}(v-u_w)\cdot n F_\eps(t,x,v')d v'\\
    := & (1 - \alpha_\eps)  L\gamma_+ F_\eps(t,x,v) + \alpha_\eps K\gamma_+ F_\eps(t,x,v).
  \end{aligned}
\end{equation}
 Here, $u_w \in \R^3$ and $\theta_w >0$ represent the velocity and temperature of the boundary, respectively. $ L \gamma_+ $ represents the specular reflection part and $K\gamma_+ $ represents the diffuse reflection part. The quantity $\alpha_\eps \in [0,1]$ denotes the accommodation coefficient, which measures the ratio between specular reflection and diffusive reflection. Throughout this paper, only the specular and almost specular cases are considered, i.e. $\alpha_\eps \ll 1$.
 
 Let $\Sigma:=\partial \Omega \times \mathbb{R}^3$ be the phase-space boundary of $\Omega \times \mathbb{R}^3$. The phase boundary $\Sigma$ can be split by outgoing boundary $\Sigma_+$, incoming boundary $\Sigma_-$, and grazing boundary $\Sigma_0$:
$$\Sigma_+=\{(x,v):x \in \partial\Omega,(v - u_w)\cdot n(x)>0\},$$
$$\Sigma_-=\{(x,v):x \in \partial\Omega,(v - u_w)\cdot n(x)<0\},$$
$$\Sigma_0=\{(x,v):x \in \partial\Omega,(v - u_w)\cdot n(x)=0\}.$$
We denote by $\gamma F$ the trace of $F$ over $\Sigma$ and $\gamma_{\pm}F=\mathbf{1}_{\Sigma_{\pm}}\gamma F$.

\subsection{Chapman-Enskog expansion}
 The Chapman-Enskog expansion is a well-known technique to derive fluid dynamics from kinetic equations. It can be summarized in the following on the derivation of the CNS (without boundary) by the Chapman-Enskog expansion (see \cites{bardos1991fluid,golse2005hydrodynamic} for details). The first-order Chapman-Enskog expansion takes the form
 \begin{equation}\label{chap}
   F_\varepsilon=M+\varepsilon G+O(\varepsilon^2).
 \end{equation}
 Here, the leading-order term $M$ denotes a local Maxwellian distribution, given by
 \begin{equation}
  M= M_{(\rho,u,\theta)}=\frac{\rho}{(2\pi\theta)^{\frac{3}{2}}}e^{-\frac{|v-u|^2}{2\theta}}.
 \end{equation}
 The the first-order correction term $G$ takes the form
 \begin{equation}\label{G}
  G = G_{(\rho,u,\theta)}=-M \big(\frac{1}{2}\hat{A}( V ):\sigma(u)+\hat{B}( V )\cdot \frac{\nabla_x\theta}{\sqrt \theta}\big),
 \end{equation}
 where $ V $ denotes the scaled velocity variable
 \begin{equation}\label{xi}
   \qquad  V =\frac{v-u}{\sqrt \theta},
 \end{equation}
 and $\sigma(u)$ is the traceless part of the deformation tensor of $u$, given by
 \begin{equation}
   \sigma(u)=\nabla_x u+(\nabla_x u)^{\mathrm{T}}-\frac{2}{3}\nabla_x\cdot u I.
 \end{equation}
 The matrix-valued function $\hat{A}( V )$ and vector-valued function $\hat{B}( V )$ are defined as
 $$\hat{A}( V )=\frac{1}{\rho \sqrt{\theta}} a(\theta, | V |)A( V ), \quad \hat{B}( V )=\frac{1}{\rho \sqrt{\theta}} b(\theta, | V |)B( V ),$$
 with 
 $$A( V )= V \otimes  V -\frac{| V |^2}{3}I,\quad B( V )= \frac{(| V |^2-5) V }{2} ,$$
 and $a(\theta, | V |)$ and $b(\theta, | V |)$ are two scalar functions depends only on $\theta$ and $| V |$. Moreover, it holds that
 \begin{equation}\label{positivity}
  \int_{\R^3 } \hat{A}_{ij}(v) A_{ij}(v) M(v) d v  >0, \quad \int_{\R^3 } \hat{B}_i(v) B_i(v) M(v) d v  >0.
 \end{equation}
 
 The macroscopic variables $(\rho,u,\theta)$ in the leading order satisfy the full compressible Navier-Stokes system with $O(\varepsilon)$ dissipation terms:
 \begin{equation}\label{cns}
 \begin{cases}
 \partial_t \rho+\nabla_x \cdot (\rho u)=0, \\
 \rho(\partial_tu+u\cdot \nabla_xu)+\nabla_x(\rho\theta)=\varepsilon\nabla_x\cdot[\mu(\theta)\sigma(u)],\\
 \frac{3}{2}\rho(\partial_t \theta+u \cdot \nabla_x\theta)+(\rho\theta)\nabla_x \cdot u=\varepsilon\nabla_x\cdot[\kappa(\theta)\nabla_x \theta]+\frac{\varepsilon}{2}\mu(\theta)\sigma(u):\sigma(u),\\
 \end{cases}
 \end{equation}
 where the viscosity coefficient $\mu(\theta)$ and the heat conductivity coefficient $\kappa(\theta)$ are given by
 \begin{equation}
 \begin{aligned}\label{def_mu_kappa}
 &\mu(\theta)= \frac{\theta}{10} \int_{\R^3 } A(V): \hat{A}(V) M dv =   \frac{2}{15}\sqrt{\theta}\int_{0}^{\infty} a(\theta,r)r^6e^{-\frac{r^2}{2}}\frac{dr}{\sqrt{2\pi}}>0,\\
 &\kappa(\theta)= \frac{\theta}{3} \int_{\R^3 } B(V) \cdot \hat{B}(V) M dv = \frac{1}{6}\sqrt{\theta}\int_{0}^{\infty} b(\theta, r)r^2(r^2-5)^2e^{-\frac{r^2}{2}}\frac{dr}{\sqrt{2\pi}}>0.
 \end{aligned}
 \end{equation}
 In particular, $a $ and $b $ are functions independent of $\theta$ for hard-sphere molecules.

 \subsection{Statement of the main results} 
 The main results of this paper include two parts. The first is on the formal derivation of the boundary conditions for the CNS from the Boltzmann equation with the Maxwell reflection boundary condition when the accommodation coefficient $\alpha_\eps=0$ or $\alpha_\eps=\chi \eps^\beta$. The second is on the rigorous justification of the CNS for the specular reflection case. 
 
 \subsubsection{Derivation of slip boundary conditions} When the Boltzmann equation is endowed with the specular or almost specular reflection boundary conditions, various slip boundary conditions for the CNS can be derived. The boundary conditions for the CNS can be summarized as follows. Recall that $n = n(x)$ is the normal unit vector outward.
\begin{theorem}
      Let $\Omega$ be any rigid body with smooth boundary $\partial\Omega$. Let the Boltzmann equation \eqref{be} be imposed with the general Maxwell reflection boundary condition \eqref{general_Maxwell_reflection}.
      \begin{itemize}
          \item  When the accommodation coefficient $\alpha_\eps$ satisfies $\alpha_\eps=0 $ or $\alpha_\eps = \chi \eps^\beta  $ for $\beta> 1$, where the positive constant $\chi = O(1)$, the boundary conditions for the CNS are the complete slip boundary conditions
    \begin{equation}\label{slip_theorem}
        \left\{
            \begin{aligned}
                & (u-u_w)\cdot n=0\,, \\
                & [(\nabla_x u+\nabla_x u^\mathrm{T})\cdot n]^{\tan}=0 \,, \\
                & \nabla_x \theta \cdot n=0 \,.
            \end{aligned}
        \right.
    \end{equation} 
    
     \item When the accommodation coefficient $\alpha_\eps=\chi \eps^\beta ,$ for $0<\beta < 1$, the boundary conditions for the CNS are
    \begin{equation}\label{Navier-slip1-theorem}
        \left\{
            \begin{aligned}
                & (u-u_w)\cdot n=0 ,\\
                & [u-u_w]^{\tan}=\frac{\eps^{1-\beta} b^I_u}{ \chi \rho}  [(\nabla_x u+\nabla_x u^\mathrm{T}) \cdot n]^{\tan} , \\
                & \theta-\theta_w=\frac{\eps^{1-\beta} b^I_\theta}{\chi \rho}  \nabla_x \theta \cdot n . \\
            \end{aligned}
        \right.
    \end{equation}
     Here, $b^I_u < 0$ and $b^I_\theta < 0$ denote the so-called slip coefficients, given by \eqref{slip-coefficients1}.

      \item When the accommodation coefficient $\alpha_\eps=\chi \eps$, the boundary conditions for the CNS are
    \begin{equation}\label{Navier-slip2-theorem}
        \left\{
            \begin{aligned}
                & (u-u_w)\cdot n=0 ,\\
                & [u-u_w]^{\tan}=\frac{ c^I_u}{ \chi \rho}  [(\nabla_x u+\nabla_x u^\mathrm{T}) \cdot n]^{\tan} , \\
                & \theta-\theta_w=\frac{ c^I_\theta}{\chi \rho}  \nabla_x \theta \cdot n + \frac{|u - u_w|^2}{4}. \\
            \end{aligned}
        \right.
    \end{equation}
     Here, the slip coefficients $c^I_u < 0$ and $c^I_\theta < 0$ are given by \eqref{slip-coefficients2}. 
      \end{itemize}
\end{theorem}
     
     \begin{remark}    When the accommodation coefficient $ \alpha_\eps=O(1),\,$  the boundary conditions for the CNS were derived by Aoki et al. in \cite{aoki2017slip}
    \begin{equation}\label{Aoki_slip}
        \left\{
            \begin{aligned}
                & (u-u_w)\cdot n=0 \,,\\
                & [u-u_w]^{\tan}=\frac{\eps a^{I}_u }{\rho } [(\nabla_x u+\nabla_x u ^\mathrm{T})\cdot n]^{\tan} + \frac{\eps a^{I}_\theta }{\rho \sqrt{\theta_w}} [\nabla_x\theta]^{\tan}\,, \\
                & \theta-\theta_w=\frac{\eps a^{II}_u \sqrt{\theta_w}}{\rho} \partial_i u_j n_i n_j+ \frac{\eps a^{II}_\theta }{\rho }\nabla_x \theta \cdot n \,. \\
            \end{aligned}
        \right.
    \end{equation}
    The coefficients $a^I_u$, $a^I_\theta, a^{II}_u $ and $a^{II}_\theta$  are the slip coefficients. For hard-sphere molecules and $\alpha_\eps = 1$, their specific values are summarized in \cite[Subsection 5.3]{aoki2017slip}.
     \end{remark} 
\begin{remark}
 The above boundary conditions are consistent when the values of $\beta$ change continuously in the following sense: 
    \begin{itemize}
        \item when $\beta\rightarrow +\infty$ (this is equivalent to $\alpha_\eps=0$), the boundary condition \eqref{Navier-slip1-theorem} tends to \eqref{slip_theorem}.
        \item when $\beta\rightarrow 1$ (this is equivalent to $\alpha_\eps=\chi \eps$), the boundary condition \eqref{Navier-slip1-theorem} tends to \eqref{Navier-slip2-theorem}.
        \item when $\beta\rightarrow 0$ (this is equivalent to $\alpha_\eps=O(1)$), the boundary condition \eqref{Navier-slip1-theorem} tends to \eqref{Aoki_slip}.
    \end{itemize}
  
\end{remark}
 \begin{remark}
        In the scenarios where $\alpha_\eps = 0$ and $\alpha_\eps = O(\eps^\beta)$ for $ \beta>1$, the Knudsen layer is not needed. The boundary conditions can be directly determined from the Boltzmann equation. In all other cases, the Knudsen layer equation must be used to determine the boundary conditions of the CNS.
 \end{remark}
 
 \subsubsection{Rigorous justification for specular reflection} The second result is a rigorous justification for the approximation from the Boltzmann equation with the specular reflection boundary condition to the CNS with the complete slip boundary condition. For technical simplicity, only the hard-sphere model is considered in this paper, i.e.
$$b(|v_1-v|,\omega )=|(v_1-v)\cdot\omega|.$$ 
It could be noted that for the more general collision kernels, there have been mature techniques to treat them. This paper focuses on the fluid limit, so the technical complexity on the collision kernels will not be emphasized.

The Boltzmann equation is imposed with the specular reflection with $u_w=0$, that is, 
\begin{equation}\label{spec}
\gamma_-F_\varepsilon(t,x,v)=F_\eps (t,x, v-2(v\cdot n)n)
 \quad \textrm{ on} \quad \! \mathbb{R}^+\times \Sigma_-.
\end{equation}

The domain $\Omega$ is assumed to be any smooth bounded domain. It will be considered for the solution to the CNS close to $(1,0,1)$, so the following ansatz for the distribution function is introduced,
 \begin{equation}\label{chap2}
   F_\varepsilon=M+\varepsilon G+\varepsilon^2\sqrt{\mu}R,
 \end{equation}
 where $\mu  $ denotes the global Maxwellian distribution defined in \eqref{mu}, $R$ represents the second order correction term, i.e., the remainder term. We now introduce the appropriate functional space for our analysis
   \begin{equation}\label{Rnorm}
     \begin{aligned}
         \|R(t)\|_{\mathbf{X}}^2
         =&\sup_{\tau\in[0,t]}\sum_{\alpha_0\leq 1}\left\|\partial^{\alpha_0}_tR\right\|^2_{L^2_{x,v}}
         +\varepsilon\sup_{\tau\in[0,t]}\left\|\mathbf{P}R\right\|^2_{L^6_{x,v}}
         +\varepsilon^2 \sup_{\tau\in[0,t]}\left\|w_kR\right\|^2_{L^\infty_{x,v}}\\
         &+\varepsilon\sum_{\alpha_0\leq 1}\int_{0}^{t}\left\|\mathbf{P}\partial^{\alpha_0}_tR\right\|_{L^2_{x,v}}^2d\tau
         +\frac{1}{\varepsilon}\sum_{\alpha_0\leq 1}\int_{0}^{t}\left\|(\mathbf{I-P})\partial^{\alpha_0}_tR\right\|_{L^2_{x,v}}^2d\tau,
     \end{aligned}
     \end{equation}
 where
 \begin{equation}\label{wk}
   w_k(v)= {(|v|^2+1)}^{\frac{k}{2}},
 \end{equation}
 the linearized Boltzmann operator $\mathcal{L}$ around $\mu$ is given by \eqref{linearbol}, and $\mathbf{P}R$ denotes the projection of $R$ onto the kernel of $\mathcal{L}$, as described in Subsection \ref{Subse:notation}. We will derive the uniform estimates for $R$ in the space $\|\cdot\|_\mathbf{X}$, see Section \ref{Se:uniform}.
 
 To solve \eqref{cns} with boundary conditions \eqref{slip}, we impose the initial data
 \begin{equation}\label{cnsinitial}
   (\rho,u,\theta)(0,x)=(\rho_0,u_0,\theta_0)(x).
 \end{equation}
 Moreover, we impose the well-prepared initial data on the scaled Boltzmann equation \eqref{be}
 \begin{equation}\label{beinitial}
   F_\varepsilon(0,x,v)=M^{in}+\varepsilon G^{in}+\varepsilon^2\sqrt{\mu}R^{in}.
 \end{equation}
 Now, we state our main results on the CNS approximation from the Boltzmann equation.
 \begin{theorem}\label{Tm:main}
  Let $w_k(v)$ be defined in \eqref{wk} with $k\geq \frac{7}{2}$. Let $T>0$ be any prescribed constant. There is a small constant $\varepsilon_0>0$ such that if for $\varepsilon\in(0,\varepsilon_0)$
  \begin{equation}
   \|( {\rho}_0,u_0, {\theta}_0)\|_{\mathbb{X}}\leq \lambda_0\varepsilon^{3/2},
  \end{equation}
 where the norm $\|\cdot\|_{\mathbb{X}}$ is defined in \eqref{cnsnorm}. Then the scaled Boltzmann equation \eqref{be} with the specular reflection boundary condition \eqref{spec} and well-prepared initial data \eqref{beinitial} admits a unique solution for $\varepsilon\in(0,\varepsilon_0)$ over the time interval $t\in[0,T]$ with expanded form \eqref{chap2}, that is,
   \begin{equation}
     \begin{aligned}
         F_\varepsilon=M+\varepsilon G+\varepsilon^2\sqrt{\mu}R\geq 0,
     \end{aligned}
     \end{equation}
 where $M=M_{(\rho, u,\theta)}$, and $(\rho,u,\theta)$ represents the unique solution described in Theorem \ref{Thmcns}. The remainder term $R$ is constructed in Lemma \ref{Lm:reuniform}, satisfying
 \begin{equation}\label{main}
  \|R\|^2_{\mathbf{X}}\lesssim \delta,
 \end{equation}
 for $\delta$ suitably small. Moreover, it holds that
 \begin{equation}\label{main1}
   \sup_{\tau\in[0,T]}\left\|\frac{F_\eps-M}{\sqrt{\mu}}\right\|_{L^2_{x,v}}\lesssim \varepsilon^2, \qquad \sup_{\tau\in [0,T]}\left\|w_k(v)\frac{F_\eps-M}{\sqrt{\mu}}\right\|_{L^\infty_{x,v}}\lesssim \varepsilon.
 \end{equation}
  \end{theorem}
   Therefore, we justified the compressible Navier-Stokes approximation for the scaled Boltzmann equation with the specular boundary condition in bounded domains.

\subsection{Notations}\label{Subse:notation}

We employ the following notation for function spaces and norms. Let $H^s$ denote the standard Sobolev space $W^{s,2}(\Omega)$ equipped with the norm $\|\cdot\|_{H^s}$. We denote by $|\cdot|_{H^s}$ the standard $H^s(\partial \Omega)$ norm. For $1\leq p \leq \infty$, $\|\cdot \|_{L^p_{x,v}}$ and $\|\cdot \|_{L^p_x}$ denote the standard norms on $L^p(\Omega \times \R^3)$ and $L^p(\Omega )$, respectively. When unambiguous, we use $\|\cdot\|_p$ to denote $\|\cdot \|_{L^p_x}$. 
 Moreover, $(\cdot,\cdot)$ denotes the $L^2$-inner product over $x\in\Omega$ and $\left<\cdot,\cdot\right>$ denotes the $L^2$-inner product over $(x,v)\in\Omega\times\mathbb{R}^3$. The notation $A\lesssim B$ indicates that  $A\leq CB$ for some harmless constant $C>0$.

We collect some well-known results on the Boltzmann collision operator, referring the reader to the foundational work \cites{caflisch1980fluid,Cercignani-Illner-Pulvirenti} for further details. We introduce the linearized Boltzmann operator $\mathcal{L}$ defined as
\begin{equation}\label{linearbol}
   \mathcal{L} f(v)=-\frac{1}{\sqrt{\mu}}\{Q(\mu,\sqrt{\mu }f)+Q(\sqrt{\mu}f,f)\},
  \end{equation}
  where $\mu$ denotes the global Maxwellian distribution given by
\begin{equation}\label{mu}
   \mu(v)=\frac{1}{(2\pi)^\frac{3}{2}}e^{-\frac{|v|^2}{2}}.
\end{equation}
The null space $\operatorname*{ker}(\mathcal{L})$ of the operator $\mathcal{L}$ is generated by the collision invariants
\begin{equation}\label{1.3.4}
    \sqrt{\mu}, \quad
    v_i \sqrt{\mu} \, (i=1,2,3), \quad 
   \frac{|v|^2-3}{2}\sqrt{\mu}. 
\end{equation}
For cutoff hard sphere molecules, $ \mathcal{L}$ can be decomposed as $$ \mathcal{L}=\nu(v)-K,$$
where $K$ is a self-adjoint operator defined as
\begin{equation}\label{K}
Kf=\int_{\mathbb{R}^3}k(v,v')f(v')dv',
\end{equation}
and $\nu(v)$ denotes the collision frequency. It holds that
\begin{equation*}
  \nu(v)=\int_{\mathbb{R}^3}\int_{\mathbb{S}^2}|(v_1-v)\cdot \omega|\mu(v_1)d\omega dv_1,
\end{equation*}
and
\begin{equation}\label{1.3.5}
  1+|v|  \lesssim \nu(v)  \lesssim 1+|v|.
\end{equation}
Furthermore, we introduce the $L^2_v$ projection $\mathbf{P}$ associated with $\mathcal{L}$. Then, we can split any admissible function $f$ into 
\begin{equation}
  f = \mathbf{P} f + (\mathbf{I-P})f.
\end{equation}
The coercivity of operator $\mathcal{L}$ implies
\begin{equation}\label{hyperco}
  \left<\mathcal{L}g,g\right> \gtrsim \|(\mathbf{I-P})g\|_\nu^2,
\end{equation}
with the weight $L^2$-norm $\|g\|_\nu^2:=\int_{\Omega\times\mathbb{R}^3} g^2(x,v)\nu dxdv$. We further introduce the bilinear term
\begin{equation}\label{def_gamma}
  \Gamma(f,g)=\frac{1}{\sqrt{\mu}}Q(\sqrt{\mu}f,\sqrt{\mu}g),
\end{equation}
which belongs to $\operatorname{ker}^\perp(\mathcal{L})$.

\subsection{Difficulties and Methodologies} In this series of work, we aim to connect the boundary conditions for the compressible Navier-Stokes system from the Boltzmann equations with the Maxwell reflection boundary.

\underline{\bf Part 1. Derivation of the slip boundary conditions}.  In the work of Aoki et al. \cite{aoki2017slip}, the Maxwell reflection boundary condition was studied under the assumption that the accommodation coefficient satisfies $\alpha_\varepsilon = O(1)$. In this paper, inspired by the idea of \cite{aoki2017slip}, and some new idea about the Chapman-Enskog expansion, the left cases are treated. Thus, the derivation of the boundary conditions of the CNS from the Maxwell reflection condition for the Boltzmann equation is complete. 

In the following, we illustrate the idea on how the structure of the asymptotic ansatz depends on the value of $\alpha_\varepsilon$. We emphasize that this ansatz differs significantly from the Hilbert expansion considered in \cite{jiang2021compressiblea}, where the expansion takes the form
\begin{align}\label{eqn: short-F_eps}
F_\varepsilon = \sum_{m,n\geq 0} \varepsilon^{\frac{1}{2}n + \beta m} F_{mn}.
\end{align}
This highlights a key difference between the Hilbert expansion and the Chapman–Enskog expansion.

First, consider the cases where $\alpha_\varepsilon = 0$ or $\alpha_\varepsilon = O(\varepsilon^\beta)$ for $\beta > 1$. In these regimes, the diffusive component is sufficiently weak that it does not contribute to the first-order Chapman–Enskog expansion. We find that, in such cases, the ansatz \eqref{chap} satisfies the Maxwell reflection boundary condition \eqref{general_Maxwell_reflection} up to $O(\varepsilon)$ if and only if the fluid variables $(\rho, u, \theta)$ satisfy the complete slip boundary conditions \eqref{slip_theorem}; see Subsection~\ref{subsec_specular} for more details. Therefore, the Knudsen layer is not needed at the first order. In this sense, the Knudsen layer is too weak to be captured by the first-order Chapman–Enskog expansion.

For the intermediate regime where $\alpha_\eps = O(\eps^\beta)$ for $0 <\beta <1$, the construction of the asymptotic ansatz becomes significantly more intricate. In this case, the diffusive effects are no longer negligible. To directly satisfy the Maxwell reflection boundary condition \eqref{general_Maxwell_reflection}, one would need to impose
$$u - u_w =0, \quad \theta - \theta_w = 0,$$
as well as
\begin{equation*}
    [(\nabla_x u + \nabla_x u^\mathrm{T})\cdot n]^{\tan} = 0 , \quad \nabla_x \theta \cdot n =0.
\end{equation*}
These constraints are overly restrictive as boundary conditions for the CNS, implying the necessity of the Knudsen layer to resolve this mismatch. To ensure the consistency up to $O(\eps^\beta)$, we assume that 
\begin{equation}\label{assumption}
    u -u_w =O(\eps^\iota), \quad \theta - \theta_w = O(\eps^\iota)
\end{equation}
for some exponent $\iota >0$ to be determined. Our subsequent analysis shows that this assumption is indeed consistent and yields $\iota = 1 - \beta$. Under this assumption, the Maxwell reflection boundary condition \eqref{general_Maxwell_reflection} can be satisfied up to $O(\eps^\beta)$. To match the boundary at $O(\eps)$, we introduce a Knudsen layer of thickness $O(\varepsilon)$ adjacent to the boundary\footnote{For incoming boundary condition for the Boltzmann equation, the Knudsen layer may appear at the leading order, see \cite{1991-PFA-Aoki}}. Consequently, our ansatz takes the form 
\begin{equation}
\begin{aligned}
    F_\eps = M &+\eps G +\cdots \\
     &+ \eps F^{bb}
\end{aligned}
\end{equation}
 where $F^{bb}$ denotes the Knudsen layer correction. Combining assumption \eqref{assumption} with the solvability condition of the Knudsen layer equation, we derive the slip boundary conditions with the slip part of $O(\eps^{1-\beta})$, namely, \eqref{Navier-slip1-theorem}.

 Finally, consider the critical case $\alpha_\eps = O(\eps)$. The analyze in this regime is similar to that for $\alpha_\eps = O(\eps^\beta)$ for $0<\beta<1$. We point out that the relation \eqref{assumption} now implies $u - u_w =O(1), \theta - \theta_w  =O(1)$. As a result, the addition term $|u - u_w|^2 $ appears in \eqref{Navier-slip2-theorem}. This term, regarded as a higher-order correction in the case $0 <\beta < 1$, was previously neglected.

\underline{\bf Part 2. Rigorous justification}
In \cite{duan2021compressible}, the diffuse reflection boundary condition was investigated for the Boltzmann equation \eqref{be}.  In their work, the boundary conditions for the CNS were given directly as $(u,\theta)=(0,0)$ in $\partial\Omega$, not by systematic derivation. However, in this case, to derive the corresponding boundary conditions of the CNS, boundary layer equations must be employed. The complete formal derivation was recently provided by Aoki et al. \cite{aoki2017slip}. In this sense, the boundary layer effects were neglected in \cite{duan2021compressible}.  For the CNS approximation from the Boltzmann equation endowed Maxwell reflection condition with accommodation coefficient $\alpha_\eps\in [0,1]$, the corresponding boundary conditions for the CNS must be derived from the Boltzmann equation. In the present paper, we rigorously justify the CNS approximation from Boltzmann equation with specular reflection boundary. As derived in the first part of this paper, the boundary conditions for the CNS are the complete slip boundary conditions, derived by our formal analysis. Furthermore, we justify the expansion with optimal scaling 
$$F_\eps(t,x,v)= M + \eps G +\eps^2 \sqrt{\mu} R.$$
Our proof mainly has two parts: uniform estimates for the CNS and uniform estimates for the remainder equation.

\underline{(1). Uniform estimates for the CNS solution}. For the CNS with complete slip boundary conditions, the techniques in \cites{masmoudi2012uniform,wang2016uniform,duan2021compressible} will be employed to derive uniform estimates. As in \cite{duan2021compressible}, we introduce the energy norm $\mathbb{X}$ with conormal derivatives. In contrast with the Dirichlet boundary, some new difficulties arise. To overcome these, we employ the Helmholtz decomposition and elliptic estimates from Agmon-Douglis-Nirenberg \cite{agmon1964estimates}. We aim to derive the following uniform estimates by a bootstrap argument
\begin{equation*}
  \begin{aligned}
    \|(\rho &, u, \theta)(t) \|_{\mathbb{X}}^{2} \\
    =&\sup _{\tau\in[0,t]} \sum_{\alpha_{0} \leq 3}\left\|\partial_{t}^{\alpha_{0}}(\rho-1, u, \theta-1)(\tau)\right\|_{2}^{2}
    +\sup _{\tau\in[0,t]} \sum_{\alpha_{0} \leq 2}\left\|\partial_{t}^{\alpha_{0}} \nabla_{x}(\rho, u, \theta) (\tau)\right\|_{2}^{2} \\
    & +\sup _{\tau\in[0,t]} \sum_{\alpha_{0} \leq 2} \varepsilon^{2}\left\|\partial_{t}^{\alpha_{0}} \nabla_{x}^{2}(\rho,u, \theta)(\tau)\right\|_{2}^{2} +\sup _{\tau\in[0,t]} \sum_{\alpha_{0} \leq 2} \varepsilon\left\|\partial_{t}^{\alpha_{0}}({\rho}, u, {\theta}) (\tau)\right\|_{H^2_{co}}^{2}\\
    &+\sup _{\tau\in[0,t]} \sum_{\alpha_{0} \leq 1} \varepsilon^{2}\left\|\partial_{t}^{\alpha_{0}} \nabla_{x}(\rho, u, \theta) (\tau)\right\|_{H_{c o}^{2}}^{2}
    +\sup _{\tau\in[0,t]} \sum_{\alpha_{0} \leq 1} \varepsilon^{4}\left\|\partial_{t}^{\alpha_{0}} \nabla_{x}^{2} \rho(\tau)\right\|_{H_{c o}^{1}}^{2}\\
    & +\sup _{\tau\in[0,t]}  \varepsilon^{4}\left\| \nabla_{x}^{2} \rho(\tau)\right\|_{H_{c o}^{2}}^{2}
    +\sup _{\tau\in[0,t]} \sum_{\alpha_{0} \leq 1} \varepsilon^{4}\left\|\partial_{t}^{\alpha_{0}} \nabla_{x}^{2}(u, \theta)(\tau)\right\|_{H_{c o}^{2}}^{2} \lesssim  \|(\rho, u, \theta)(0) \|_{\mathbb{X}}^{2}.
    \end{aligned}
\end{equation*}
Here, $\|\cdot\|_{H_{co}^1}$ and $\|\cdot\|_{H_{co}^2}$ denote the conormal Sobolev norms, whose precise definitions are provided in Section \ref{Sec_CNS_uniform}. To establish the desired uniform estimates, we first derive the zero-order energy estimate \eqref{L2}, which notably incorporates pure time derivatives up to third order. Higher-order energy estimates for $(\rho,u,\theta)$ require a more refined analysis, achieved through a combination of Helmholtz decomposition and elliptic regularity techniques (see Lemmas \ref{Lm:H1}, \ref{Lm:rhoH2}, and \ref{Lm:rhouthetaH3} for details). In the last step, we perform the same energy estimates by first acting the conormal differentiation $Z^\alpha(\alpha\leq2)$ and conclude the proof.

\underline{(2). Uniform estimates of the remainder $R$.} We employ the standard $L^2\mbox{-}L^6\mbox{-}L^\infty$ framework (see \cites{esposito2013non,Esposito-Guo-Kim-Marra-2018}, for instance) to obtain uniform estimates of the remainder term $R$. We also utlize the new $L^6$ estimates developed in \cite{chen2023macroscopic}, see Section \ref{Se:uniform} for details.

Since the specular reflection boundary condition is considered, the corresponding $L^3$ estimates for the macroscopic part, as discussed in \cites{Esposito-Guo-Kim-Marra-2018,duan2021compressible} based on the velocity average lemma, are not available. It is worth explaining the difference with \cite{duan2021compressible}. In their work, they considered the complete diffuse boundary $\alpha_\eps=1$ and employed the expansion:
\begin{equation*}
  F_\eps=M+\eps G+\eps^{3/2}R.
\end{equation*}
Particularly, the $3/2$ order of $\eps$ mainly comes from the singularity of the boundary. They choose $\eps^{3/2}$ instead of $\eps^2$ to reduce the singularity at the boundary. For the diffuse boundary condition, the Knudsen layer must appear, and the most precise boundary conditions of the CNS are typically $\eps\mbox{-}$related slip boundary conditions, as formally illustrated in \cite{aoki2017slip}. In \cite{duan2021compressible}, the simple homogeneous Dirichlet boundary condition $(u,\theta)=(0,0)$ for the CNS was given, neglecting the boundary layer, thereby leading to the singularity on the boundary of the remainder.

In contrast, we consider the specular reflection boundary condition and derive the complete slip boundary condition for the CNS. Based on our derivation, the specular reflection boundary condition \eqref{spec} is fulfilled up to $O(\eps)$. Therefore, this kind of singularity in \cite{duan2021compressible} actually does not arise. Consequently, we take the expansion $F=M+\eps G+\eps^2 \sqrt{\mu} R$. It turns out that the $L^3$ estimates are unnecessary, allowing us to conclude our proof. In our analysis, we choose the global Maxwellian \eqref{mu} as our reference equilibrium state, rather than the local Maxwellian $M_{(\rho,u,\theta)}$. This choice yields significant simplifications in the subsequent calculations. The validity of this approach stems from our solution regime. The solutions to the compressible Navier-Stokes system under consideration remain sufficiently close to the constant equilibrium state $(1,0,1)$. 

More specifically, the remainder equation reads
\begin{equation*}
  \partial_tR+v\cdot\nabla_xR+\frac{1}{\varepsilon}\mathcal{L}R=\varepsilon\Gamma(R,R) + l.o.t.,
\end{equation*}
with the specular boundary condition $\gamma_-R=L\gamma_+R$. Here and below $l.o.t.$ refers to lower-order term and can be treated easily. The linearized Boltzmann collision operator $\mathcal{L}$ is defined in \eqref{linearbol}, and the bilinear term $\Gamma$ is introduced in \eqref{def_gamma}.

In establishing the $L^2$ estimates, the most key point is to treat the trilinear term $ \varepsilon\int_{0}^{t} \int_{\Omega \times \mathbb{R}^3} \Gamma(R,R)Rdxdvd\tau$. We note that $\Gamma(R,R) \in \operatorname{ker}^\perp(\mathcal{L})$. Decompose $R$ into macroscopic part and microscopic part $R=\mathbf{P}R +\mathbf{(I-P)}R$. Then H{\"o}lder's inequality gives
\begin{equation}
  \varepsilon\int_{0}^{t} \int_{\Omega \times \mathbb{R}^3} \Gamma(R,R)Rdxdvd\tau\lesssim \varepsilon^{-1}\int_{0}^{t}\left\|\nu^{\frac{1}{2}}(\mathbf{I-P})R\right\|^2_{L^2_{x,v}}d\tau
  +\varepsilon^3\left\|w_kR\right\|^2_{L^\infty_{x,v}}\int_{0}^{t}\|\mathbf{P} R\|_{L^2_{x,v}}^2d\tau+l.o.t.,
\end{equation}
where $\left\|w_kR\right\|_{L^\infty_{x,v}}$ can be estimated as
\begin{equation*}
  \varepsilon\left\|w_kR\right\|_{L^\infty_{x,v}}\lesssim \varepsilon^{1/2}\left\|\mathbf{P}R\right\|_{L^6_{x,v}}+\varepsilon^{-1/2}\left\|(\mathbf{I-P})R\right\|_{L^2_{x,v}}+l.o.t..
\end{equation*}
The second term on the right-hand side can be estimated as
\begin{multline*}
  \varepsilon^{-1/2}\left\|(\mathbf{I-P})R(\tau)\right\|_{L^2_{x,v}}\lesssim \varepsilon^{-1/2}\left\|(\mathbf{I-P})R(0)\right\|_{L^2_{x,v}}\\
  +2\varepsilon^{-1}\int_{0}^{t}\left\|(\mathbf{I-P})R(\tau)\right\|_{L^2_{x,v}}\left\|(\mathbf{I-P})\partial_tR(\tau)\right\|_{L^2_{x,v}}d\tau.
\end{multline*}
The $L^2$ estimate in Lemma \ref{Lm:reL2} indicates that
\begin{equation*}
  \int_{0}^{t}\left\|(\mathbf{I-P})R\right\|^2_{L^2_{x,v}}d\tau \sim O(\varepsilon).
\end{equation*}
The pivotal macroscopic estimate, as stated in Lemma \ref {Lm:maL2L6}, leads to
\begin{equation*}
  \int_{0}^{t}\left\|\mathbf{P}R\right\|^2_{L^2_{x,v}}d\tau \sim O(\varepsilon^{-1}),
  \quad \text{ and }\quad  \left\|\mathbf{P}R\right\|^2_{L^6_{x,v}} \sim O(\varepsilon^{-1}).
\end{equation*}
We can finally conclude our proof by acting $\partial_t$ to the remainder equation, see Section \ref{Se:uniform} for details.




\subsection{Arrangement of the paper}

The paper is arranged as follows. In Section 2, we derive the boundary conditions for the CNS from the Boltzmann equation with the Maxwell reflection boundary condition.  Section 3 is devoted to establishing the uniform regularity for the full compressible Navier-Stokes system.  In Section 4, we establish uniform estimates for the remainder term in the framework $L^2\mbox{-}L^6\mbox{-}L^\infty$. In Appendix \ref{Se_useful_estimate}, we introduce some useful estimates. In Appendix \ref{App.A}, we verify the Complementing Condition of Agmon-Douglis-Nirenberg \cite{agmon1964estimates} for an elliptic system, which is essential to establish the uniform regularity of the compressible Navier-Stokes system.

\section{Boundary conditions for the Compressible Navier-Stokes system}\label{Sec-boundary}
In this section, the boundary conditions for the compressible Navier-Stokes system derived from the Boltzmann equation will be derived.  The general Maxwell reflection boundary condition with a moving boundary \eqref{general_Maxwell_reflection} with $\alpha_\eps \ll 1$ will be considered. It is worthy to note that the boundary condition \eqref{general_Maxwell_reflection} implies that there is no instantaneous mass flow across the boundary, i.e.,
\begin{equation}\label{flow_across_boundary}
  \int_{\R^3} (v - u_w ) \cdot n F_\eps (t,x,v) d v= 0.
\end{equation}
Throughout this section, the domain $\Omega$ can represent any rigid body, such as bounded domains, half-spaces, exterior domains, etc. 
  
 \subsection{$\alpha_\eps=0$ or $\alpha_\eps= \chi \eps^\beta$ with $\beta>1$ }\label{subsec_specular} In this subsection, the accommodation coefficient is settled as $\alpha_\eps= 0$, or $\chi \eps^\beta$ with $\beta>1$, i.e., specular reflection or extremely nearly specular reflection. To be consistent with the Maxwell reflection boundary condition \eqref{general_Maxwell_reflection} up to $O(\eps)$ considered in the Chapman-Enskog expansion \eqref{chap}, it is needed to satisfy the boundary conditions up to $O(\eps)$.  By substituting \eqref{chap} into the boundary condition \eqref{general_Maxwell_reflection},  and rearranged by the order of $\eps$, we obtain
 \begin{equation*}
     \gamma_- M = L \gamma_+ M.
 \end{equation*}
 That is
 \begin{equation*}
  \frac{\rho}{({2\pi\theta})^\frac{3}{2}}e^{-\frac{|v - u |^2}{2 \theta}} =  \frac{\rho}{({2\pi\theta})^\frac{3}{2}}e^{-\frac{|v - 2[(v - u_w) \cdot n]n - u |^2}{2 \theta}}
 \end{equation*}
for $(v - u_w) \cdot n <0$. Combining with the fact $|v - u | = | v - u - 2 [(v - u) \cdot n] n |$, it infers that
 \begin{equation}\label{u-u_w-n}
   (u - u_w )\cdot n=0,\quad \textrm{ on } \partial\Omega.
 \end{equation}
In order that the first-order term $G$ defined in \eqref{G} satisfies \eqref{general_Maxwell_reflection} at $O(\eps)$, we then have
 \begin{equation*}
  (\gamma_--L\gamma_+)G=0.
 \end{equation*}
Define $V = \frac{v - u}{\sqrt{\theta}}$. Hence, it holds by using \eqref{u-u_w-n}, $L\gamma_+ f(V) = f(\frac{v -2[(v-u_w)\cdot n]n - u}{\sqrt{\theta}}  ) = f(R_x V)$, where $R_x V = V - 2(V \cdot n)n$. We thereby get
 \begin{equation}\label{G-LgammaG}
   \frac{\rho}{{(2\pi\theta)}^\frac{3}{2}}e^{-\frac{|V|^2}{2}}\big(\frac{a(\theta,|V|)}{2 \rho\sqrt{\theta}} {[A(V)-A(R_xV)]:\sigma(u)}+\frac{b(\theta,|V|)}{ \rho\sqrt{\theta}}[B(V)-B(R_xV)]\cdot\frac{\nabla_x \theta}{\sqrt{\theta}}\big)=0.
 \end{equation}
Direct calculations show that
 \begin{equation}
 \begin{aligned}
     \left[A(V)-A(R_xV)\right]\colon\sigma(u)
     &=\sum_{i,j=1}^{3}\left(V_iV_j-\left(V_i-2(V\cdot n)n_i\right)(V_j-2(V\cdot n )n_j)\right)\sigma_{ij}(u)\\
     &=\sum_{i,j=1}^{3}\left(4(V\cdot n)\left(V_i-(V\cdot n)n_i\right)n_j\right)\sigma_{ij}(u)\\
     &=4(V\cdot n)[V]^{\tan}\cdot\left(\sigma(u) n\right)\\
     &=4(V\cdot n)[(\nabla_x u+\nabla_x u^\mathrm{T})\cdot n]^{\tan} \cdot V,
 \end{aligned}
 \end{equation}
 where $[V]^{\tan}=V-(V\cdot n)n$ represents the tangential part of $V$. And we used the fact that $\sigma(u)$ is symmetric and
 \begin{equation*}
   \sum_{i,j=1}^{3}\left(4(V\cdot n)\left(V_i-(V\cdot n)n_i\right)n_j\right)\delta_{ij}\nabla_x\cdot u
   =\left(4(V\cdot n)\left(V-(V\cdot n)n\right)\cdot n\right)\nabla_x\cdot u
   =0.
 \end{equation*} 
 Then \eqref{G-LgammaG} reduces to
 \begin{equation}\label{G-LgammaG2}
   \frac{1}{({2\pi})^\frac{3}{2}\theta ^2}e^{-\frac{|V|^2}{2}}\Big\{2 a(\theta,|V|)(V\cdot n)[(\nabla_x u+\nabla_x u^\mathrm{T})\cdot n]^{\tan}V + b(\theta,|V|)(|V|^2-5)(V\cdot n)\frac{\nabla_x \theta\cdot n}{\sqrt{\theta}}\Big\}=0,
 \end{equation}
  for all $V\in \mathbb{R}^3$. It is direct to check that \eqref{G-LgammaG2} holds if and only if
 \begin{equation*}
   [S(u)n]^{\tan}=0, \quad \textrm{ and } \quad\nabla_x\theta\cdot n=0,  \quad \text{ on } \partial\Omega,
 \end{equation*}
 with $S(u)=\frac{1}{2}(\nabla_xu+\nabla_xu^\mathrm{T})$.
 Consequently, we obtain the complete slip type boundary conditions for the compressible Navier-Stokes system:
 \begin{equation}\label{slip1_formal}
   (u - u_w)\cdot n=0,\quad [S(u)n]^{\tan}=0, \quad \textrm{ and } \quad \nabla_x\theta\cdot n=0,  \quad \text{ on } \partial\Omega.
 \end{equation}
 Furthermore, we conclude that no boundary layer appears in this order. This observation underscores that specular reflection represents the most idealized form of reflection, devoid of boundary layer effects in this regime. As discussed in Subsection \ref{Subse:conclude} below, the Knudsen layer may appear in high order if we adopt the second-order Chapman-Enskog expansion.

\subsection{$\alpha_\eps= O(\eps^\beta)$ with $0< \beta < 1$}\label{subsec_almost_specular}
This case will be more complicate. As we will see below, the diffuse reflection part plays a crucial role, and the Knudsen layer must appear. Plug \eqref{chap} into the Maxwell reflection boundary condition \eqref{general_Maxwell_reflection} and rearrange by the order of $\eps$, we have 
\begin{align*}
   & O(1): \quad \gamma_ - M = L \gamma_+ M ,\\
   & O(\eps^\beta):\quad  L \gamma_+ M = K \gamma_+ M ,\\
   & O(\eps): \quad \gamma_ - G = L \gamma_+ G.
\end{align*}
We note that if $\beta = 1$, the last two equalities reduce to $\gamma_ - G  + L \gamma_+ M = K \gamma_+ M + L \gamma_+ G $. In both cases, we obtain that
\begin{equation}\label{dirichlet}
  u=u_w, \theta= \theta_w \,,
\end{equation}
as well as
\begin{equation}
  [(\nabla_x u + \nabla_x u ^{\mathrm{T}}) n ]^{\tan} =0\,, \quad \textrm{and}\quad \!  \nabla_x \theta \cdot n =0 .
\end{equation}
However, these constraints on the boundary are too many as the boundary conditions for the CNS. Therefore, this approach does not work. However, the relation \eqref{dirichlet} suggests that 
\begin{equation}\label{approximate}
  u - u_w =O(\eps^\iota) , \quad  \theta - \theta_w = O(\eps^\iota)
\end{equation}for some $\iota \geq 0$. Indeed, as shown below, $\iota = 1- \beta $ for $0< \beta \leq 1 $. In order to obtain the solution satisfying the boundary condition, we need to introduce the kinetic boundary layer, i.e., the Knudsen layer, with thickness $O(\eps)$ adjacent to the boundary \cite{aoki2017slip,sone2002kinetic,sone2007molecular}. Since the boundary condition \eqref{general_Maxwell_reflection} can be satisfied by the Chapman-Enskog solution \eqref{chap} up to $O(\eps)$ by the choice of \eqref{approximate}, we assume the the kinetic boundary layer $F^{bb}$ starts at the order of $\eps$. Hence, let 
\begin{equation}\label{Expansion_layer}
  \begin{aligned}
    F_\eps = M  &+ \eps G + O(\eps^2)\\
     & + \eps F^{bb} + O(\eps^2),
  \end{aligned}
\end{equation}
where $F^{bb}$ is appreciable only in the Knudsen layer, and it is assumed to decay rapidly away from the boundary.
\subsubsection{Knudsen layer equation}
First, we express a point $x_w$ on the boundary as a function of time $t$ and of coordinates $\chi_1$ and $\chi_2$ fixed on the surface of the boundary:
\begin{equation}
  x_w = x_w (t, \chi_1 , \chi_2).
\end{equation}
Then the velocity of the boundary $u_{wi}$ and the normal outward unit vector $n $ can be expressed by
\begin{align*}
  u_w(t,\chi_1,\chi_2)  &= \p_t x_w (t, \chi_1, \chi_2),\\
  n (t,\chi_1,\chi_2)  &= \pm \big( \p_{\chi_1} x_w \times  \p_{\chi_2} x_w \big) | \p_{\chi_1} x_w \times  \p_{\chi_2} x_w |^{-1},
\end{align*}
where $\p_{\chi_i} = \tfrac{\p  }{\p \chi_i } \, (i=1,2)$, and the sign $\pm$ is chosen in the way that $n$ points out the domain. In order to analyze the Knudsen layer, we need to introduce a new coordinate system that is local near the boundary and appropriate to describe the rapid change of the physical quantities in the direction normal to the boundary. Define the new variables $\tilde{t},\eta$ by the following relations:
\begin{equation}\label{3.7}
\begin{aligned}
 & t= \tilde{t},\\
 & x= - \eps  \eta n(\tilde{t},\chi_1,\chi_2)+x_w(\tilde{t},\chi_1,\chi_2),\\
  & v=v_w+u_w(\tilde{t},\chi_1,\chi_2).
\end{aligned}
\end{equation}
Here, $\eta$ is a stretched normal coordinate and $v_w$ is the molecular velocity relative to the velocity of the boundary. We assume that $F^{bb}$ is a function of $(\tilde{t},\eta,\chi_1,\chi_2,v_w)$ and vanishes as $\eta \rightarrow \infty$:
\begin{equation}\label{far_field_condition}
  \begin{aligned}
    F^{bb} &=F^{bb}(\tilde{t},\chi_1,\chi_2,\eta,v_w), \\
    F^{bb}  &\to 0 \quad \text{as} \quad \eta\to 0.
  \end{aligned}
\end{equation}

Put the ansatz \eqref{Expansion_layer} into the Boltzmann equation \eqref{be} and note that $M + \eps G $ already solves the Boltzmann equation, we obtain the following equation for $F^{bb}$:
  \begin{equation}\label{3.6}
    \eps \p_t F^{bb} + \eps  v \cdot \nabla_x F^{bb} = Q(M,F^{bb})+Q(F^{bb},M)+O(\eps ).
  \end{equation}

  Now, let us consider the region inside the Knudsen layer, that is, $\eta=O(1)$ or $|(x-x_w)\cdot n|=O(\eps )$. The dependence of $M$ on the variable $x$ is through the quantities $(\rho,u,\theta)$. Consequently, inside the Knudsen layer, the quantities $(\rho,u,\theta)$ can be expanded in a Taylor expansion around $x=x_w$, that is,
  \begin{equation}
    \rho=\rho_B+O(\eps \eta),\quad u=u_B+O(\eps  \eta),\quad \theta=\theta_B+O(\eps \eta),
  \end{equation}
  where the subscript B indicates the value on the boundary. Because $u_B=u_w+O(\eps^\iota )$ and $\theta_B=\theta_w+O(\eps^\iota )$, we can write\footnote{When $\iota = 0$, we straightforwardly skip this step.}
  \begin{equation}
    \rho=\rho_B+O(\eps \eta),\quad u=u_w+O(\eps^\iota \eta),\quad  \theta=\theta_w+O(\eps^\iota \eta).
  \end{equation}
  With these expressions, the local Maxwellian $M$ inside the Knudsen layer is expanded as
  $$M=M_w\big(1+O(\eps^\iota \eta)\big),$$
  where
  $$M_w=\frac {\rho_B}{{(2\pi \theta_w)}^{\frac {3}{2}} }\exp({-\frac{|v-u_w|^2}{2\theta_w}}).$$
  Therefore, the first and second term of the right-hand side of \eqref{3.6} can be written as follows:
  \begin{equation}\label{3.9}
    Q(M,F^{bb})+Q(F^{bb},M)=Q(M_w,F^{bb})+Q(F^{bb},M_w)+O(\eps^\iota \eta F^{bb}).
  \end{equation}
  Note that $\eta F^{bb}$ decay rapidly because of the rapid decay of $F^{bb}$. Now we need to express the left-hand side of \eqref{3.6} in terms of the new variables $(\tilde{t},\eta,\chi_1,\chi_2, v_w)$. From \eqref{3.7}, we have
\begin{equation}\label{3.11}
  \begin{aligned}
    \frac{\partial}{\partial t} & =\frac{\partial \tilde{t}}{\partial t} \frac{\partial}{\partial \tilde{t}}+\frac{\partial \eta}{\partial t} \frac{\partial}{\partial \eta}+\frac{\partial \chi_{1}}{\partial t} \frac{\partial}{\partial \chi_{1}}+\frac{\partial \chi_{2}}{\partial t} \frac{\partial}{\partial \chi_{2}}+\frac{\partial v_{wi}}{\partial t} \frac{\partial}{\partial v_{wi}}, \\
    \frac{\partial}{\partial x_{i}} & =\frac{\partial \tilde{t}}{\partial x_{i}} \frac{\partial}{\partial \tilde{t}}+\frac{\partial \eta}{\partial x_{i}} \frac{\partial}{\partial \eta}+\frac{\partial \chi_{1}}{\partial x_{i}} \frac{\partial}{\partial \chi_{1}}+\frac{\partial \chi_{2}}{\partial x_{i}} \frac{\partial}{\partial \chi_{2}}+\frac{\partial v_{wj}}{\partial x_{i}} \frac{\partial}{\partial v_{wj}}.
    \end{aligned}
\end{equation}
The formula \eqref{3.7} gives the following expressions for the coefficients in \eqref{3.11}(see \cite[Appendix A]{aoki2017slip}):
\begin{equation}\label{3.12}
  \begin{aligned}
    \frac{\partial \tilde{t}}{\partial {t}}=1, \quad \frac{\partial \eta}{\partial {t}}= \frac{1}{\varepsilon} u_{w} \cdot  n+O(1), \quad \frac{\partial \chi_{1,2}}{\partial {t}}=O(1), \quad \frac{\partial v_{wi}}{\partial {t}}=O(1),\\
    \frac{\partial \tilde{t}}{\partial x_{i}}=0, \quad \frac{\partial \eta}{\partial x_{i}}= - \frac{1}{\varepsilon} n_{i}+O(1), \quad \frac{\partial \chi_{1,2}}{\partial x_{i}}=O(1), \quad \frac{\partial v_{wj}}{\partial {x_i}}=O(1).
    \end{aligned}
\end{equation}
When we substitute \eqref{3.9} into \eqref{3.6} and combine \eqref{3.11} with \eqref{3.12}, we obtain the following equation:
\begin{equation}  \label{3.10}
  -v_w \cdot n \partial_{\eta}F^{bb}=Q(M_w,F^{bb})+Q(F^{bb},M_w)+O(\eps^\iota ).
  \end{equation}

Define $\xi=\frac{v-u_w}{\sqrt{\theta_w}} = \frac{v_w}{\sqrt{\theta_w}}$. Then $M_w$ can be expressed as
\begin{equation}\label{M_w}
  M_w=\frac{\rho_B}{\theta_w^{\frac{3}{2}}}\mu(\xi),
\end{equation}
where $\mu$ is the global Maxwellian defined in \eqref{mu}. If we let
\begin{equation}\label{3.15}
  F^{bb}(\tilde{t},\chi_1,\chi_2,\eta,v_w)=\frac{\rho_B}{\theta_w^{\frac{3}{2}}}\sqrt{\mu(\xi)} \tilde{f}^{bb}(\tilde{t},\chi_1,\chi_2,\eta,\xi),
\end{equation}
the first and second term on the right-hand side of \eqref{3.10} becomes
\begin{equation}\label{L_theta}
  \begin{aligned}
  &Q(M_w,F^{bb})+Q(F^{bb},M_w)\\
  =& \frac{\rho^2_B}{\theta_w} \int_{{\mathbb{R}^3} \times \mathbb{S}^2} \big(\mu(\xi'_1)\sqrt{\mu}(\xi')\tilde{f}^{bb}(\xi') +\mu(\xi')\sqrt{\mu}(\xi'_1)\tilde{f}^{bb}(\xi'_1) \\
  &- \mu(\xi_1)\sqrt{\mu}(\xi)\tilde{f}^{bb}(\xi)+\mu(\xi)\sqrt{\mu}(\xi_1)\tilde{f}^{bb}(\xi_1)\big) \frac{b( \sqrt{\theta_w}|\xi|, \omega)}{\sqrt{\theta_w}} d\xi_1 d \omega\\
  :=&-\frac{\rho^2_B}{\theta_w}\sqrt{\mu(\xi)} \L_{\theta_w}\tilde{f}^{bb}(\tilde{t},\chi_1,\chi_2,\eta,\xi).
  \end{aligned}
\end{equation}
In particular, for hard-sphere molecules, $\L_{\theta_w} = \L$, which is defined by \eqref{linearbol}. Combining \eqref{3.10}, \eqref{3.15} and \eqref{L_theta}, we finally get the Knudsen layer equation:
\begin{equation}\label{3.13}
  (\xi \cdot n) \partial_\eta \tilde{f}^{bb}(\xi)=\rho_B \L_{\theta_w} \tilde{f}^{bb}(\xi)+O(\eps^\iota).
\end{equation}
In order to get rid of $\rho_B$ in \eqref{3.13}, we introduce the new coordinate $y = \rho_B \eta   $ in place of $\eta$, and let 
\begin{equation}\label{fbb_tilde_fbb}
\tilde{f}^{bb}(\tilde{t},\chi_1,\chi_2,\eta,\xi)=f^{bb}(\tilde{t},\chi_1,\chi_2,y,\xi).    
\end{equation}
We have
\begin{equation}\label{3.4}
(\xi \cdot n) \partial_y f^{bb}(\xi)=\L_{\theta_w} f^{bb}(\xi)+O(\eps^\iota).
\end{equation}
If we neglect the terms of $O(\eps^\iota  )$, we obtain the equation for $f^{bb}$, i.e., that for $F^{bb}$:
\begin{equation}\label{3.5}
(\xi \cdot n) \partial_y f^{bb}(\xi)= \L_{\theta_w}f^{bb}(\xi).
\end{equation}
\subsubsection {Knudsen Layer Boundary Condition}
Now we consider the boundary condition. The solution $F_\eps$ in \eqref{Expansion_layer} should satisfy the general Maxwell boundary condition \eqref{general_Maxwell_reflection} up to $O(\eps)$. If we put \eqref{Expansion_layer} in \eqref{general_Maxwell_reflection}, we obtain that
\begin{equation}\label{4.25}
  (L^R-\alpha_\eps L^D)F^{bb}=-\frac{(L^R-\alpha_\eps L^D)M}{\eps }-(L^R-\alpha_\eps L^D)G+O(\eps ),
\end{equation}
where 
\begin{equation}\label{4.1}
L^R=\gamma_--L\gamma_+,
\end{equation}
with $L\gamma_+F(v)=F(v-2[(v-u_w)\cdot n]n)$ and
\begin{equation}\label{4.19}
  L^DF=K \gamma_+F-L\gamma_+F,
\end{equation}
with
\begin{equation}\label{4.23}
  K \gamma_+F(v)=\frac {1}{2\pi \theta_w^2}\exp({-\frac{|v-u_w|^2}{2\theta_w}})\int_{(v-u_w)\cdot n>0}{(v-u_w)\cdot nF}\,{\rm d}v .
\end{equation}
Note that 
$$L\gamma_+F(\xi)=F(\frac{v-2[(v-u_w)\cdot n]n-u_w}{\sqrt{\theta_w}})=F(\xi-2(\xi \cdot n)n).$$
It should be noted that the terms $\alpha_\eps L^D F^{bb}, \alpha_\eps L^D G$ are high-order term since $\alpha_\eps = \eps^\beta $ for $\beta >0$, and we should have neglected them. To compare the distinguishes with the case $\alpha_\eps= O(1)$, which has already been investigated in \cite{aoki2017slip}, we still retain these terms here and will eventually drop them. In accordance with \eqref{approximate}, we set 
\begin{equation}\label{u-u_w_theta-theta_w}
  u-u_w = \eps^\iota \hat{u}, \quad \textrm{and} \quad \! \theta - \theta_w = \eps^\iota \hat{\theta},
\end{equation}
where $\hat{u}$ and $\hat{\theta}$ are quantities of $O(1)$.

{\bf{Step 1. Calculation of $(L^R-\alpha_\eps L^D)M.$}} Notice that
$$\partial_uM=\frac{v-u}{M}, \quad \partial_\theta M=(\frac{|v-u|^2}{2\theta^2}-\frac{3}{2\theta})M.$$
Performing a Taylor's expansion around $u_w,\theta_w$ results in
\begin{equation}\label{4.22}
  M=\frac{\rho}{{(2\pi \theta)}^{\frac {3}{2}}} e^{{-\frac{|v-u|^2}{2\theta}}}
=M_w\big( 1+\eps ^\iota [\frac{ \hat{u} \cdot \xi}{\sqrt{\theta_w}}+(\frac{|\xi|^2-3}{2})\frac{\hat{\theta}}{\theta_w}]+O(\eps ^{2 \iota })\big).
\end{equation}
By \eqref{4.1}, it infers that
\begin{equation}\label{4.2}
  L^RM=(\gamma_--L\gamma_+)M= M_w \big( 2\eps^\iota \frac{(\xi \cdot n)(\hat{u} \cdot n)}{\sqrt{\theta_w}} + O(\eps^{2 \iota})  \big)
\end{equation}
Note that the right-hand side term $O(\eps^{2 \iota}) \sim \eps^{2 \iota} \hat{u} \cdot n O(1)$.
Moreover, we have
\begin{equation}\label{4.15}
\begin{aligned}
K \gamma_+M&=\frac{1}{2\pi\theta_w^2}e^{-\frac{|\xi|^2}{2}}\int_{\xi\cdot n>0} \theta_w^2 (\xi \cdot n)M d\xi\\
&=\frac{1}{2\pi}e^{-\frac{|\xi|^2}{2}} \int_{\xi\cdot n>0}(\xi\cdot n )\frac{\rho_B}{(2 \pi \theta_w )^{\frac{3}{2}}}e^{-\frac{|\xi|^2}{2}} \big( 1+\eps^\iota [\frac{\xi \cdot \hat{u}}{\sqrt{\theta_w}}+(\frac{|\xi|^2-3}{2})\frac{\hat{\theta}}{\theta_w}]\big)  + O(\eps^{2\iota})\\
&=\frac {1}{2\pi} M_w\int_{\xi \cdot n>0}{(\xi \cdot n)e^{-\frac{|\xi|^2}{2}}
\big( 1+\eps^\iota [\frac{\xi \cdot \hat{u}}{\sqrt{\theta_w}}+(\frac{|\xi|^2-3}{2})\frac{\hat{\theta}}{\theta_w}]\big)
} d \xi +O(\eps ^{2\iota}).
\end{aligned}
\end{equation}
For the right-hand side of \eqref{4.15}, easy computations show that:
\begin{equation}\label{4.16}
  \int_{\xi \cdot n>0}(\xi\cdot n)e^{-\frac{|\xi|^2}{2}}d\xi=\int_{-\infty}^{\infty}\int_{-\infty}^{\infty}\int_{0}^{\infty}\xi_3e^{-\frac{|\xi|^2}{2}}d\xi_1d\xi_2d\xi_3=2\pi,
\end{equation}

\begin{equation}\label{4.17}
  \begin{aligned}
  \int_{\xi\cdot n>0}(\xi \cdot n)(\xi \cdot \hat{u})e^{-\frac{|\xi|^2}{2}}d\xi &=\int_{\xi \cdot >0}(\xi \cdot n)[(\xi \cdot t)(\hat{u}\cdot t)+(\xi \cdot s)(\hat{u}\cdot s)+(\xi \cdot n)(\hat{u}\cdot n)]e^{-\frac{|\xi|^2}{2}}d\xi\\
  &=(\hat{u}\cdot n)\int_{\xi \cdot n>0}(\xi\cdot n)^2e^{-\frac{|\xi|^2}{2}}d\xi\\
  &=\frac{ {(2\pi)}^{\frac{3}{2}}}{2}(\hat{u}\cdot n),
\end{aligned}
\end{equation}
where $t,s$ are two unit vectors on the plane tangent to the boundary, which are orthogonal to each other. Moreover,
\begin{equation}\label{4.18}
  \int_{\xi \cdot n>0 }(\xi \cdot n)(\frac{|\xi|^2-3}{2})\frac{\hat{\theta}}{\theta_w}e^{-\frac{|\xi|^2}{2}}d\xi=\pi \frac{\hat{\theta}}{\theta_w}.
\end{equation}
Therefore, using \eqref{4.16}-\eqref{4.18}, \eqref{4.15} becomes
\begin{equation}\label{4.3}
  K \gamma_+M=M_w\big( 1+\frac{\sqrt{2\pi}\eps^{\iota}  (\hat{u}\cdot n)}{2\sqrt{\theta_w}}+\frac{\eps^\iota  \hat{\theta}}{2\theta_w}\big)+O(\eps ^{2\iota}).
\end{equation}
Using \eqref{4.22} and \eqref{4.3}, one has,
\begin{equation}
  \begin{aligned}
  L^DM=&(K\gamma_+-L\gamma_+)M\\
  =&M_w \big(1+\frac{\sqrt{2\pi}\eps^\iota  (\hat{u}\cdot n)}{2\sqrt{\theta_w}}+\frac{\eps^\iota  \hat{\theta}}{2\theta_w} \big)-M_w \big( 1+\frac{\eps^\iota [\xi-2(\xi \cdot n)n]\hat{u}}{\sqrt{\theta_w}}+\frac{|\xi|^2-3}{2}\frac{\eps^\iota \hat{\theta}}{\theta_w}\big)+O(\eps ^{2\iota})\\
  =&\eps^\iota  M_w \big( \frac{\sqrt{2\pi}\hat{u}\cdot n}{2\sqrt{\theta_w}}+\frac{\hat{\theta}}{2\theta_w}-\frac{\xi-2(\xi\cdot n)n}{\sqrt{\theta_w}}\hat{u}-\frac{|\xi|^2-3}{2}\frac{\hat{\theta}}{\theta_w}\big)+O(\eps ^{2\iota}).
\end{aligned}
\end{equation}
Together with \eqref{4.2}, one finally obtains that,
\begin{equation}\label{4.20}
  \begin{aligned}
    (L^R-\alpha_\eps L^D)M=M_w \eps^\iota \big(-\alpha_\eps \big(\frac{\sqrt{2\pi}\hat{u}\cdot n}{2\sqrt{\theta_w}}-\frac{[\xi-2(\xi \cdot n)n]\cdot \hat{u}}{\sqrt{\theta_w}}+(2-\frac{|\xi|^2}{2})\frac{\hat{\theta}}{\theta_w}\big)\\
    +\frac{2(\xi\cdot n)(\hat{u}\cdot n)}{\sqrt{\theta_w}}\big)+O(\eps ^{\beta+ 2\iota}).
  \end{aligned}
\end{equation}
If we integrate \eqref{3.10} over the whole space of $\xi$, we have 
\begin{equation}
\frac{\partial }{\partial \eta }\int_{\mathbb{R}^3 }^{}  (\xi\cdot n)F^{bb}d\xi =O(\eps^\iota ).
\end{equation}
Since $F^{bb}$ and the right-hand side term $O(\eps^\iota)$ vanishes rapidly as $\eta \to \infty$, the integration of the above equation with respect to $\eta$ from $0$ to $\infty$ leads to 
\begin{equation}\label{4.50}
\int_{\mathbb{R}^3 }   (\xi\cdot n)F^{bb}(\eta=0)d\xi =O(\eps^\iota ).
\end{equation}
On the other hand, the condition of no net flow on the boundary \eqref{flow_across_boundary} should be satisfied by $F_\eps$ of \eqref{Expansion_layer}. Substituting \eqref{Expansion_layer} into \eqref{flow_across_boundary} and using \eqref{4.50}, we have
\begin{equation}
\int_{\mathbb{R}^3 }  (v-u_w)\cdot n( M + \eps G)  (\eta=0)dv = O(\eps^{1 + \iota}) .
\end{equation}
Together with \eqref{4.22}, one has that
\begin{equation}\label{u_dot_n}
\rho_B\hat{u}\cdot n=O(\eps ).
\end{equation}
Therefore, \eqref{4.20} becomes
\begin{equation}\label{4.21}
  (L^R-\alpha_\eps L^D)M=-\alpha_\eps M_w\eps^\iota \big(-\frac{[\xi-(\xi \cdot n)n]\cdot \hat{u}}{\sqrt{\theta_w}}+(2-\frac{|\xi|^2}{2})\frac{\hat{\theta}}{\theta_w}\big)+O(\eps ^{\beta + 2 \iota} ).
\end{equation}

{\bf Step2. Calculation of $(L^R-\alpha_\eps L^D)G.$} Taylor's expansion around $u_w$ and $\theta_w$ yields that
\begin{equation}\label{4.8}
\begin{aligned}
G&=M_{(\rho,u,\theta)}\big(-\frac 12\frac{a(\theta,|V|)}{\rho\sqrt{\theta}}A(V):\sigma(u)-\frac{b(\theta,|V|)}{\rho{\theta}}B(V)\cdot \nabla_x \theta\big)\\
&=M_w \big(-\frac 12\frac{a(\theta_w, |\xi|)}{\rho_B\sqrt{\theta_w}}A(\xi):\sigma(u)_B-\frac{b(\theta_w, |\xi|)}{\rho_B{\theta_w}}B(\xi)\cdot (\nabla_x \theta)_B\big) +O(\eps^\iota ).
\end{aligned}  
\end{equation} 
{\bf Step 2.1. Calculation of $L^RG.$} Using \eqref{4.8}, one has
$$\begin{aligned}
L^RG=(\gamma_--L\gamma_+)G=M_w \big(&-\frac 12\frac{a(\theta_w, |\xi|)}{\rho_B\sqrt{\theta_w}}[A(\xi)-A(R_x\xi)]:\sigma(u)_B\\
&-\frac{b(\theta_w, |\xi|)}{\rho_B{\theta_w}}[B(\xi)-B(R_x\xi)]\cdot {(\nabla_x \theta)_B}\big)+O(\eps^\iota ).
\end{aligned}$$
Noting that $R_x\xi=\xi-2(\xi \cdot n)n$, by performing calculations analogous to those in the derivation of \eqref{G-LgammaG2}, we obtain
$$\begin{aligned} 
L^RG=M_w \big(&-2\frac{a(\theta_w,|\xi|)}{\rho_B\sqrt{\theta_w}}(\xi \cdot n)[\xi-(\xi \cdot n)n]\otimes n:[\nabla_xu+\nabla_xu^\mathrm{T}]_B \\
&-\frac{b(\theta_w,|\xi|)}{\rho_B{\theta_w}}(\xi \cdot n)(|\xi|^2-5) {n\cdot (\nabla_x \theta)_B}\big)+O(\eps^\iota ).
\end{aligned}
$$
{\bf Step 2.2 Calculation of $L^DG.$} This term is a high order term since $\alpha_\eps= O(\eps^\beta)$. In order to make a comparison with the case $\alpha_\eps= O(1)$, we also carry out the calculation of this term here.

Noting \eqref{4.8} and \eqref{4.23}, one has that,
\begin{equation}\label{4.7}
\begin{aligned} 
K \gamma_+G=&M_w\frac {1}{2\pi}\int_{\xi \cdot n>0}{(\xi \cdot n) e^{-\frac{|\xi|^2}{2}}[-\frac 12\frac{a(\theta_w , |\xi|)}{\rho_B\sqrt{\theta_w}}A(\xi):\sigma(u)]} d \xi + \\
&M_w\frac {1}{2\pi}\int_{\xi \cdot n>0}{(\xi \cdot n) e^{-\frac{|\xi|^2}{2}}[-\frac{b(\theta_w, |\xi|)}{\rho_B{\theta_w}}B(\xi)\cdot (\nabla_x \theta)_B]} d \xi +O(\eps^\iota ).
\end{aligned}
\end{equation}
Note that $\int_{\R^3} \hat{B}_3 \xi_3 \mu(\xi) d \xi =0$,  the second term on the right-hand side of \eqref{4.7} can be calculated as follows:
\begin{equation}
\begin{aligned}
 &M_w\frac {1}{2\pi}\int_{\xi \cdot n>0}{(\xi \cdot n ) e^{-\frac{|\xi|^2}{2}}[-\frac{b(\theta_w, |\xi|)}{\rho_B{\theta_w}}B(\xi)\cdot (\nabla_x \theta)_B]} d \xi \\
= & -\frac{M_w}{2\pi}\frac{(\nabla_x \theta)_B}{\rho_B\theta_w} \int_{\xi\cdot n>0}(\xi\cdot n)e^{-\frac{|\xi|^2}{2}}b(\theta_w, |\xi|)\frac{\xi(|\xi|^2-5)}{2}d\xi\\
= & -\frac{M_w}{2\pi}\frac{(\nabla_x \theta)_B \cdot n}{\rho_B\theta_w} \int_{\xi\cdot n>0}(\xi\cdot n)e^{-\frac{|\xi|^2}{2}}b(\theta_w, |\xi|)\frac{(\xi\cdot n)(|\xi|^2-5)}{2}d\xi\\
= & 0.
\end{aligned}
\end{equation}

Since $\xi$ can be expressed as $\xi=(\xi \cdot n)n+(\xi \cdot t)t+(\xi \cdot s) s$, where $t,s$ are two unit vectors on the plane tangent to the boundary, which are orthogonal to each other. For the first term on the right-hand side of \eqref{4.7}, one has,
\begin{equation}\label{4.9}
\begin{aligned}
  &M_w\frac {1}{2\pi}\int_{\xi \cdot n>0}{(\xi \cdot n) e^{-\frac{|\xi|^2}{2}}[-\frac 12\frac{a(\theta_w,|\xi|)}{\rho_B\sqrt{\theta_w}}A(\xi):\sigma(u)]} d \xi \\
=&-\frac{M_w}{4\pi\rho_B \sqrt{\theta_w}}\int_{\xi \cdot n>0}{(\xi \cdot n) e^{-\frac{|\xi|^2}{2}}a(\theta_w,|\xi|)(\xi_i\xi_j-\frac13\delta_{ij}|\xi|^2)} d \xi  \cdot(\partial_iu_j+\partial_ju_i)_B\\
=&-\frac{M_w}{4\pi\rho_B \sqrt{\theta_w}}\big(\int_{\xi \cdot n>0}{(\xi \cdot n)^3 e^{-\frac{|\xi|^2}{2}} a(\theta_w,|\xi|)} d \xi  \cdot n_in_j\\
&+(s_is_j+t_it_j)\int_{\xi \cdot n>0}{(\xi \cdot n)(\xi \cdot t)^2 e^{-\frac{|\xi|^2}{2}} a(\theta_w,|\xi|)} d \xi  \\
&-\frac{\delta_{ij}}{3}\int_{\xi \cdot n>0}{(\xi \cdot n) |\xi|^2 e^{-\frac{|\xi|^2}{2}} a(\theta_w,|\xi|)} d \xi  \big)(\partial_iu_j+\partial_ju_i)_B.\\
\end{aligned} 
\end{equation}
For the right-hand side of \eqref{4.9}, polar coordinate transformation yields that,
\begin{equation}\label{4.10}
\begin{aligned}
\int_{\xi \cdot n>0}{(\xi \cdot n)^3 e^{-\frac{|\xi|^2}{2}} a(\theta_w,|\xi|)} d \xi &=
\int_{-\infty}^{+\infty}\int_{-\infty}^{+\infty}\int_{0}^{+\infty}{\xi_3^3 e^{-\frac{|\xi|^2}{2}} a(\theta_w,|\xi|)} d \xi _1{\rm d}\xi_2{\rm d}\xi_3\\
&=\int_{0}^{2\pi}d\theta\int_{0}^{\frac{\pi}{2}}d\phi\int_{0}^{+\infty}(r\cos \phi)^3 e^{-\frac{|r|^2}{2}} a(\theta_w,|r|) r^2 \sin \phi dr\\
&=\frac{\pi}{2}\int_{0}^{+\infty}{r^5 e^{-\frac{|r|^2}{2}} a(\theta_w , |r|)}\,{\rm d}r  \\
&:= \frac{\pi}{2} \mathcal{F}_{\theta_w} ,
\end{aligned}
\end{equation}
Similarly, one gets
\begin{equation}\label{4.11}
  \begin{aligned}
  &\int_{\xi \cdot n>0}{(\xi \cdot n)(\xi \cdot t)^2 e^{-\frac{|\xi|^2}{2}} a(\theta_w,|\xi|)} d \xi\\ 
  = & \int_{0}^{2\pi}d\theta\int_{0}^{\frac{\pi}{2}}d\phi\int_{0}^{+\infty}(r\cos \phi)(r \sin \phi \cos \theta)^2 e^{-\frac{|r|^2}{2}}a(|r|)  r^2 \sin \phi dr \\
  = & \frac{\pi}{4}\mathcal{F}_{\theta_w} ,
\end{aligned}
\end{equation}
and 
\begin{equation}\label{4.12}
  \begin{aligned}
  \int_{\xi \cdot n>0}{(\xi \cdot n) |\xi|^2 e^{-\frac{|\xi|^2}{2}} a(\theta_w,|\xi|)} d \xi 
  &=\int_{0}^{2\pi}d\theta\int_{0}^{\frac{\pi}{2}}d\phi\int_{0}^{+\infty}(r\cos \phi)r^2 e^{-\frac{|r|^2}{2}}a(\theta, |r|)  r^2 \sin \phi dr\\
  &=\pi\mathcal{F}_{\theta_w} .
\end{aligned}
\end{equation}
As a result, \eqref{4.7} becomes 
$$K \gamma_+G=-\frac{M_w}{4\pi \rho_B \sqrt{\theta_w}}[\frac{\pi n_in_j}{6}-\frac{\pi}{12}(\delta_{ij}-n_in_j)]\cdot(\partial_iu_j+\partial_ju_i)_B \mathcal{F}_{\theta_w},$$
and therefore
\begin{equation}\label{4.13}
\begin{aligned}
L^DG=&K \gamma_+G-L\gamma_+G\\
=&-\frac{M_w}{4\pi \rho_B \sqrt{\theta_w}}[\frac{\pi n_in_j}{6}-\frac{\pi}{12}(\delta_{ij}-n_in_j)]\cdot(\partial_iu_j+\partial_ju_i)_B \mathcal{F}_{\theta_w}\\
&+ \frac{M_w}{\rho_B \sqrt{\theta_w}}[\frac 12 a(\theta_w,|\xi|)A(R_x\xi):\sigma(u)_B
+ b(\theta_w,|\xi|)B(R_x\xi)\cdot \frac{(\nabla_x \theta)_B}{\sqrt \theta_w}].
\end{aligned}  
\end{equation}
In order to deal with the second term on the right-hand side of \eqref{4.13}, we consider $[(\partial_iu_j)_B+(\partial_ju_i)_B](\delta_{ij}-n_in_j)$ and $[(\partial_iu_j)_B+(\partial_ju_i)_B](\delta_{ij}-n_in_j)(\delta_{ik}-n_in_k)(\delta_{jl}-n_jn_l).$
Since the derivatives are all in the tangential directions to the boundary, we can replace $u$ inside the derivatives with $u_w+\eps \hat{u}$. On the other hand, $u_w$ is the velocity of rigid body motion, so that the relation $\partial_iu_{wj}+\partial_ju_{wi}=0$ holds for any $i$ and $j$. Therefore, we have the following relations:
\begin{equation}\label{eq.2}
\begin{aligned} 
 &(\partial_iu_j+\partial_ju_i)_B(\delta_{ij}-n_in_j)=O(\eps ), \\
& (\partial_iu_j+\partial_ju_i)_B(\delta_{ik}-n_in_k)(\delta_{jl}-n_jn_l)=O(\eps ). 
  \end{aligned} 
\end{equation}
In consequence of \eqref{eq.2}, the following relation holds:
\begin{equation}
\begin{aligned}
&A(R_x\xi):\sigma(u)_B\\
=&\big([\xi_i-(\xi \cdot n)n_i-(\xi \cdot n)n_i][\xi_j-(\xi \cdot n)n_j-(\xi \cdot n)n_j]-\frac{\delta_{ij}}{3}|\xi|^2\big)(\partial_iu_j+\partial_ju_i)_B\\
=& \big([\xi_i-(\xi \cdot n)n_i][\xi_j-(\xi \cdot n)n_j]+[(\xi \cdot n)^2-\frac13|\xi|^2]n_in_j+\frac{|\xi|^2}{3}(n_in_j-\delta_{ij})\\
&-\big([\xi_i-(\xi \cdot n)n_i](\xi\cdot n)n_j+[\xi_j-(\xi \cdot n)n_j](\xi\cdot n)n_i\big)\big)(\partial_iu_j+\partial_ju_i)_B\\
=& 2n_in_j[(\xi \cdot n)^2-\frac13|\xi|^2](\partial_iu_j)_B-2(\xi \cdot n)[\xi-(\xi \cdot n)n]\otimes n:[\nabla_xu+\nabla_xu^\mathrm{T}]_B+O(\eps ).
\end{aligned}  
\end{equation}
We finally obtain that
\begin{equation}\label{4.24}
\begin{aligned}
(L^R-\alpha_\eps L^D)G= & \frac{M_w}{\rho_B \sqrt{\theta_w}}\big((-2+\alpha_\eps) a(\theta_w,|\xi|)(\xi \cdot n)[\xi-(\xi \cdot n)n]\otimes n:[\nabla_xu+\nabla_xu^\mathrm{T}]_B\\
&+(-2+\alpha_\eps)b(\theta_w,|\xi|)\frac{(\xi\cdot n)(|\xi|^2-5)}{2} \frac{(\nabla_x \theta)_B\cdot n}{\sqrt \theta_w}\\
&+\alpha_\eps\big(\frac{\mathcal{F}_{\theta_w}}{12}- a(\theta_w,|\xi|)[(\xi\cdot n)^2-\frac13|\xi|^2]\big)\partial_iu_jn_in_j\\
&-\alpha_\eps b(\theta_w,|\xi|)\frac{|\xi|^2-5}{2}\cdot[\xi-(\xi\cdot n)n] \frac{(\nabla_x \theta)_B}{\sqrt \theta_w}\big)+O(\eps ).
\end{aligned}
\end{equation}
Substituting \eqref{4.21} and \eqref{4.24} into \eqref{4.25}, the Knudsen layer boundary condition reduces to the following form
\begin{align*}
(L^R-\alpha_\eps L^D)F^{bb}= & M_w\big(-\alpha_\eps(\frac{|\xi|^2}{2}-2)\frac{\hat{\theta}}{\theta_w}-\alpha_\eps[\xi-(\xi\cdot n)n]\frac{\hat{u}}{\sqrt{\theta_w}}\big)\\
&+ \frac{M_w}{\rho_B \sqrt{\theta_w}}(2-\alpha_\eps) a(\theta_w,|\xi|)(\xi \cdot n)[\xi-(\xi \cdot n)]\otimes n:[\nabla_xu+\nabla_xu^\mathrm{T}]_B\\
&+ \frac{M_w}{\rho_B \sqrt{\theta_w}}(2-\alpha_\eps) b(\theta_w,|\xi|)\frac{(\xi\cdot n)(|\xi|^2-5)}{2} \frac{(\nabla_x \theta)_B\cdot n}{\sqrt \theta_w}\\
&+ \frac{M_w}{\rho_B \sqrt{\theta_w}}\alpha_\eps\big( a(\theta_w,|\xi|)[(\xi\cdot n)^2-\frac13|\xi|^2]-\frac{\mathcal{F}_{\theta_w}}{12}\big)\partial_iu_jn_in_j\\
&+ \frac{M_w}{\rho_B \sqrt{\theta_w}}\alpha_\eps b(\theta_w,|\xi|)\frac{|\xi|^2-5}{2} [\xi-(\xi\cdot n)n] \cdot \frac{(\nabla_x \theta)_B}{\sqrt \theta_w}+O(\eps^\iota ).
\end{align*}
Using \eqref{3.15} and \eqref{fbb_tilde_fbb}, we finally get the boundary condition of $f^{bb}$ 
\begin{equation}\label{4.100}
\begin{aligned}
\gamma_-f^{bb} = &(1-\alpha_\eps)L\gamma_+f^{bb}+\alpha_\eps \sqrt{2\pi }\sqrt{\mu}(\xi)\int_{\xi\cdot n>0}{(\xi\cdot n)f^{bb} \sqrt{\mu}(\xi)} d \xi \\
&+\sqrt{\mu(\xi)}\big(-\frac{\alpha_\eps \eps^\iota}{\eps}(\frac{|\xi|^2}{2}-2)\frac{\hat{\theta}}{\theta_w}-\frac{\alpha_\eps \eps^\iota}{\eps}[\xi-(\xi\cdot n)n]\frac{\hat{u}}{\sqrt{\theta_w}}\\
&+\frac{1}{\rho_B \sqrt{\theta_w}}(2-\alpha_\eps)  a(\theta_w,|\xi|)(\xi \cdot n)[\xi-(\xi \cdot n)]\otimes n:[\nabla_xu+\nabla_xu^\mathrm{T}]_B\\
&+\frac{1}{\rho_B \sqrt{\theta_w}}(2-\alpha_\eps) b(\theta_w,|\xi|)\frac{(\xi\cdot n)(|\xi|^2-5)}{2} \frac{(\nabla_x \theta)_B\cdot n}{\sqrt \theta_w}\\
&+\frac{1}{\rho_B \sqrt{\theta_w}}\alpha_\eps\big( a(\theta_w,|\xi|)[(\xi\cdot n)^2-\frac13|\xi|^2]-\frac{\mathcal{F}_{\theta_w}}{12}\big)(\partial_iu_j)_B n_in_j\\
&+\frac{1}{\rho_B \sqrt{\theta_w}}\alpha_\eps b(\theta_w,|\xi|)\frac{|\xi|^2-5}{2}[\xi-(\xi\cdot n)n] \cdot \frac{(\nabla_x \theta)_B}{\sqrt \theta_w}\big)+O(\eps^\iota ).
\end{aligned}
\end{equation}
\subsubsection{Summary and slip boundary conditions} 
Now we omit the high-order terms in \eqref{4.100} and take into account \eqref{3.5} and the far-field condition in \eqref{far_field_condition}. Then we obtain the following boundary value problem in the half-space $(y>0)$ for $f^{bb}$.
 \begin{equation}
  \left\{ \begin{aligned}\label{Knudsen1}
     (\xi \cdot n) \partial_y & f^{bb}(\tilde{t}, \chi_1,\chi_2, y, \xi) = \L_{\theta_w}f^{bb}(\tilde{t}, \chi_1,\chi_2, y,\xi) , \quad (y>0), \\
   \gamma_-f^{bb} = & L\gamma_+f^{bb} +\sqrt{\mu(\xi)}\big(-\frac{\alpha_\eps \eps^\iota}{\eps}(\frac{|\xi|^2}{2}-2)\frac{\hat{\theta}}{\theta_w}-\frac{\alpha_\eps \eps^\iota}{\eps}[\xi-(\xi\cdot n)n]\frac{\hat{u}}{\sqrt{\theta_w}}\\
   &+\frac{2}{\rho_B \sqrt{\theta_w}}  a(\theta_w, |\xi|)(\xi \cdot n)[\xi-(\xi \cdot n) n]\otimes n:[\nabla_xu+\nabla_xu^\mathrm{T}]_B\\
   &+\frac{2}{\rho_B \sqrt{\theta_w}} b(\theta_w, |\xi|)\frac{(\xi\cdot n)(|\xi|^2-5)}{2} \frac{(\nabla_x \theta)_B\cdot n}{\sqrt \theta_w} \big) ,\\ 
  \lim_{y \to \infty}  f^{bb}( &\tilde{t}, \chi_1,\chi_2, y , \xi) =0 ,
   \end{aligned}\right.
 \end{equation}
where $\L_{\theta_w}$ is defined in \eqref{L_theta}. Here we recall the subscript $B$ indicates the value on the boundary. The system \eqref{Knudsen1} is over-determined under the far-field condition  $\lim_{y\to \infty}f^{bb} = 0$, see \cites{bardos1986milne,golse1988boundary}. They are solvable only under certain solvability conditions. These solvibality conditions indeed give the boundary conditions for the compressible Navier-Stokes system. 
More precisely, consider the following Knudsen layer problem:
\begin{equation}\label{Kn}
  \left\{\begin{aligned}
   & ( \xi \cdot n) \p_y f = \L_{\theta_w}f,\\
   & \gamma_- f = L\gamma_+ f + h(\xi) ,\\
   & \lim_{y \to \infty} f =0.
  \end{aligned}\right.
\end{equation}
Golse, Perthame and Sulem \cite{golse1988boundary} proved that the solvability conditions of \eqref{Kn} is
\begin{equation}\label{solvability_condition}
  \int_{\xi \cdot n <0} (\xi \cdot n)   \begin{pmatrix}
   1\\
   \xi  \cdot t\\
   \xi  \cdot s\\
   \frac{|\xi|^2 -3}{2} 
   \end{pmatrix} h(\xi ) \sqrt{\mu}(\xi ) d \xi    =0,
\end{equation}
where $t,s$ are two unit vectors on the plane tangent to the boundary, which are orthogonal to each other.
From substituting the source term in \eqref{Knudsen1} into \eqref{solvability_condition}, we immediately derive that $\iota = 1 - \beta$, and  
\begin{equation}
  \int_{\xi \cdot n <0} (\xi \cdot n) (\xi \cdot t)^2 \mu(\xi) d \xi \frac{\chi (\hat{u} \cdot t)}{\sqrt{\theta_w}} = \int_{\xi \cdot n <0} \frac{2 a(\theta_w,|\xi|)}{\rho_B \sqrt{\theta_w}}{(\xi \cdot n)}^2 (\xi \cdot t)^2 \mu(\xi) d \xi [(\nabla_x u + \nabla_x u^{\mathrm{T}}) \cdot n]\cdot t,
\end{equation}
\begin{equation}
  \int_{\xi \cdot n <0} (\xi \cdot n) (\frac{|\xi|^2}{2} - 2)  (\frac{|\xi|^2 - 5}{2} ) \mu(\xi) d \xi \frac{\chi  \hat{\theta}  }{ {\theta_w}} = \int_{\xi \cdot n <0} \frac{2 b(\theta_w,|\xi|)}{\rho_B \sqrt{\theta_w}}{(\xi \cdot n)}^2 (\frac{|\xi|^2 - 5}{2})^2 \mu(\xi) d \xi \frac{(\nabla_x \theta)_B \cdot n}{\sqrt{\theta_w}}.
\end{equation}
Together with \eqref{u-u_w_theta-theta_w} and \eqref{u_dot_n}, we derive the slip boundary conditions for the CNS
\begin{equation}\label{Navier-slip}
 \left\{ \begin{aligned}
    (u-u_w )\cdot n &=0,\\
    [u - u_w]^{\tan} &= \tfrac{b^I_u}{ \chi \rho } \eps^{1 - \beta } [(\nabla_x u + \nabla_x u^{\mathrm{T}})\cdot n]^{\tan},\\
    \theta- \theta_w &= \tfrac{b^I_\theta}{ \chi \rho } \eps^{1 - \beta } (\nabla_x \theta \cdot n),
  \end{aligned}\right.
\end{equation}
where 
\begin{equation}\label{slip-coefficients1}
  \begin{aligned}
    & b^I_u = - \sqrt{2 \pi} \int_{\R^3}  { a(\theta_w,|\xi|)} {(\xi \cdot n)}^2 (\xi \cdot t)^2 \mu d \xi ,\\
    & b^I_\theta = - \tfrac{\sqrt{2 \pi}}{2} \int_{\R^3 }  {b(\theta_w,|\xi|)}{(\xi \cdot n)}^2 (\frac{|\xi|^2 - 5}{2})^2 \mu d \xi .
  \end{aligned}
\end{equation}
 We note that the slip boundary conditions \eqref{Navier-slip} are nonlinear due to the factor $\frac{1}{\rho}$.
\begin{remark}
 It should be noted that from \eqref{positivity}, $b^I_u$ and $b^I_\theta$ are negative constants depending on $\theta_w$. This is crucial in proving the well-posedness for \eqref{cns} with the boundary conditions \eqref{Navier-slip}.
\end{remark}

\begin{remark}
   In \eqref{4.100}, owing to the fact that $\alpha_\eps = O(\eps^\beta)$, the terms concerning $(\p_i u_j)_B n_in_j $ and $(\nabla_x \theta)_B$ in the last two lines are regarded as high-order terms and are ultimately omitted. This treatment differs fundamentally from the case where $\alpha = O(1)$: in the latter scenario, these terms remain $O(1)$ and the boundary conditions \eqref{Aoki_slip} can be derived, as discussed in \cite{aoki2017slip}.
\end{remark}

\subsection{$\alpha_\eps= \chi \eps$} The main distinguish between this case and the case $\alpha = \chi \eps^{\beta} (0<\beta<1)$ is that $u - u_w =O(1), \theta-\theta_w = O(1)$. So, Taylor's expansion around $(u_w,\theta_w)$ does not work. We directly consider the variable $V_B =\frac{v - u_B}{\theta_B}$. We define
$$F^{bb}(\tilde{t}, \chi_1, \chi_2 , \eta, v_w) = \frac{\rho_B}{\theta_B^{\frac{3}{2}}} \sqrt{\mu(V_B)}f^{bb}(\tilde{t}, \chi_1, \chi_2 , y, V_B). $$
By similar calculations as above, we obtain the following Knudsen layer problem:
\begin{equation}
  \left\{\begin{aligned}\label{Knudsen2}
    ( V_B  \cdot n) \partial_y & f^{bb}(\tilde{t}, \chi_1,\chi_2, y, V_B) = \L_{\theta_B}f^{bb}(\tilde{t}, \chi_1,\chi_2, y,V_B) , \quad (y>0), \\
  \gamma_-f^{bb} = & L\gamma_+f^{bb} +   \frac{\chi}{\sqrt{\mu(V_B)}} \big(- \mu(V_B) + \tfrac{ {\theta_B}^2}{ \theta_w^2} \mu(\xi) \big) \\
  &+\frac{2}{\rho_B \sqrt{\theta_B}} a(\theta_B, | V_B |)( V_B  \cdot n)[ V_B -( V_B  \cdot n) n]\otimes n:[\nabla_xu+\nabla_xu^\mathrm{T}]_B \sqrt{\mu( V_B )}\\
  &+\frac{2}{\rho_B \sqrt{\theta_B}} b(\theta_B, | V_B |)\frac{( V_B \cdot n)(| V_B |^2-5)}{2} \frac{(\nabla_x \theta)_B\cdot n}{\sqrt \theta_B}  \sqrt{\mu( V_B )}  ,\\ 
 \lim_{ y  \to \infty}  f^{bb}( &\tilde{t}, \chi_1,\chi_2, y ,  V_B ) =0 .
  \end{aligned}\right.
\end{equation}
Note that $V_B = \tfrac{\sqrt{\theta_w}\xi + u_w - u_B}{\sqrt{\theta_B}}$, by a change of variable, we have
\begin{equation*}
    \begin{aligned}
        &\int_{V_B \cdot n <0} (V_B\cdot n)(V_B\cdot t)\big( - \mu(V_B) + \tfrac{\theta^2_B}{\theta_w^2} \mu(\xi) \big) d V_B \\
        =&  \int_{\xi \cdot n <0} (\xi \cdot n)(\tfrac{(\sqrt{\theta_w}\xi + u_w - u_B)\cdot t}{\sqrt{\theta_B}} ) \mu(\xi)   d \xi \\
        = & - \tfrac{( u_w - u_B) \cdot t}{\sqrt{2 \pi \theta_B}},
    \end{aligned}
\end{equation*}
and 
\begin{equation*}
    \begin{aligned}
        &\int_{V_B \cdot n <0} (V_B\cdot n) \tfrac{|V_B|^2 }{2} \big( - \mu(V_B) + \tfrac{\theta^2_B}{\theta_w^2} \mu(\xi) \big) d V_B \\
        = &  \tfrac{\sqrt{2 \pi}}{\pi} + \tfrac{1}{2} \int_{\xi \cdot n <0} (\xi \cdot n) |\tfrac{(\sqrt{\theta_w}\xi + u_w - u_B)}{\sqrt{\theta_B}}|^2   \mu(\xi)   d \xi \\
        = &  \tfrac{\sqrt{2 \pi}}{\pi} \big(\tfrac{\theta_B - \theta_w   }{{\theta_B}} - \tfrac{ |u_w - u_B|^2 }{4 \theta_B} \big).
    \end{aligned}
\end{equation*}
Then the solvability conditions \eqref{solvability_condition} yield
\begin{equation}\label{Navier-slip1}
 \left\{ \begin{aligned}
    (u-u_w )\cdot n &=0,\\
    [u - u_w]^{\tan} &= \tfrac{c^I_u}{ \chi \rho }  [(\nabla_x u + \nabla_x u^{\mathrm{T}})\cdot n]^{\tan},\\
    \theta- \theta_w &= \tfrac{c^I_\theta}{ \chi \rho }  (\nabla_x \theta \cdot n) + \tfrac{|u - u_w|^2}{4},
  \end{aligned}\right.
\end{equation}
where the slip coefficients depends on $\theta_B$:
\begin{equation}\label{slip-coefficients2}
  \begin{aligned}
    & c^I_u = - \sqrt{2 \pi} \int_{\R^3}  { a(\theta_B,|V_B|)} {(V_B \cdot n)}^2 (V_B \cdot t)^2 \mu(V_B) d V_B ,\\
    & c^I_\theta = - \tfrac{\sqrt{2\pi}}{2} \int_{\R^3 }  {b(\theta_B,|V_B|)}{(V_B \cdot n)}^2 (\frac{|V_B|^2 - 5}{2})^2 \mu(V_B) d V_B .
  \end{aligned}
\end{equation}

\subsection{Concluding Remarks}\label{Subse:conclude}
In this study, we derive the slip boundary conditions for the CNS by taking into account the first-order Chapman-Enskog expansion. Our objective is to match the boundary conditions up to $O(\eps)$. In this context, for the cases where $\alpha_\eps = 0$ and $\alpha_\eps = O(\eps^\beta) \, (\beta>1)$, the Knudsen layer does not appear in $O(\eps)$. It appears in higher-order and we put it in the remainder term.  We say that the Knudsen layer is so weak that it is invisible in the first-order Chapman-Enskog expansion. If we adopt the second-order Chapman-Enskog expansion, the Knudsen layer always appears, at least in the high-order term. Nevertheless, the second-order expansion gives rise to the Burnett equations \cite{1970-Chapman}. These equations, however, are associated with unfavorable properties \cites{Bobylev-2006-JSP} and involve higher-order derivatives. 

\section{Uniform regularity for the compressible Navier-Stokes system with the complete slip boundary conditions}\label{Sec_CNS_uniform}
The complete slip boundary conditions \eqref{slip_theorem} with $u_w = 0$ reads
 \begin{equation}\label{slip}
   u \cdot n=0,\quad [S(u)n]^{\tan}=0, \quad \textrm{and} \quad \! \nabla_x\theta\cdot n=0,  \quad \text{ on } \partial\Omega,
 \end{equation}
 where $S(u) = \frac{1}{2} ( \nabla_x u + \nabla_x u ^{\mathrm{T}})$. For smooth solutions, the boundary conditions \eqref{slip} can be rewritten in the form of (see \cite{xiao20133d} for details)
    \begin{equation}\label{slip2}
       u \cdot n=0,\quad n\times \omega=- [2S(n)u]^{\tan}, \quad \textrm{ and }\quad\! \nabla_x\theta\cdot n=0,  \quad \text{ on } \partial\Omega,
    \end{equation}
 where $\omega=\nabla_x\times u$ denotes the vorticity.

This section is devoted to establishing the existence and higher-regularity of solutions to the compressible Navier-Stokes system with the complete slip boundary conditions \eqref{cns}, \eqref{cnsinitial}, and \eqref{slip}. The analysis presents significant challenges due to the weak dissipation effects on the right-hand side of \eqref{cns}, as well as the influence of boundary conditions. Obtaining uniform higher-regularity estimates for these solutions requires careful treatment. To overcome these difficulties, we employ techniques developed in previous work, particularly those in \cites{duan2021compressible, wang2016uniform, masmoudi2012uniform}. Key tools include elliptic estimates, Helmholtz decomposition, and the use of conormal derivatives. These methods allow us to systematically address the inherent analytical obstacles posed by the system.

For notational simplicity, and since the variable $v$ does not appear in this section, we adopt the following conventions: The spatial gradient $\nabla_x$ is simply denoted as $\nabla$, and the integral $\int_{\Omega} \cdot dx$ is abbreviated as $\int \cdot dx$. For $1\leq p \leq \infty$, we denote by $\|\cdot \|_p$ the standard $L^p(\Omega)$ norm.

To obtain the uniform regularity of the compressible Navier-Stokes system, we introduce the following conormal Sobolev space as in \cites{masmoudi2012uniform,wang2016uniform,duan2021compressible}. We consider a domain $\Omega \subseteq\mathbb{R}^3$ has a cover with the form $\Omega\in \Omega_0\cup_{i=1}^n\Omega_i$ where $\bar{\Omega}_0\subset\Omega$. In each subdomain $\Omega_i$, there exists a function $\psi _i$ such that
\begin{equation*}
  \begin{aligned}
    \Omega&\cap\Omega_i=\{x=(x_1,x_2,x_3)|x_3\geq \psi_i(x_1,x_2)\}\cap\Omega_i,\\
    \partial\Omega&\cap\Omega_i=\{x_3=\psi_i(x_1,x_2)\}\cap\Omega_i.
  \end{aligned}
\end{equation*}
Since $\partial\Omega$ is locally given by $x_3=\psi(x_1,x_2)$, we adopt the following corrdinate transformation for convenience:
\begin{equation*}
  \Psi: (y,z)\longmapsto (y,\psi(y)+z)=x,
\end{equation*}
where we omit the subscript $j$ of $\psi_j$ for notational simplicity. We thereby define $e_{y^1}=(1,0,\partial_1\psi)^\mathrm{T},e_{y^1}=(0,1,\partial_2\psi)^\mathrm{T}$ and $e_z=(0,0,-1)^\mathrm{T}$. Note that $e_{y^1}$ and $e_{y^2}$ are tangential to $\partial\Omega$, while $e_z$ is not a normal vector field in general. We then define the conormal derivatives
\begin{equation}\label{def_Z1Z2Z3}
  Z_i=\partial_{y^i}=\partial_i+\partial_i\psi\partial_z,\quad i=1,2, \quad Z_3=\varphi(z)\partial_z,
\end{equation}
where 
\begin{equation}\label{def_varphi}
  \varphi(z)=\frac{z}{1+z}
\end{equation}
is smooth, supported in $\mathbb{R}_+$ satisfying $\varphi(0)=0,\varphi'(0)\neq 0,\varphi(z)>0$ for $z>0$.
The unit outward normal $n$ is locally given by
\begin{equation}\label{N}
  n(x)\equiv n (\Psi(y,z))=\frac{1}{\sqrt{1+|\nabla_x \psi(y)|^2} }\begin{pmatrix}
    \partial _1\psi(y) \\
    \partial _2\psi (y) \\
    -1
    \end{pmatrix}
    := \frac{-N(y)}{\sqrt{1+|\nabla_x \psi(y)|^2} }.
\end{equation}
We introduce the notation
 \begin{equation}\label{def_mathcalZ}
  \mathcal{Z}^\alpha=\partial_t^{\alpha_0}Z^{\alpha_1}=\partial^{\alpha_0}_tZ_1^{\alpha_{11}}Z_2^{\alpha_{12}}Z_3^{\alpha_{13}},
\end{equation}
where $\alpha,\alpha_0,\alpha_1$ are multi-indices with $\alpha=(\alpha_0,\alpha_1)$, $\alpha_1=(\alpha_{11},\alpha_{12},\alpha_{13})$. We define the functional space
\begin{equation*}
  H^m_{co}=\{f\in L^2(\Omega),\quad Z^{\alpha_1} f\in L^2(\Omega),\quad |\alpha_1|\leq m\},
\end{equation*}
where $|\alpha_1|=|\alpha_{11}| + |\alpha_{12} |+ |\alpha_{13}|$ is a multi-index. We also set for smooth function $f(t,x)$
\begin{equation*}
  \|f\|^2_{H^m_{co}}=\sum_{|\alpha_1|\leq m}\|Z^{\alpha_1}f\|^2_{L^2_x},
 \quad
  \|f\|^2_{H^{\infty}_{co}}=\sum_{|\alpha_1|\leq m}\|Z^{\alpha_1}f\|^2_{L^\infty_x}.
\end{equation*}

To establish uniform regularity estimates for the solutions $(\rho,u,\theta)$ of the full compressible Navier-Stokes system, and particularly for applications in the remainder term analysis, we introduce a suitably chosen functional space. Our construction is motivated by the work of \cite{duan2021compressible}, where analogous structural and analytical challenges arise. Specifically, we define
$$\mathbb{X}(t)=\left\{(\rho,u,\theta);\|(\rho,u,\theta)(t)\|_{\mathbb{X}}\lesssim  \varepsilon^{\frac{3}{2}}\right\},$$
 where the norm $\|(\cdot,\cdot,\cdot)\|_{\mathbb{X}}$ is given by
\begin{equation}\label{cnsnorm}
  \begin{aligned}
    \|(\rho &, u, \theta)(t) \|_{\mathbb{X}}^{2} \\
    =&\sup _{\tau\in[0,t]} \sum_{\alpha_{0} \leq 3}\left\|\partial_{t}^{\alpha_{0}}(\bar{\rho}, u, \bar{\theta})(\tau)\right\|_{2}^{2}
    +\sup _{\tau\in[0,t]} \sum_{\alpha_{0} \leq 2}\left\|\partial_{t}^{\alpha_{0}} \nabla(\rho, u, \theta) (\tau)\right\|_{2}^{2} \\
    & +\sup _{\tau\in[0,t]} \sum_{\alpha_{0} \leq 2} \varepsilon^{2}\left\|\partial_{t}^{\alpha_{0}} \nabla^{2}(\rho,u, \theta)(\tau)\right\|_{2}^{2} +\sup _{\tau\in[0,t]} \sum_{\alpha_{0} \leq 2} \varepsilon\left\|\partial_{t}^{\alpha_{0}}({\rho}, u, {\theta}) (\tau)\right\|_{H^2_{co}}^{2}\\
    &+\sup _{\tau\in[0,t]} \sum_{\alpha_{0} \leq 1} \varepsilon^{2}\left\|\partial_{t}^{\alpha_{0}} \nabla(\rho, u, \theta) (\tau)\right\|_{H_{c o}^{2}}^{2}
    +\sup _{\tau\in[0,t]} \sum_{\alpha_{0} \leq 1} \varepsilon^{4}\left\|\partial_{t}^{\alpha_{0}} \nabla^{2} \rho(\tau)\right\|_{H_{c o}^{1}}^{2}\\
    & +\sup _{\tau\in[0,t]}  \varepsilon^{4}\left\| \nabla^{2} \rho(\tau)\right\|_{H_{c o}^{2}}^{2}
    +\sup _{\tau\in[0,t]} \sum_{\alpha_{0} \leq 1} \varepsilon^{4}\left\|\partial_{t}^{\alpha_{0}} \nabla^{2}(u, \theta)(\tau)\right\|_{H_{c o}^{2}}^{2},
    \end{aligned}
\end{equation}
where $\bar{\rho}= \rho -1$ and $\bar{\theta} = \theta -1$. We now state our main result concerning the existence, uniqueness, and uniform regularity of the solution to the full compressible Navier-Stokes system.
\begin{theorem}\label{Thmcns}
  Let $T>0$ be an arbitrary finite time. There exist $\lambda_0>0$ and $\varepsilon_0>0$ small such that if
  $$\|(\rho_0,u_0,\theta_0)\|^2_{\mathbb{X}}\leq \lambda_0^2\varepsilon^3,$$
  for $0<\varepsilon<  \varepsilon_0$, then the compressible Navier-Stokes system \eqref{cns} with the slip boundary condition \eqref{slip} admits a unique smooth solution $(\rho,u,\theta)$ on $[0,T]$. Moreover, the solution satisfies the uniform estimate:
  \begin{equation}\label{apriori}
    \Upsilon(t)\lesssim \left\|(\rho_{0}, u_{0}, \theta_{0})\right\|_{\mathbb{X}}^{2},
  \end{equation}
where
\begin{equation}
\begin{aligned}
\Upsilon(t):= &\|({\rho},u,{\theta})(t)\|^2_{\mathbb{X}} +\varepsilon \sum_{\alpha_{0} \leq 3} \int_{0}^{t}\left\|\partial_{t}^{\alpha_{0}} \nabla(u,\theta)(\tau)\right\|_{2}^{2}  d\tau +\varepsilon \sum_{\alpha_{0} \leq 2} \int_{0}^{t}\left\|\partial_{t}^{\alpha_{0}} \nabla \rho(\tau)\right\|_{2}^{2}  d\tau  \\
&+\varepsilon \sum_{\alpha_{0} \leq 2} \int_{0}^{t}\left\|\partial_{t}^{\alpha_{0}} \nabla^{2}(\rho, u, \theta)(\tau)\right\|_{2}^{2}  d\tau + \varepsilon^3 \sum_{\alpha_{0} \leq 2}\int_{0}^{t}\|\partial_{t}^{\alpha_{0}} \nabla^3 (u,\theta)(\tau)\|_2^2d\tau\\
&+\varepsilon^{2} \sum_{\alpha_{0} \leq 1} \int_{0}^{t}\left\|\partial_{t}^{\alpha_{0}} \nabla  {\rho} (\tau)\right\|_{H_{c o}^{2}}^{2}  d\tau   +\varepsilon^{2} \sum_{\alpha_{0} \leq 2} \int_{0}^{t}\left\|\partial_{t}^{\alpha_{0}} \nabla (u,  {\theta}) (\tau)\right\|_{H_{c o}^{2}}^{2}  d\tau \\
&+\varepsilon^{3} \sum_{\alpha_{0} \leq 1} \int_{0}^{t}\left\|\partial_{t}^{\alpha_{0}} \nabla^{2}( u,  {\theta}) (\tau)\right\|_{H_{c o}^{2}}^{2}  d\tau  
+\varepsilon^{3} \sum_{\alpha_{0} \leq 1} \int_{0}^{t}\left\|\partial_{t}^{\alpha_{0}} \nabla^{2}\rho(\tau)\right\|_{H_{c o}^{1}}^{2}  d\tau\\
&+\varepsilon^{3}  \int_{0}^{t}\left\|\nabla^{2}\rho(\tau)\right\|_{H_{c o}^{2}}^{2}  d\tau 
 +\varepsilon^{5} \sum_{\alpha_{0} \leq 1} \int_{0}^{t}\left\|\partial_{t}^{\alpha_{0}} \nabla^{3}(u,\theta)(\tau)\right\|_{H_{c o}^{2}}^{2}  d\tau.
\end{aligned}
\end{equation}
\end{theorem}
\begin{remark}
  $\p_t (\rho,u, \theta)$ does not make sense in general. However, it can be determined by the time evolution equations. The same principle applies to $\p_t R_0$ as discussed in Section \ref{Se:uniform}.
\end{remark}

We now proceed to establish the \emph{a priori} estimate \eqref{apriori} under the following bootstrap assumption:
\begin{equation}\label{assum}
  \Upsilon(t)\leq \lambda_0^2\varepsilon^3.
\end{equation}
Under the bootstrap assumption \eqref{assum}, we immediately derive the following key estimates:
\begin{lemma}\label{Lemma:rhoutheta_Linfty}
  Under the \emph{a priori} assumption \eqref{assum}, it holds that
  \begin{equation}\label{rhoinfty}
    \sup_{\tau\in[0,t]}\|\partial^{\alpha_0}_t(\bar{\rho},u,\bar{\theta})(\tau)\|_\infty\lesssim
    \sup_{\tau\in[0,t]}\|\partial^{\alpha_0}_t (\bar{\rho},u,\bar{\theta})(\tau)\|^{\frac{1}{2}}_{H^2}\|\partial^{\alpha_0}_t(\bar{\rho},u,\bar{\theta})(\tau)\|^{\frac{1}{2}}_{H^1}\lesssim
    \lambda_0\varepsilon,
  \end{equation}
  for $\alpha_0\leq 2$, and
  \begin{equation}\label{nablainf}
   \sup_{\tau\in[0,t]} \|\partial^{\alpha_0}_t\nabla(\bar{\rho},u,\bar{\theta})(\tau)\|_\infty\lesssim
   \sup_{\tau\in[0,t]} \|\partial^{\alpha_0}_t\nabla^2 (\bar{\rho},u,\bar{\theta})(\tau)\|^{\frac{1}{2}}_{H^1_{co}}\|\partial^{\alpha_0}_t\nabla(\bar{\rho},u,\bar{\theta})(\tau)\|^{\frac{1}{2}}_{H^2_{co}}\lesssim
    \lambda_0,
  \end{equation}
  for $\alpha_0\leq 1$.
\end{lemma}
\begin{proof}
  \eqref{rhoinfty} and \eqref{nablainf}  follow directly from \eqref{assum}, Proposition \ref{Pro:Agmon} and Proposition \ref{Pro:infty}.
\end{proof}
To obtain the high-order time and conormal energy estimates, we apply $\mathcal{Z}^\alpha$ (defined in \eqref{def_mathcalZ}) to  the compressible Navier-Stokes system \eqref{cns}. This yields
\begin{equation}\label{cocns}
  \begin{cases}
    \mathcal{Z}^\alpha\partial_t\rho+\mathcal{Z}^\alpha\nabla \cdot (\rho u)=0,\\
      \rho\mathcal{Z}^\alpha \partial_tu +\rho u\cdot \nabla\mathcal{Z}^\alpha u+\mathcal{Z}^\alpha\nabla(\rho\theta)
      =\frac{4}{3} \varepsilon \mu(\theta)\mathcal{Z}^\alpha \nabla\cdot u
      -\varepsilon\mu(\theta)\mathcal{Z}^\alpha\nabla\times \omega+\sum_{i=1}^{3}C^\alpha_i,\\
    \frac{3}{2}\rho\mathcal{Z}^\alpha\partial_t\theta+\frac{3}{2}\rho u \cdot \nabla \mathcal{Z}^\alpha\theta+\rho\theta\mathcal{Z}^\alpha\nabla \cdot u-\varepsilon\kappa(\theta)\mathcal{Z}^\alpha\Delta \theta=\sum_{i=4}^{7}C^\alpha_i,
  \end{cases}
\end{equation}
where $\omega=\nabla\times u$ is the vorticity and
\begin{equation}
\begin{cases}
      &C^\alpha_1=-[\mathcal{Z}^\alpha,\rho]\partial_t u=-\sum_{|\beta|\geq1,\beta+\gamma=\alpha}C_{\alpha,\beta}\mathcal{Z}^\beta\rho\mathcal{Z}^\gamma\partial_tu,\\
     & C^\alpha_2=-[\mathcal{Z}^\alpha,\rho u\cdot\nabla]u=-\sum_{|\beta|\geq1,\beta+\gamma=\alpha}C_{\alpha,\beta}\mathcal{Z}^\beta(\rho u)\mathcal{Z}^\gamma\nabla u,\\
     & C^\alpha_3=\eps \mathcal{Z}^\alpha[\nabla\mu(\theta)\sigma(u)]+\eps\sum_{|\beta|\geq1,\beta+\gamma=\alpha}C_{\alpha,\beta}\mathcal{Z}^\beta(\mu(\theta))\mathcal{Z}^\gamma(\frac{4}{3}\nabla\nabla\cdot u-\nabla\times \omega),\\
     & C^\alpha_4=-[\mathcal{Z}^\alpha,\rho]\partial_t \theta=-\sum_{|\beta|\geq1,\beta+\gamma=\alpha}C_{\alpha,\beta}\mathcal{Z}^\beta\rho \mathcal{Z}^\gamma\partial_t\theta,\\
     & C^\alpha_5=-[\mathcal{Z}^\alpha,\rho u\cdot\nabla]\theta=-\sum_{|\beta|\geq1,\beta+\gamma=\alpha}C_{\alpha,\beta}\mathcal{Z}^\beta(\rho u)\mathcal{Z}^\gamma\nabla \theta-\rho u\cdot[\mathcal{Z}^\alpha,\nabla]\theta,\\
     & C^\alpha_6=-[\mathcal{Z}^\alpha,\rho\theta]\nabla\cdot u=-\sum_{|\beta|\geq1,\beta+\gamma=\alpha}C_{\alpha,\beta}\mathcal{Z}^\beta(\rho \theta)\mathcal{Z}^\gamma\nabla\cdot u,\\
     & C^\alpha_7=\eps\mathcal{Z}^\alpha[\nabla\kappa(\theta)\cdot\nabla\theta]
      +\eps\sum_{|\beta|\geq1,\beta+\gamma=\alpha}C_{\alpha,\beta}\mathcal{Z}^\beta
      (\kappa(\theta))\mathcal{Z}^\gamma\Delta\theta
      +\frac{\varepsilon}{2}\mathcal{Z}^\alpha[\mu(\theta)\sigma(u):\sigma(u)].
\end{cases}
\end{equation}
Here, the square bracket notation denotes the standard commutator: $[A,B]=AB-BA$. Since the proof of \eqref{apriori} is complicated, we will divide our proof into several subsections.
\subsection{Classical Energy estimates for $\rho,u$ and $\theta$}
This subsection is devoted to deduce the estimates for the classical derivatives of $\rho,u$ and $\theta$. The weak dissipation mechanism in system \eqref{cns} presents particular challenges for deriving energy estimates, especially when dealing with higher-order derivatives. To overcome these difficulties, our analysis employs the following key techniques: the Helmholtz decomposition, elliptic estimates and careful treatment of nonlinear coupling terms. We begin by deriving the following $L^2$ estimates.
\begin{lemma}\label{Lm:L2}
Under the \emph{a priori} assumption \eqref{assum}, it holds that
\begin{multline}\label{L2}
      \sup _{\tau\in[0,t]} \sum_{\alpha_{0} \leq 3}\left\|\partial_{t}^{\alpha_{0}}(\rho-1, u, \theta-1)(\tau)\right\|_{2}^{2}+\varepsilon\sum_{\alpha_0\leq3}\int_{0}^{t}\|\partial^{\alpha_0}_t\nabla(u,\theta)\|_2^2d\tau\\
      +\varepsilon\sum_{\alpha_0\leq2}\int_{0}^{t}\|\partial^{\alpha_0}_t\nabla\rho\|^2_2d\tau\lesssim \Upsilon(0)+(\lambda_0+\varepsilon)\Upsilon(t).
\end{multline}
\end{lemma}
\begin{proof} Recall $(\bar{\rho},\bar{\theta})=(\rho-1,\theta-1)$. By the boundary conditions \eqref{slip2}, one has  
\begin{equation}\label{slip2t}
  \partial^{\alpha_0}_t u\cdot n=0,\qquad n\times \partial^{\alpha_0}_t\omega=-[2S(n)\partial^{\alpha_0}_tu]^{\tan},\qquad \nabla\partial^{\alpha_0}_t\theta \cdot n=0.
\end{equation}
Multiplying $\eqref{cocns}_1,\eqref{cocns}_2,$ and $\eqref{cocns}_3$ (with $\mathcal{Z}^\alpha=\partial^{\alpha_0}_t,\,  \alpha_0\leq 3)$ by $\partial_t^{\alpha_0}\bar{\rho},\partial_t^{\alpha_0}u$, and $\partial_t^{\alpha_0}\bar{\theta}$, respectively, one obtains that
\begin{equation}\label{rhoL2}
  \frac{1}{2}\frac{d}{dt}\int|\partial^{\alpha_0}_t\bar{\rho}|^2dx+\int\partial^{\alpha_0}_t\nabla\cdot(\bar{\rho}u)\partial^{\alpha_0}_t\bar{\rho} dx+\int\partial^{\alpha_0}_t\nabla\cdot u\partial_t^{\alpha_0}\bar{\rho} dx=0,
\end{equation}
\begin{multline}\label{uL2}
      \frac{1}{2}\frac{d}{dt}\int\rho|\partial^{\alpha_0}_tu|^2dx+\int\partial^{\alpha_0}_t\nabla(\bar{\rho}\bar{\theta}) \cdot \partial^{\alpha_0}_tudx+\int\partial^{\alpha_0}_t\nabla(\bar{\rho}+\bar{\theta}) \cdot \partial^{\alpha_0}_tudx\\
      =\frac{4}{3}\varepsilon\int\mu(\theta)\partial^{\alpha_0}_t\nabla\nabla\cdot u \cdot \partial^{\alpha_0}_tudx
      -\varepsilon\int\mu(\theta)\partial^{\alpha_0}_t\nabla\times\omega\cdot\partial^{\alpha_0}_tudx+\sum_{i=1}^{3}\int C^\alpha_i\partial^{\alpha_0}_t u dx,
\end{multline}
and
\begin{multline}\label{thetaL2}
    \frac{3}{4}\frac{d}{dt}\int\rho|\partial^{\alpha_0}_t\bar{\theta}|^2 dx+\int(\rho\theta-1)\partial^{\alpha_0}_t\nabla\cdot u\partial^{\alpha_0}_t\bar{\theta} dx+\int\partial^{\alpha_0}_t\nabla\cdot u\partial^{\alpha_0}_t\bar{\theta} dx\\
    =\eps\int\kappa(\theta)\partial^{\alpha_0}_t\Delta \theta \partial^{\alpha_0}_t\bar{\theta} dx+\sum_{i=4}^{7}\int C^\alpha_i\partial^{\alpha_0}_t\bar{\theta} dx.
\end{multline}
Recall \eqref{def_mu_kappa} for the definitions of $\mu(\theta)$ and $\kappa(\theta)$. For the dissipation of the velocity field $u$, we combine the boundary condition \eqref{slip2t} and Lemma \ref{Lemma:rhoutheta_Linfty}, and then by integrating by parts, we obtain
\begin{equation}\label{divu}
  \frac{4}{3}\varepsilon\int\mu(\theta)\partial^{\alpha_0}_t\nabla\nabla\cdot u\partial^{\alpha_0}_tudx
  \lesssim -\varepsilon\|\partial^{\alpha_0}_t\nabla\cdot u\|_2^2+\lambda_0\varepsilon\int|\partial^{\alpha_0}_t\nabla u|^2+|\partial^{\alpha_0}_t  u|^2dx,
\end{equation}
and
\begin{multline}\label{curlu}
      -\varepsilon\int\mu(\theta)\partial^{\alpha_0}_t\nabla\times\omega\partial^{\alpha_0}_tudx
      \lesssim-\varepsilon\|\partial^{\alpha_0}_t\nabla\times u\|_2^2+(\lambda_0+\delta_0)\varepsilon\int|\partial^{\alpha_0}_t\nabla u|^2dx\\
      +C_{\delta_0}\varepsilon\int |\partial^{\alpha_0}_t u|^2dx,
\end{multline}
for some $\delta_0>0$ small. Similarly, we derive the dissipation of the temperature
\begin{equation}\label{nablatheta}
  \varepsilon\int\kappa(\theta)\partial^{\alpha_0}_t\Delta \theta\partial^{\alpha_0}_t\bar{\theta} dx\lesssim
  -\varepsilon\int\kappa(\theta)|\partial^{\alpha_0}_t\nabla\bar{\theta}|^2dx+\lambda_0\varepsilon\|\partial^{\alpha_0}_t\bar{\theta}\|^2_{H^1}.
\end{equation}
Consequently, taking the summation of \eqref{rhoL2}, \eqref{uL2}, \eqref{thetaL2} and using \eqref{divu}, \eqref{curlu}, \eqref{nablatheta}, then integrating the resultant inequality over $[0,t]$, one obtains
\begin{equation}\label{L2*}
\begin{aligned}
  &\|\partial^{\alpha_0}_t (\bar{\rho},u,\bar{\theta})(t)\|^2_2+\varepsilon\int_{0}^{t}\|\nabla\partial^{\alpha_0}_t (u,\theta)\|_2^2d\tau+\int_{0}^{t}\int\partial^{\alpha_0}_t\bar{\rho}\partial^{\alpha_0}_t \nabla\cdot(\bar{\rho}u)dxd\tau\\
  &+\int_{0}^{t}\int\partial^{\alpha_0}_t u\nabla(\bar{\rho}\bar{\theta})dxd\tau+\int_{0}^{t}\int(\rho\theta-1)\partial^{\alpha_0}_t\bar{\theta}\nabla\cdot udxd\tau\\
  \lesssim&
  \|\partial^{\alpha_0}_t (\bar{\rho}_0,u_0,\bar{\theta}_0)\|^2_2
+\eps \int_{0}^{t}\|\partial^{\alpha_0}_t(u,\bar{\theta})\|_2^2d\tau
  +\sum_{i=1}^{3}\int_{0}^{t}\int C^\alpha_i\partial^{\alpha_0}_tudxd\tau
  \\
  &+\sum_{i=4}^{7}\int_{0}^{t}\int C^\alpha_i\partial^{\alpha_0}_t\bar{\theta}dxd\tau,
\end{aligned}
\end{equation}
where we used Proposition \ref{Pro:div} that
\begin{equation}
  \|\partial^{\alpha_0}_t\nabla u\|^2_2\leq \left(\|\partial^{\alpha_0}_t\nabla\times u\|^2_2+\|\partial^{\alpha_0}_t\nabla\cdot u\|^2_2+|\partial^{\alpha_0}_t u\cdot n|_{H^{\frac{1}{2}}}^2\right),
\end{equation}
as well as the trace theorem and H{\"o}lder's inequality. Now we estimate the third term on the left-hand side of \eqref{L2*}. Note that Lemma \ref{Lemma:rhoutheta_Linfty} only provides a second-order time derivative estimate for  $\|(\rho,u,\theta)\|_\infty$ and a first-order estimate for $\|\nabla(\rho,u,\theta)\|_\infty$. We cannot take $L^\infty$ norm for those high-order time derivative terms. Integrating by parts, using the boundary condition \eqref{slip2t} and Lemma \ref{Lemma:rhoutheta_Linfty}, one obtains
\begin{equation}
\begin{aligned}
    &\int_{0}^{t}\int\partial^{\alpha_0}_t\bar{\rho}\partial^{\alpha_0}_t \nabla\cdot(\bar{\rho}u)dxd\tau\\
    =&\int_{0}^{t}\int\partial^{\alpha_0}_t\bar{\rho}\partial^{\alpha_0}_t \nabla\bar{\rho}\cdot u dxd\tau    +\sum_{\gamma_0\geq 1,\beta_0+\gamma_0=\alpha_0}\int_{0}^{t}\int \partial^{\alpha_0}_t\bar{\rho}\partial^{\beta_0}_t\nabla\bar{\rho} \cdot \partial^{\gamma_0}_t u  dxd\tau \\
    &+ \quad\int_{0}^{t}\int\partial^{\alpha_0}_t\bar{\rho}\partial^{\alpha_0}_t (\bar{\rho}\nabla\cdot u)dxd\tau\\
    \gtrsim&  -\frac{1}{2}\int_{0}^{t}\int|\partial^{\alpha_0}_t\rho|^2 \nabla\cdot u dxd\tau-\lambda_0\int_{0}^{t}\left(\|\partial^{\alpha_0}_t(\bar{\rho},u)(\tau)\|_2^2+\sum_{\beta_0\leq\alpha_0} \varepsilon\|\nabla\partial^{\beta_0}_t(\bar{\rho},u)(\tau)\|_2^2  \right)d\tau\\
    \gtrsim&  -\lambda_0\int_{0}^{t}\left(\|\partial^{\alpha_0}_t(\bar{\rho},u)(\tau)\|_2^2+\sum_{\beta_0\leq\alpha_0} \varepsilon\|\nabla\partial^{\beta_0}_t(\bar{\rho},u)(\tau)\|_2^2  \right)d\tau.
\end{aligned}
\end{equation}
For the fourth term on the left-hand side of \eqref{L2*}, it follows from integrating by parts and H{\"o}lder's inequality that
\begin{multline}
      \int_{0}^{t}\int\partial^{\alpha_0}_t u \cdot \nabla(\bar{\rho}\bar{\theta})dxd\tau
      =-\int_{0}^{t}\int\partial^{\alpha_0}_t \nabla\cdot u\partial^{\alpha_0}_t(\bar{\rho}\bar{\theta})dxd\tau\\
      \gtrsim -\lambda_0\varepsilon  \int_{0}^{t}\|\partial^{\alpha_0}_t \nabla u\|_2^2 d\tau-\lambda_0\varepsilon\int_{0}^{t}\sum_{\beta_0\leq \alpha_0}\|\partial^{\beta_0}_t(\bar{\rho},\bar{\theta})(\tau)\|_2^2d\tau,
\end{multline}
where we consistently bound lower-order time derivatives in the  $L^\infty$ norm to leverage Lemma \ref{Lemma:rhoutheta_Linfty}.

In order to complete the estimates in \eqref{L2*}, one needs to estimate term concerning $C^\alpha_i$. By applying computational procedures to those demonstrated above, we have successfully derived
\begin{equation}\label{C123L2}
      \sum_{i=1}^{3}\int_{0}^{t}\int C^\alpha_i\partial^{\alpha_0}_tudxd\tau\leq \sum_{\beta_0\leq \alpha_0} \lambda_0\varepsilon\int_{0}^{t}\|\partial^{\beta_0}_t\nabla(u,\bar{\theta})\|_2^2d\tau+\lambda_0\int_{0}^{t}\|\partial^{\alpha_0}_t(u,\theta)(\tau)\|_2^2d\tau,
\end{equation}
and
\begin{equation}\label{C4-7L2}
    \sum_{i=4}^{7}\int_{0}^{t}\int C^\alpha_i\partial^{\alpha_0}_t\bar{\theta}dxd\tau
    \lesssim \sum_{\beta_0\leq \alpha_0} \lambda_0\varepsilon\int_{0}^{t}\|\partial^{\beta_0}_t\nabla(u,\bar{\theta})\|_2^2d\tau
    +\sum_{\beta_0\leq \alpha_0}(\lambda_0+\varepsilon)\int_{0}^{t}\|\partial^{\beta_0}_t(u,\theta)(\tau)\|_2^2d\tau.
\end{equation}
The details of the calculations are tedious but trivial, so we omit the details here. Substituting the estimates above into \eqref{L2} and choosing $\lambda_0$ and $\varepsilon_0$ suitably small, one proves
\begin{equation}\label{L2**}
  \sum_{\alpha_0\leq 3}\|\partial^{\alpha_0}_t(\bar{\rho},u,\bar{\theta})(t)\|_2^2+\sum_{\alpha_0\leq 3}\varepsilon\int_{0}^{t}\|\partial^{\alpha_0}_t\nabla(u,\bar{\theta})(\tau)\|_2^2d\tau
   \lesssim \Upsilon(0)+(\lambda_0+\varepsilon)\Upsilon(t).
\end{equation}
We then take the summation of $\varepsilon\left(\eqref{cocns}_1,\partial^{\alpha_0}_t\nabla\cdot u\right)$ and $\varepsilon\left(\eqref{cocns}_2,\partial^{\alpha_0}_t\nabla\rho/\rho\right)$ for $\alpha_0\leq 2$ to obtain the dissipation of $\partial_t^{\alpha_0}\nabla\bar{\rho}$, namely
\begin{multline}
    -\eps\int\partial^{\alpha_0}_t \nabla\cdot u\partial^{\alpha_0}_t\bar{\rho}dx+\eps\int_{0}^{t}\|\partial^{\alpha_0}_t\nabla\bar{\rho}\|_2^2d\tau\lesssim\eps\int\partial^{\alpha_0}_t \nabla\cdot u_0\partial^{\alpha_0}_t\bar{\rho}_0dx\\
    + \eps^3\int_{0}^{t}\|\nabla^2 u\|_2^2d\tau+\eps\int_{0}^{t}\|\partial^{\alpha_0}_t\nabla(u,\theta)\|_2^2d\tau,
\end{multline}
where we utilized the boundary condition \eqref{slip2t}. Hence, we have for $\alpha_0 \leq 2$
\begin{equation}\label{nablarho}
    \eps\int_{0}^{t}\|\partial^{\alpha_0}_t\nabla\bar{\rho}\|_2^2d\tau\lesssim \Upsilon(0)+\eps \Upsilon(t).
\end{equation}
As a consequence, \eqref{L2} follows from \eqref{L2**} and \eqref{nablarho}, thereby completing the proof.
\end{proof}

Next, we aim to derive the $H^1$ estimates of $(\rho,u,\theta)$. This is done by Lemma \ref{Lm:H1}.
\begin{lemma}\label{Lm:H1}
Under the a priori assumption \eqref{assum}, it holds that
\begin{multline}\label{H1}
  \sup _{\tau\in[0,t]} \sum_{\alpha_{0} \leq 2}\left\|\partial_{t}^{\alpha_{0}} \nabla(\rho, u, \theta)(\tau)\right\|_{2}^{2}
        +\sum_{\alpha_{0} \leq 2} \varepsilon\int_{0}^{t}\left\|\partial_{t}^{\alpha_{0}} \nabla^{2}( u, \theta)(\tau)\right\|_{2}^{2}d\tau \lesssim \Upsilon(0)+(\lambda_0+\varepsilon)\Upsilon(t).
\end{multline}
\end{lemma}
\begin{proof}  
  Multiplying  $\nabla\eqref{cocns}_1$ by $\partial^{\alpha_0}_t\nabla\rho$ (with $\mathcal{Z}^\alpha=\partial^{\alpha_0}_t,\alpha_0\leq 2$) and integrating yields that
  \begin{equation}\label{rhoH1}
    \frac{1}{2}\frac{d}{dt}\int|\partial^{\alpha_0}_t\nabla\rho|^2dx+\int\partial^{\alpha_0}_t\nabla\nabla\cdot(\bar{\rho}u) \cdot \partial^{\alpha_0}_t\nabla\rho dx +\int\partial^{\alpha_0}_t\nabla\nabla\cdot u \cdot \partial^{\alpha_0}_t\nabla\rho dx=0.
  \end{equation}
  Multiplying $\eqref{cocns}_2$ by $\partial^{\alpha_0}_t\nabla\nabla\cdot u$, we obtain
  \begin{equation}\label{uH1}
  \begin{aligned}
     & \frac{1}{2}\frac{d}{dt}\int\rho|\partial^{\alpha_0}_t\nabla\cdot u|^2dx
     +\int\nabla\rho\cdot\partial^{\alpha_0+1}_t u\partial^{\alpha_0}_t\nabla\cdot udx
     -\frac{1}{2}\int\partial_t\rho|\partial^{\alpha_0}_t\nabla\cdot u|^2dx\\
      &\quad- \int\rho u\cdot \partial^{\alpha_0}_t \nabla u\partial^{\alpha_0}_t\nabla\nabla\cdot udx
      - \int\partial^{\alpha_0}_t\nabla(\bar{\rho} \bar{\theta} + \bar{\rho} + \bar{\theta} )\partial^{\alpha_0}_t\nabla\nabla\cdot u dx   \\
      &\quad+\frac{4}{3}\varepsilon \int\mu(\theta)|\partial^{\alpha_0}_t\nabla\nabla\cdot u|^2dx
      - \int \varepsilon\mu(\theta)\partial^{\alpha_0}_t\nabla\times \omega\cdot\partial^{\alpha_0}_t\nabla\nabla\cdot udx\\
      &\quad+\sum_{i=1}^{3} \int C^\alpha_i\partial^{\alpha_0}_t\nabla\nabla\cdot udx =0.
  \end{aligned}
  \end{equation}
To handle the term $\int \varepsilon\mu(\theta)\partial^{\alpha_0}_t\nabla\times \omega\cdot\partial^{\alpha_0}_t\nabla\nabla\cdot udx$, we integrate with respect to $t$ and apply integration by parts. This yields
\begin{equation}\label{nabla2u}
\begin{aligned}
     &\eps \int_{0}^{t} \int\mu(\theta)\partial^{\alpha_0}_t\nabla\times \omega \cdot \partial^{\alpha_0}_t\nabla\nabla\cdot udx  d\tau\\
     =& \eps \int_{0}^{t} \int\mu(1)\partial^{\alpha_0}_t\nabla\times \omega \cdot \partial^{\alpha_0}_t\nabla\nabla\cdot udx  d\tau + \lambda_0 \eps \Upsilon(t)\\
=& \varepsilon \int_{0}^{t} \int_{\partial\Omega}\mu(1)n\times \partial^{\alpha_0}_t\omega\cdot\partial_t^{\alpha_0}\nabla\nabla\cdot u dS_x d \tau + \lambda_0 \eps \Upsilon(t)\\
=& - 2 \varepsilon \int_{0}^{t} \int_{\partial\Omega}\mu(1)[S(n)\partial_t^{\alpha_0} u]^{\tan}\cdot\partial^{\alpha_0}_t\nabla\nabla\cdot udS_x d \tau + \lambda_0 \eps \Upsilon(t)\\
\lesssim& \eps \int_{0}^{t}\|\partial^{\alpha_0}_t\nabla u \|_{H^{\frac{1}{2}}}^2 d \tau  + \lambda_0 \eps \Upsilon(t)\\
\lesssim & \delta_0 \eps \int_{0}^{t}\|\partial^{\alpha_0}_t \nabla^2 u \|_2^2 d \tau  + C_{\delta_0}(\Upsilon(0) + (\lambda_0 +\eps) \Upsilon(t)).
\end{aligned}
\end{equation}
Here $\delta_0$ represents a constant suitably small. The second last line follows because for the tangential component $[S(n)\partial_t^{\alpha_0} u]^{\tan}$, the derivative in $\nabla \nabla\cdot u$ reduces to the tangential derivative $Z_j \nabla\cdot u $ with $j=1,2$, where $Z_j$ are tangential derivatives defined in \eqref{def_Z1Z2Z3}. Based on this reduction, we can perform integration by parts. By further applying the trace theorem, the validity of the step in the second last line is established.  For the last inequality, we employ Lemma \ref{Lm:L2}.

To estimate the temperature part, we apply $\nabla$ to \eqref{cocns} and take inner product with $\partial^{\alpha_0}_t\nabla\theta$ to obtain
\begin{equation}\label{thetaH1}
  \begin{aligned}
   & \frac{3}{4} \frac{d}{dt}\int\rho|\partial^{\alpha_0}_t\nabla\theta|^2dx+  \int(\rho\theta - 1 )\partial^{\alpha_0}_t\nabla\nabla\cdot u\cdot\partial^{\alpha_0}_t\nabla\theta dx + \int \partial^{\alpha_0}_t\nabla\nabla\cdot u\cdot\partial^{\alpha_0}_t\nabla\theta dx \\
   &-\varepsilon  \int\nabla\left(\kappa(\theta)\partial^{\alpha_0}_t\Delta\theta\right)\cdot\partial^{\alpha_0}_t\nabla\theta dx  \\
    =& -\int(\frac{3}{2}\nabla\rho\partial^{\alpha_0+1}_t\theta+\nabla(\rho u)^\mathrm{T}\partial^{\alpha_0}_t\nabla\theta+\nabla(\rho\theta)\partial^{\alpha_0}_t\nabla\cdot u)\partial^{\alpha_0}_t\nabla\theta dx   +\sum_{i=4}^{7}  \int\nabla C^\alpha_i\partial^{\alpha_0}_t\nabla\theta dx.
  \end{aligned}
\end{equation}
Integrating by parts and using the boundary condition \eqref{slip2t} yields
\begin{equation}\label{nabla2theta}
  -\varepsilon \int_{0}^{t} \int\nabla\left(\kappa(\theta)\partial^{\alpha_0}_t\Delta\theta\right)\cdot\partial^{\alpha_0}_t\nabla\theta dxd\tau
  =\varepsilon \int_{0}^{t} \int\kappa(\theta)|\partial^{\alpha_0}_t\Delta\theta|^2 dx d \tau.
\end{equation}
Summing \eqref{rhoH1}, \eqref{uH1} and \eqref{thetaH1} and applying \eqref{nabla2u}, \eqref{nabla2theta}, we then obtain
\begin{equation}\label{H1'}
  \begin{aligned}
    &\|\partial_t^{\alpha_0}(\nabla\rho,\nabla\cdot u,\nabla\theta)\|^2_2
    +\varepsilon\int_{0}^{t}\|\partial^{\alpha_0}_t(\nabla\nabla\cdot u,\nabla^2\theta)\|_2^2d\tau
    + \int_{0}^{t} \int\partial^{\alpha_0}_t\nabla\nabla\cdot(\bar{\rho}u) \cdot \partial^{\alpha_0}_t\nabla\rho dx d \tau\\
    &+\int_{0}^{t} \int\nabla\rho\cdot\partial^{\alpha_0+1}_t u\partial^{\alpha_0}_t\nabla\cdot udx d \tau
    -\frac{1}{2}\int_{0}^{t} \int\partial_t\rho|\partial^{\alpha_0}_t\nabla\cdot u|^2dx d \tau - \int_{0}^{t} \int\rho u\cdot \partial^{\alpha_0}_t \nabla u\partial^{\alpha_0}_t\nabla\nabla\cdot udx d \tau \\
    & - \int_{0}^{t} \int\partial^{\alpha_0}_t\nabla(\bar{\rho} \bar{\theta} )\partial^{\alpha_0}_t\nabla\nabla\cdot u dx d \tau
     + \int_{0}^{t} \int(\rho\theta - 1 )\partial^{\alpha_0}_t\nabla\nabla\cdot u\cdot\partial^{\alpha_0}_t\nabla\theta dx d\tau  \\
      &+\sum_{i=1}^{3} \int_{0}^{t}\int C^\alpha_i\partial^{\alpha_0}_t\nabla\nabla\cdot udx d \tau -\sum_{i=4}^{7} \int_{0}^{t} \int\nabla C^\alpha_i\partial^{\alpha_0}_t\nabla\theta dx d \tau\\
      &  +\int_{0}^{t} \int(\frac{3}{2}\nabla\rho\partial^{\alpha_0+1}_t\theta+\nabla(\rho u)^\mathrm{T}\partial^{\alpha_0}_t\nabla\theta+\nabla(\rho\theta)\partial^{\alpha_0}_t\nabla\cdot u)\partial^{\alpha_0}_t\nabla\theta dx d\tau  \\
    &\lesssim \delta_0 \eps \int_{0}^{t}\|\partial^{\alpha_0}_t \nabla^2 u \|_2^2 d \tau  + C_{\delta_0}(\Upsilon(0) + (\lambda_0 +\eps) \Upsilon(t)),
  \end{aligned}
\end{equation}
where we used Proposition \ref{Pro:div}  
\begin{equation*}
  \|\nabla^2\theta\|_2\lesssim \|\Delta\theta\|_2 +|\nabla\theta\cdot n|_{H^{\frac{1}{2}}}.
\end{equation*}
Most of the terms in \eqref{H1'} can be estimated similarly to Lemma \ref{Lm:L2}. For example, we have
\begin{equation*}
  \begin{aligned}
    \int \nabla \rho \cdot \p^{\alpha_0 +1}_t u \p^{\alpha_0}_t \nabla \cdot u dxd\tau \lesssim& \lambda_0 \int_{0}^{t}\|\p^{\alpha_0+1}_t u\|_2 \|\p^{\alpha_0}_t \nabla\cdot u\|_2 d\tau\\
    \lesssim& \lambda_0 \int_{0}^{t}\|\p^{\alpha_0+1}_t u\|_2^2 d\tau + \lambda_0 \int_{0}^{t}  \|\p^{\alpha_0}_t \nabla\cdot u\|^2_2  d \tau\\
    \lesssim& \delta_0 \sup_{\tau \in [0,t]} \|\p^{\alpha_0}_t \nabla\cdot u(\tau)\|_2^2 + \Upsilon(0)   + (\eps + \lambda_0) \Upsilon(t) 
  \end{aligned}
\end{equation*}
 by taking $\lambda_0$ small. Consequently, we have for $\alpha_0 \leq 2$,
\begin{equation}\label{H1*}
  \begin{aligned}
    &\sup_{\tau \in [0,t]} \sum_{\alpha_0 \leq 2}\|\partial_t^{\alpha_0}(\nabla\rho,\nabla\cdot u,\nabla\theta)(\tau)\|^2_2 
    +\sum_{\alpha_0 \leq 2}\varepsilon\int_{0}^{t}\|\partial^{\alpha_0}_t(\nabla\nabla\cdot u,\nabla^2\theta)\|_2^2d\tau \\
    \lesssim & \delta_0 \eps \int_{0}^{t}\|\partial^{\alpha_0}_t \nabla^2 u \|_2^2 d \tau  + C_{\delta_0}(\Upsilon(0) + (\lambda_0 +\eps) \Upsilon(t)).
  \end{aligned}
\end{equation}

To obtain the estimates of $\|\partial^{\alpha_0}_t\nabla\times u\|_2^2$ and $\int_{0}^{t}\|\partial^{\alpha_0}_t\nabla^2 u\|_2^2d\tau$, we use the Helmholtz decomposition. Let $P_1$ be the curl-free projection operator and $P_2$ be the divergence-free operator defined in Proposition \ref{Pro:Helm}. We split
 \begin{equation}\label{decomposition_u}
  u=u_1+u_2,\quad \text{with} \quad u_1=P_1 u=\nabla\mathfrak{u}, \quad u_2=P_2u=\nabla\times\mathfrak{v}.
\end{equation}
Acting $P_2$ to $\eqref{cocns}_2$ (with $ \mathcal{Z}^\alpha=\partial^{\alpha_0}_t,\alpha_0\leq 2$), we obtain
\begin{equation}\label{omegaL2}
\begin{aligned}
      &\partial^{\alpha_0+1}_{t} u_{2}+P_{2}\left\{\bar{\rho} \partial^{\alpha_0+1}_{t} u\right\}
      +P_{2}\left\{\rho u \cdot \partial^{\alpha_0}_t\nabla u\right\}
      +\varepsilon \mu(1)  \partial^{\alpha_0}_{t}\nabla \times \omega \\
      =& \frac{4 \varepsilon}{3} P_{2}\left\{\left(\mu\left(\theta \right)-\mu(1)\right) \partial^{\alpha_0}_{t}\nabla \nabla \cdot u\right\}
      -\eps P_{2}\left\{\left(\mu\left(\theta \right)-\mu(1)\right)\partial^{\alpha_0}_{t} \nabla \times \omega\right\} +\sum_{i=1}^{3}P_2 C^{\alpha}_i.
\end{aligned}
\end{equation}
Taking the inner product of \eqref{omegaL2} with  $\partial^{\alpha_0+1}_{t} u$  and integrating the resultant equation with respect to  $t$. Note that
\begin{equation*}
  \begin{aligned}
        \varepsilon\int\partial^{\alpha_0}_t \nabla\times \omega\cdot\partial^{\alpha_0+1}_tudx
        &=\varepsilon\frac{1}{2}\frac{d}{dt}\int|\partial^{\alpha_0}_t\omega|^2dx-\varepsilon\int_{\partial\Omega}2[S(n)\partial_t^{\alpha_0}u]^{\tan}\partial^{\alpha_0+1}_tudS_x\\
      &=\varepsilon\frac{1}{2}\frac{d}{dt}\|\partial^{\alpha_0}_t\omega\|^2
      -\frac{d}{dt}\int_{\partial\Omega}[S(n)\partial_t^{\alpha_0}u]\cdot\partial^{\alpha_0}_tudS_x,
  \end{aligned}
\end{equation*}
since $S(n)$ is sysmmetric and $[\partial^{\alpha_0}_tu]^{\tan}=\partial^{\alpha_0}_tu$. Combined with Proposition \ref{Pro:Helm}, we obtain
\begin{equation}\label{omegaL21}
\begin{aligned}
   & \varepsilon\|\partial^{\alpha_0}_t\omega\|^2_2+ \int_{0}^{t}\left\|\partial^{\alpha_0+1}_{t} u_2\right\|_{2}^{2} d \tau
     +\eps\int_{\partial\Omega} [S(n)\partial^{\alpha_0}_tu_0]\cdot \partial_t^{\alpha_0}u_0dS_x \\
    \lesssim & \varepsilon\left\|\partial^{\alpha_0}_t\omega_0\right\|_2^2
    +\varepsilon|\partial^{\alpha_0}_t u|^2_{L^2(\partial\Omega)}
    + \int_{0}^{t}\left\|\bar{\rho} \partial^{\alpha_0+1}_{t} u\right\|_{2}^{2} d\tau
    + \int_{0}^{t}\left\|\rho u \cdot \partial^{\alpha_0}_{t}\nabla u\right\|_{2}^{2} d\tau \\
    & +\varepsilon^{2} \int_{0}^{t}\left\|\left(\mu\left(\theta \right)-\mu(1)\right) \partial^{\alpha_0}_{t}\nabla \nabla \cdot u\right\|_{2}^{2} d\tau \\
    & + \varepsilon^{2} \int_{0}^{t}\left\|\left(\mu\left(\theta \right)-\mu(1)\right)\partial^{\alpha_0}_{t} \nabla \times \omega\right\|_{2}^2 d\tau +\sum_{i=1}^{3}\int_{0}^{t}\|C^\alpha_i\|^2_2d\tau\\
    \lesssim &  \varepsilon \Upsilon(0)+ \lambda_0^{2} \varepsilon \Upsilon(t),
\end{aligned}
\end{equation}
Utilizing \eqref{assum} and the trace theorem, \eqref{omegaL21} reduces to
\begin{equation}\label{omegaL2*}
  \varepsilon\|\partial^{\alpha_0}_t\omega\|^2_2+ \int_{0}^{t}\left\|\partial^{\alpha_0+1}_{t} u_2\right\|_{2}^{2} d \tau
  \lesssim  \varepsilon \Upsilon(0)+ \lambda_0^{2} \varepsilon \Upsilon(t).
\end{equation}  

Next, we reformulate \eqref{cns} as the following elliptic system:
\begin{equation}
 \begin{cases}
     \varepsilon \nabla \cdot u= \mathfrak{g}_1, \\
     -\varepsilon \mu(1) \Delta u+\nabla \bar{P} = \mathfrak{g}_2, \\
     \left.u\cdot n\right|_{\partial \Omega}=0,\\
     \left.[S(u) n]^{\tan}\right|_{\partial \Omega}=0,
 \end{cases}
\end{equation}
with
\begin{equation*}
\begin{aligned}
  \mathfrak{g}_1 = &-\varepsilon \partial_{t} \rho-\varepsilon \bar{\rho} \nabla \cdot u-\varepsilon \nabla \rho \cdot u,   \\
  \mathfrak{g}_2 = & -\partial_{t} u_2-\bar{\rho} \partial_{t} u-\rho u \cdot \nabla u+\varepsilon\left(\mu\left(\theta \right)-\mu(1)\right) \Delta u \\
  & +\varepsilon \nabla \mu\left(\theta \right) \cdot \sigma(u)+\frac{\varepsilon}{3}\left(\mu\left(\theta \right)-\mu(1)\right) \nabla \nabla \cdot u,  \\
    \bar{P} = & \partial_t \mathfrak{u}+\rho\theta-\frac{\varepsilon}{3}\mu(1)\nabla\cdot u,
\end{aligned}
\end{equation*}
where $\mathfrak{u}$ is defined in \eqref{decomposition_u}. Then, using Lemma \ref{Lm:ADN}, we have for $0\leq \alpha_0\leq2$, 
\begin{equation}
\begin{aligned}
  \varepsilon^2\|\nabla^2\partial^{\alpha_0}_t u\|^2_2
 &\lesssim\|\partial^{\alpha_0}_t \mathfrak{g}_1\|^2_{H^1}+\|\partial^{\alpha_0}_t \mathfrak{g}_2\|^2_2
  \\ &\lesssim \varepsilon^{2}\left\|\partial_{t}^{\alpha_{0}+1} \rho\right\|_{H^{1}}^{2}
      +\lambda_0^{2} \varepsilon^2 \sum_{\beta_{0}  \leq \alpha_0}\left\|\nabla \partial_{t}^{\beta_{0}}(\rho, u,\theta) \right\|_{2}^{2}
      +\lambda_0^{2} \varepsilon^{4} \sum_{\beta_{0}\leq \alpha_0}\left\|\nabla^{2} \partial_{t}^{\beta_{0}} (\bar{\rho}, u) \right\|_{2}^{2}\\
     & \quad+\left\|\partial_{t}^{\alpha_{0}+1} u_2\right\|_{2}^{2}
      +\lambda_0^2\eps^2 \sum_{\beta_{0}\leq \alpha_0}\|\partial^{\beta_0+1}_t u\|_2^2 .
\end{aligned}
\end{equation}
We note that 
\begin{equation}
  \begin{aligned}
    \int_{0}^{t}\left\|\partial^{\alpha_0+1}_t\rho\right\|^2_{H^1} d \tau=& \int_{0}^{t}\left\|\partial^{\alpha_0}_t\nabla\cdot(\rho u)\right\|^2_{H^1} d \tau \\
    \lesssim& \Upsilon(0)+(\lambda_0+\eps)\Upsilon(t) + \int_{0}^{t} \|\partial^{\alpha_0}_t \nabla \nabla\cdot u \|_2^2 d \tau \\
    \lesssim& \delta_0 \eps \int_{0}^{t}\|\partial^{\alpha_0}_t \nabla^2 u \|_2^2 d \tau  + C_{\delta_0}(\Upsilon(0) + (\lambda_0 +\eps) \Upsilon(t)),
  \end{aligned}
  \end{equation}
where we have applied \eqref{H1*}. Combining the above inequalities with \eqref{H1} and \eqref{omegaL2}, and taking $\delta_0$ sufficiently small, we obtain
\begin{equation}\label{H1**}
\sup_{\tau \in [0,t]} \|\partial^{\alpha_0}_t\nabla u(\tau)\|_2^2+ \varepsilon\int_{0}^{t}\|\nabla^2 \partial^{\alpha_0}_t u\|^2_2 d\tau\lesssim \Upsilon(0)+(\lambda_0+\varepsilon)\Upsilon(t).
\end{equation}
Again, we used Proposition \ref{Pro:div}. Thereupon, \eqref{H1} follows from \eqref{H1*} and \eqref{H1**}.  
\end{proof}

To obtain the dissipation of $\partial^{\alpha_0}_t \nabla^2\rho$ for $\alpha_0\leq 2$, we apply $P_1$ to $\eqref{cocns}_2$ to get:
\begin{equation}\label{P1u}
  \begin{array}{l}
  P_{1}\left\{\rho \left(\partial^{\alpha_0+1}_{t} u+u \cdot \nabla \partial^{\alpha_0}_{t}u\right) \right\}+\partial^{\alpha_0}_{t}\nabla\left(\rho \theta  \right) \\
  \quad= \frac{4 \varepsilon}{3} P_{1}\left\{\left(\mu\left(\theta  \right)-\mu(1)\right) \partial^{\alpha_0}_{t}\nabla \nabla \cdot u\right\}+\frac{4 \varepsilon}{3} \mu(1) \partial^{\alpha_0}_{t}\nabla \nabla \cdot u \\
  \qquad
  - \varepsilon P_{1}\left\{\left(\mu\left(\theta  \right)-\mu(1)\right)\partial^{\alpha_0}_{t}   \nabla \times \omega\right\}+\sum_{i=1}^{3}\varepsilon P_1C^{\alpha}_i.
  \end{array}
\end{equation}
Then we have the following Lemma.
\begin{lemma}\label{Lm:rhoH2} Under the $a$ $prior$ assumption \eqref{assum}, it holds that
\begin{equation}\label{rhoH2}
  \sup_{\tau \in[0,t]} \sum_{\alpha_0 \leq 2}\varepsilon^2\|\partial^{\alpha_0}_t\nabla^2\rho(\tau)\|^2_2+ \sum_{\alpha_0 \leq 2} \varepsilon\int_{0}^{t}\|\partial^{\alpha_0}_t\nabla^2 \rho\|_2^2d\tau
  \lesssim \Upsilon(0)+(\lambda_0+\eps)\Upsilon(t).
\end{equation}
\end{lemma}
\begin{proof}  $\frac{4 \mu(1) \varepsilon^{2}}{3}\left(\nabla^{2}\eqref{cocns}_{1}, \partial^{\alpha_0}_{t}\nabla^{2} \rho\right)+\varepsilon\left(\nabla\eqref{P1u},  \partial^{\alpha_0}_{t}\nabla^{2} \rho\right)$  yields 
\begin{equation}\label{nabla2rho}
  \begin{aligned}
  \frac{2 \mu(1) \varepsilon^{2}}{3} & \frac{d}{d t}\left\|\partial^{\alpha_0}_t\nabla^{2} \rho\right\|_{2}^{2}+\frac{4 \mu(1) \varepsilon^{2}}{3}\left(\nabla^{2}\partial^{\alpha_0}_t\left(\nabla \rho \cdot u\right), \partial^{\alpha_0}_t\nabla^{2} \rho\right)\\
  &+\frac{4 \mu(1) \varepsilon^{2}}{3}\left(\partial^{\alpha_0}_t\nabla^{2}\left(\bar{\rho} \nabla \cdot u\right), \partial^{\alpha_0}_t\nabla^{2} \rho\right) \\
  & +\varepsilon\left(\nabla P_{1}\left\{\rho\left(\partial^{\alpha_0+1}_{t} u+u \cdot \nabla \partial^{\alpha_0}_t u\right)\right\}, \partial^{\alpha_0}_t\nabla^{2} \rho\right)+\varepsilon\left(\partial^{\alpha_0}_t\nabla^2(\rho\theta), \partial^{\alpha_0}_t\nabla^{2} \rho\right)\\
  &-\frac{4 \varepsilon^{2}}{3}\left(P_{1}\left\{\left(\mu\left(\theta  \right)-\mu(1)\right) \partial^{\alpha_0}_t \nabla \nabla \cdot u\right\}, \partial^{\alpha_0}_t\nabla^{2} \rho\right)  -\varepsilon\left(P_1C^\alpha_i, \partial^{\alpha_0}_t\nabla^{2} \rho\right) \\
  & +\frac{4 \varepsilon^{2}}{3}\left(\nabla P_{1}\left\{\left(\mu\left(\theta  \right)-\mu(1)\right) \partial^{\alpha_0}_t\nabla \times \nabla \times u\right\}, \partial^{\alpha_0}_t\nabla^{2} \rho\right) =0.
  \end{aligned}
\end{equation}
We need to be more careful for the second term on the left-hand side of \eqref{nabla2rho}, since we do not have the dissipation of $\nabla^3\rho$. Note that
\begin{equation}
  \begin{aligned}
    &\varepsilon^2\int_{0}^{t}\int\partial^{\alpha_0}_t\nabla^2(\nabla \rho \cdot u)\partial^{\alpha_0}_t\nabla^2\rho dxd\tau\\
    =&\underbrace{\varepsilon^2\int_{0}^{t}\int\partial^{\alpha_0}_t\nabla^3\rho\cdot u\partial^{\alpha_0}_t\nabla^2\rho dx d\tau}_{J_1}
    +\underbrace{\sum_{\substack{\gamma_0 \geq 1 \\
    \beta_0+\gamma_0=\alpha_0}}\varepsilon^2\int_{0}^{t}\int\partial^{\beta_0}_t\nabla^3\rho \cdot \partial_t^{\gamma_0}u \partial^{\alpha_0}_t\nabla^2\rho dxd\tau}_{J_2}\\
    &+\underbrace{\varepsilon^2\int_{0}^{t}\int\partial^{\alpha_0}_t(\nabla^2\rho\nabla u+\nabla\rho\cdot \nabla^2 u)\partial^{\alpha_0}_t\nabla^2\rho dxd\tau}_{J_3}.
  \end{aligned}
\end{equation}
For $J_1$, after integrating by parts and applying the boundary condition \eqref{slip2t}, we obtain
$$
| J_1 |=\varepsilon^2 \int_0^\mathrm{T} \int \frac{\nabla\left|\partial_t^{\alpha_0} \nabla^2 \rho\right|^2}{2} \cdot u d x d \tau=-\varepsilon^2 \int_0^\mathrm{T} \int \frac{1}{2} \nabla \cdot u\left|\partial_t^{\alpha_0} \nabla^2 \rho\right|^2 d x d \tau \lesssim \lambda_0 \varepsilon \Upsilon(t).
$$
Recall the definition of the conormal derivatives given by \eqref{def_Z1Z2Z3}, and the definition of $\varphi(z)$ specified in \eqref{def_varphi}. Then, by H{\"o}lder's inequality, the term $J_2$ can be estimated as
$$\begin{aligned}
| J_2 | =&\sum_{\substack{\gamma_0 \geq 1 \\
\beta_0+\gamma_0=\alpha_0}} \varepsilon^2 \int_0^\mathrm{T} \int \partial_t^{\beta_0} \nabla^2 \partial_{y^i} \rho \partial_t^{\gamma_0} u_i \partial_t^{\alpha_0} \nabla^2 \rho d x d \tau \\
& +\sum_{\substack{\gamma_0 \geqslant 1 \\
\beta_0+\gamma_0=\alpha_0}} \varepsilon^2 \int_0^\mathrm{T} \int \partial_t^{\beta_0}\left(\varphi(z) \nabla^2 \partial_z \rho\right)\left(\frac{\partial_t^{\gamma_0} u \cdot N}{\varphi(z)}\right) \partial_t^{\alpha_0} \nabla^2 \rho d x d \tau \\
& \lesssim \sum_{\substack{\gamma_0 \geq 1 \\\beta_0+\gamma_0=\alpha_0}}\eps^2 \int_{0}^{t}\|\p^{\beta_0}_t \nabla^2 \rho \|_{H^1_{co}}\|\p^{\gamma_0}_t u_i \|_\infty \|\p^{\alpha_0}_t \nabla^2 \rho\|_2 d\tau\\
 &\quad +\sum_{\substack{\gamma_0 \geq 1 \\
\beta_0+\gamma_0=\alpha_0}} \varepsilon^2 \int_0^\mathrm{T}\left\|\frac{\partial_t^{\gamma_0} u \cdot N}{\varphi(z)}\right\|_{\infty}\left\|\partial_t^{\beta_0} \nabla^2\rho\right\|_{H^1_{co}}\left\|\partial_t^{\alpha_0} \nabla^2 \rho\right\|_2 d \tau \\
& \lesssim \lambda_0 \varepsilon \Upsilon(t)+\sum_{\substack{\gamma_0 \geq 1 \\
\beta_0+\gamma_0=\alpha_0}} \varepsilon^2 \int_0^\mathrm{T}\left\|\partial_t^{\gamma_0} \nabla^3 u\right\|_2^{1 / 2}\left\|\partial_t^{\gamma_0} \nabla^2 u\right\|_2^{1/2}\left\|\partial_t^{\beta_0} \nabla^2 \rho\right\|_{H_{co}^{1}}\left\|\partial_t^{\alpha_0} \nabla^2 \rho\right\|_2 d \tau \\
& \lesssim \lambda_0 \Upsilon(t),
\end{aligned}
$$
where $N$ is defined in \eqref{N} and we used the Agmon's inequality \eqref{agmon}, as well as Hardy's inequality:
\begin{equation}\label{hardy}
  \left\|\frac{u \cdot N}{\varphi(z)}\right\|_\infty \lesssim\|u\|_{W^{1, \infty}}.
\end{equation}
This inequality has been rigorously established on pp.543 of \cite{masmoudi2012uniform}. Similarly, we have
\begin{equation}
 | J_3 |\lesssim \lambda_0\Upsilon(t).
\end{equation}
For the remaining terms in \eqref{nabla2rho}, standard estimates under the assumption \eqref{apriori} suffice to establish control. While the calculations are routine, they involve multiple steps that we omit here for brevity. This completes the proof of \eqref{rhoH2}.
\end{proof}
The final lemma of this subsection establishes second-order spatial derivative estimates. This analysis primarily exploits the elliptic structure inherent in the compressible Navier-Stokes system, where higher-order spatial derivatives can be bounded in terms of lower-order spatial derivatives and higher-order time derivatives. The precise formulation is presented in the following lemma.
\begin{lemma}\label{Lm:rhouthetaH3}  Under the $a$ $prior$ assumption \eqref{assum}, it holds that
\begin{equation}\label{rhouthetaH3}
  \sup _{\tau\in[0,t]} \sum_{\alpha_{0} \leq 2} \varepsilon^{2}\left\|\partial_{t}^{\alpha_{0}} \nabla^{2} (\rho,u,\theta)(\tau)\right\|_{2}^{2}+\varepsilon^3\sum_{\alpha_0\leq 2}\int_{0}^{t}\|\partial^{\alpha_0}_t\nabla^3 (u,\theta)(\tau)\|_2^2d\tau \lesssim \Upsilon(0)+(\lambda_0+\eps)\Upsilon(t).
\end{equation}
\end{lemma}
\begin{proof} 
We rewrite the compressible Navier-Stokes system in the following form:
\begin{multline*}
      -\varepsilon \mu(1) \Delta u-\frac{\varepsilon}{3} \mu(1) \nabla \nabla \cdot u=-\rho \theta_{t} u-\rho u \cdot \nabla u-\nabla p+\varepsilon\left(\mu\left(\theta \right)-\mu(1)\right) \Delta u \\
      +\varepsilon \nabla \mu\left(\theta \right) \cdot \sigma(u)+\frac{\varepsilon}{3}\left(\mu\left(\theta \right)-\mu(1)\right) \nabla \nabla \cdot u,
\end{multline*}
    \begin{multline*}
      -\varepsilon \kappa(1) \Delta \theta =-\frac{3}{2} \rho\left(\partial_{t} \theta +u \cdot \nabla \theta \right)-\rho \theta  \nabla-u+\varepsilon \nabla \kappa\left(\theta \right) \cdot \nabla \theta  \\
      +\varepsilon\left(\kappa\left(\theta \right)-\kappa(1)\right) \Delta \theta +\frac{\varepsilon}{2} \mu\left(\theta \right) \sigma(u): \sigma(u).
    \end{multline*}
By applying the standard Agmon-Douglis-Nirenberg elliptic estimates \cite{agmon1964estimates} (with verification of the Complementing Condition provided in \cite{chen2023macroscopic}), we conclude that for all $\alpha_{0} \leq2$, the following estimate holds:
\begin{equation}\label{nabla2u1}
      \varepsilon^{2}\left\|\nabla^{2} \partial_{t}^{\alpha_{0}} u\right\|_{2}^{2}
      \lesssim \left\|  \partial_{t}^{\alpha_0+1} u\right\|_{2}^{2}+\sum_{\beta_{0}  \leq \alpha_{0}}\left\|\nabla \partial_{t}^{\beta_0}\left(\bar{\rho}, \bar{\theta} \right)\right\|_{2}^{2}
      +\lambda_0^{2} \varepsilon \sum_{\beta_0 \leq \alpha_{0}}\left\|\nabla \partial_{t}^{\beta_0} u\right\|_2^{2}  
      +\lambda_0^{2}\eps^{2} \sum_{\alpha_{0} \leq \beta_{0}}\left\|\nabla^{2} \partial_{t}^{\beta_0} u\right\|_{2}^{2},
\end{equation}
\begin{multline}\label{nabla3u}
      \varepsilon^{4}\left\|\nabla^{3} \partial_{t}^{\alpha_0} u\right\|_{2}^{2}
      \lesssim \varepsilon^{2}\left\| \partial_{t}^{\alpha_{0}+1}u\right\|_{H^{1}}^{2}+\varepsilon^{2} \sum_{\beta_{0}  \leq  {\alpha}_{0}}\left\|\nabla \partial_{t}^{\beta_0} (\bar{\rho}, \bar{\theta} )\right\|_{H^{1}}^{2} \\
      +\lambda_0^{2} \varepsilon^{3} \sum_{\beta_0\leq \alpha_0}\left\|\nabla \partial_{t}^{\beta_{0} } u\right\|_{H^{1}}^{2}
      +\lambda_0^{2} \varepsilon^{4} \sum_{\beta_0\leq \alpha_0}\left\|\nabla^{3} \partial_{t}^{\beta_{0}} u\right\|_{2}^{2},
\end{multline}
and
\begin{multline}\label{nabla2theta1}
      \varepsilon^{2}\left\|\nabla^{2} \partial_{t}^{\alpha_{0}} \theta \right\|_{2}^{2}
       \lesssim\left\|\left( \partial_{t}^{\alpha_{0}+1} \theta, \nabla \partial_{t}^{\alpha_{0}} u\right) \right\|_{2}^{2} +\lambda_0^{2} \varepsilon^{2} \sum_{\alpha_{0} \leq \alpha_{0}}\left\|  \partial_{t}^{\alpha_{0}} (u, \theta )\right\|_{H^1}^{2}+\lambda_0^{2} \eps^{3} \sum_{\beta_{0} \leq \alpha_{0}}\left\|\partial_{t}^{\beta_{0} } \nabla^{2} \theta \right\|_{2}^{2},
\end{multline}
\begin{multline}\label{nabla3theta}
      \varepsilon^{4}\left\|\nabla^{3} \partial_{t}^{\alpha_{0}} \theta \right\|_{2}^{2} \leq \varepsilon^{2}\left\|\left(\partial_{t} \partial_{t}^{\alpha_{0}} \theta, \nabla \partial_{t}^{\alpha_{0}} u\right)\right\|_{H^{1}}^{2}
      +\eps^2\|\partial^{\alpha_0}_t\nabla(\rho,u,\theta)\|^2_2\\
      \quad+\lambda_0^{2} \varepsilon^{4} \sum_{\beta_0 \leq \alpha_{0}}\left\|\nabla^{2} \partial_{t}^{\beta_0} (\theta, u) (\tau)\right\|_{H^1}^{2}.
\end{multline}
Consequently, \eqref{rhouthetaH3} follows from \eqref{L2}, \eqref{H1}, \eqref{rhoH2} and the above estimates.
\end{proof}

\begin{corollary}\label{Co:classicalder}
  Under the \emph{a priori} assumption \eqref{assum}, it holds that
  \begin{equation}\label{classicalder}
    \begin{aligned}
     &\sup _{\tau\in[0,t]} \sum_{\alpha_{0} \leq 3}\left\|\partial_{t}^{\alpha_{0}}(\rho-1, u, \theta-1)(\tau)\right\|_{2}^{2}
+\sup _{\tau\in[0,t]} \sum_{\alpha_{0} \leq 2}\left\|\partial_{t}^{\alpha_{0}} \nabla(\rho, u, \theta) (\tau)\right\|_{2}^{2} \\
&\qquad +\sup _{\tau\in[0,t]} \sum_{\alpha_{0} \leq 2} \varepsilon^{2}\left\|\partial_{t}^{\alpha_{0}} \nabla^{2}(\rho,u, \theta)(\tau)\right\|_{2}^{2}+\varepsilon \sum_{\alpha_{0} \leq 3} \int_{0}^{t}\left\|\partial_{t}^{\alpha_{0}} \nabla(u, \theta)\right\|_{2}^{2}  d\tau \\
&\qquad +\varepsilon \sum_{\alpha_{0} \leq 2} \int_{0}^{t}\left\|\partial_{t}^{\alpha_{0}} \nabla \rho\right\|_{2}^{2}  d\tau+\varepsilon \sum_{\alpha_{0} \leq 2} \int_{0}^{t}\left\|\partial_{t}^{\alpha_{0}} \nabla^{2}(\rho, u, \theta)\right\|_{2}^{2}  d\tau  \\
&\qquad+\varepsilon^3\sum_{\alpha_0\leq 2}\int_{0}^{t}\|\partial^{\alpha_0}_t\nabla^3 (u,\theta)\|_2^2d\tau \lesssim \Upsilon(0)+(\lambda_0+\eps)\Upsilon(t).
    \end{aligned}
  \end{equation}
\end{corollary}
\begin{proof}
  The corollary follows directly by combining the estimates from Lemma \ref{Lm:L2}, Lemma \ref{Lm:H1}, Lemma \ref{Lm:rhoH2} and Lemma \ref{Lm:rhouthetaH3}.
\end{proof}
\subsection{Conormal Energy estimates for $\rho,u$ and $\theta$}
In this subsection, utilizing \eqref{classicalder}, we aim to establish uniform estimates on the conormal derivatives of $\rho,u$ and $\theta$, which are crucial to complete the proof of \eqref{apriori}. The analysis in this subsection follows an approach closely aligned with the conormal energy estimates developed in \cite[pp.7884-7888]{duan2021compressible}.
\begin{lemma}\label{Lm:coL2}
  Under the \emph{a priori} assumption \eqref{assum}, it holds that for $ |\alpha_1|\leq 2$,
\begin{multline}\label{coL2}
      \eps \sum_{\alpha_0\leq 2}\|\partial_t^{\alpha_0}Z^{\alpha_1}(\rho,u,\theta)(\tau)\|_2^2
      +\eps^2\sum_{\alpha_0\leq 2}\int_{0}^{t}\|\nabla\partial_t^{\alpha_0}Z^{\alpha_1}(u,\theta) (\tau)\|_2^2d\tau\\
       +\eps^2\sum_{\alpha_0\leq 1}\int_{0}^{t}\|\nabla\partial_t^{\alpha_0}Z^{\alpha_1}\rho (\tau)\|_2^2d\tau
      \lesssim \Upsilon(0)+(\lambda_0+\eps)\Upsilon(t).
\end{multline}
\end{lemma}
\begin{proof}
Multiplying $\eqref{cocns}_1,\eqref{cocns}_2$ and $\eqref{cocns}_3$ by $\mathcal{Z}^\alpha\bar{\rho},\mathcal{Z}^\alpha u,\mathcal{Z}^\alpha\bar{\theta}$ (with $\mathcal{Z}^\alpha=\partial^{\alpha_0}_tZ^{\alpha_1},\alpha_0\leq 2,|\alpha_1|=1, 2$), respectively. Then, by integrating by parts, one obtains
\begin{equation}\label{rhocoL2}
  \frac{1}{2}\frac{d}{dt}\int|\mathcal{Z}^\alpha\rho|^2dx+\int\mathcal{Z}^\alpha\nabla\cdot(\bar{\rho}u)\mathcal{Z}^\alpha \bar{\rho} d x + \int\mathcal{Z}^\alpha\nabla\cdot u\mathcal{Z}^\alpha \bar{\rho} d x=0,
\end{equation}
\begin{multline}\label{ucoL2}
      \frac{1}{2}\frac{d}{dt}\int\rho|\mathcal{Z}^\alpha u|^2dx+\int\mathcal{Z}^\alpha\nabla(\bar{\rho}\bar{\theta})\cdot\mathcal{Z}^\alpha udx
      +\int\mathcal{Z}^\alpha\nabla(\bar{\rho}+\bar{\theta})\cdot\mathcal{Z}^\alpha udx\\
      \quad =\frac{4}{3}\varepsilon\int\mu(\theta)\mathcal{Z}^\alpha\nabla\nabla\cdot u \cdot \mathcal{Z}^\alpha udx
      -\varepsilon\int\mu(\theta)\mathcal{Z}^\alpha\nabla\times\omega\cdot\mathcal{Z}^\alpha udx+\sum_{i=1}^{3}\int C^\alpha_i\mathcal{Z}^\alpha u d x,
\end{multline}
and
\begin{multline}\label{thetacoL2}
    \frac{3}{4}\frac{d}{dt}\int\rho|\mathcal{Z}^\alpha\theta|^2dx+\int(\rho\theta-1)\mathcal{Z}^\alpha\nabla\cdot u\mathcal{Z}^\alpha\theta dx+\int\mathcal{Z}^\alpha\nabla\cdot u\mathcal{Z}^\alpha\theta dx\\
    =\int\kappa(\theta)\mathcal{Z}^\alpha\Delta \theta \mathcal{Z}^\alpha\theta d x+\sum_{i=4}^{7}\int C^\alpha_i\mathcal{Z}^\alpha\theta dx.
\end{multline}
Following arguments analogous to those in Lemma \ref{Lm:L2} and employing the key estimates
\begin{equation}
  \|Zf\|_2\lesssim \|\nabla f\|_2,
\quad \text{ and }
  |[Z,\nabla]f|\lesssim |\nabla f| \quad  \text{for any } f\in H^1,
\end{equation}
where $[A,B] =AB-BA$ denotes the commutator, we derive the following estimate
\begin{equation}\label{coL2*}
\begin{aligned}
    &\varepsilon\|\mathcal{Z}^\alpha(\rho,u,\theta)(t)\|_2^2+\varepsilon^2\int_{0}^{t}\|\nabla \mathcal{Z}^\alpha(u,\theta)\|_2^2 d \tau\\
    \lesssim&\varepsilon\int_{0}^{t}\|\partial^{\alpha_0}_t\nabla(\rho,u,\theta)(\tau)\|_{H^1}^2d\tau+\Upsilon(0)+(\lambda_0+\eps)\Upsilon(t)\\
    \lesssim &\Upsilon(0)+(\lambda_0+\eps)\Upsilon(t).
  \end{aligned}
  \end{equation}
The last inequality follows from Lemma \ref{Lm:L2} and Lemma \ref{Lm:H1}.

As for the dissipation of $\nabla\mathcal{Z}^\alpha\rho$ (with $\alpha_0\leq 1,|\alpha_1 |=1,2$), we get from the inner product of $\varepsilon^2(\eqref{cocns}_1,\nabla\cdot\mathcal{Z}^\alpha u)$ and $\varepsilon^2(\eqref{cocns}_2,\nabla\mathcal{Z}^\alpha \bar{\rho}/\rho) $ that
\begin{equation*}
\begin{aligned}
      -\varepsilon^{2}&\left(\mathcal{Z}^{\alpha} u, \nabla \mathcal{Z}^{\alpha} \bar{\rho}\right)(t)+ \varepsilon^{2} \int_{0}^{t}\left\|\nabla \mathcal{Z}^{\alpha} \bar{\rho}\right\|_{2}^{2} d\tau  \\
      \leq& C \varepsilon^{2}\left|\left(\mathcal{Z}^{\alpha} u, \nabla \mathcal{Z}^{\alpha} \bar{\rho}\right)(0)\right|+C \varepsilon^{4} \int_{0}^{t}\left\|\nabla^{2} \mathcal{Z}^{\alpha} u\right\|_{2}^{2} d\tau  \\
      &+\left(\lambda_0^{2} \varepsilon+\varepsilon^{2}\right) \int_{0}^{t}\left\|\partial^{\alpha_0}_t\nabla(\bar{\rho}, u, \bar{\theta})(\tau)\right\|_{H^{1}}^{2} d\tau  \\
     &+ \lambda_0 \varepsilon^{2} \sum_{|\alpha|=2} \int_{0}^{t}\left\|\nabla \mathcal{Z}^{\alpha}(u,\theta)(\tau)\right\|_2^{2} d\tau .
\end{aligned}
\end{equation*}
We thereby conclude that
\begin{equation}\label{rhoco}
  \varepsilon^2\int_{0}^{t}\|\nabla\mathcal{Z}^\alpha\rho\|_2^2d\tau
  \lesssim \Upsilon(0)+(\lambda_0+\varepsilon)\Upsilon(t).
\end{equation}
Combining the estimates \eqref{coL2} and \eqref{rhoco}, we obtain the desired assertion \eqref{coL2}. This completes the proof of Lemma \ref{Lm:coL2}.
\end{proof}

\begin{lemma}\label{Lm:coH1}
  Under the \emph{a priori} assumption \eqref{assum}, the following estimate holds for $\alpha=(\alpha_0,\alpha_1)$ with $\alpha_0\leq 1$ and $|\alpha_1|=1,2$:
  \begin{equation}\label{coH1}
    \sup_{\tau \in [0,t]}\varepsilon^2\|\nabla\mathcal{Z}^\alpha(\rho,u,\theta)(\tau)\|_2^2+\varepsilon^3\int_{0}^{t}\|\nabla^2\mathcal{Z}^\alpha(u,\theta) (\tau)\|_2^2 d \tau \lesssim \Upsilon(0)+(\lambda_0+\eps)\Upsilon(t).
  \end{equation}
\end{lemma}
\begin{proof}
  We begin by multiplying equation $\eqref{cocns}_2$ by $\varepsilon^2\nabla\times \mathcal{Z}^\alpha\omega$, which yields
\begin{equation}
  \begin{aligned}
        &\eps^2\int_{0}^{t}\int\rho\mathcal{Z}^\alpha\partial_t u \cdot \nabla\times \mathcal{Z}^\alpha\omega dxd\tau
        +\eps^2\int_{0}^{t}\int(\rho u \cdot \nabla\mathcal{Z}^\alpha u+\mathcal{Z}^\alpha\nabla(\rho\theta))\nabla\times\mathcal{Z}^\alpha\omega dxd\tau\\
        =&\eps^3\int_{0}^{t}\mu(\theta)\mathcal{Z}^\alpha\nabla\nabla\cdot u\nabla\times\mathcal{Z}^\alpha\omega dxd\tau
        -\eps^3\int_{0}^{t}\mu(\theta) \mathcal{Z}^\alpha\nabla\times\omega  \nabla\times\mathcal{Z}^\alpha\omega dxd\tau\\
        &+\int_{0}^{t}\int C^\alpha_i\nabla\times\mathcal{Z}^\alpha\omega dxd\tau,
  \end{aligned}
\end{equation}
which further implies
\begin{equation}\label{coomega}
  \eps^2\| \mathcal{Z}^\alpha \omega\|_2^2+\eps^3\int_{0}^{t}\|\nabla\times\mathcal{Z}^\alpha\omega\|_2^2 d\tau\lesssim \Upsilon(0)+(\lambda_0+\eps)\Upsilon(t).
\end{equation}
Next, $\varepsilon^2\left(\nabla\eqref{cocns}_1,\nabla\mathcal{Z}^\alpha\rho\right)-\varepsilon^2\left(\eqref{cocns}_2,\nabla\mathcal{Z}^\alpha\nabla\cdot u\right)+\varepsilon^2\left(\nabla\eqref{cocns}_3,\nabla\mathcal{Z}^\alpha\theta\right)$ yields that
\begin{equation}
  \begin{aligned}
    \frac{\varepsilon^{2}}{2}& \frac{d}{d t}  \left\|\nabla \mathcal{Z}^{\alpha} \rho\right\|_{2}^{2}+\frac{\varepsilon^{2}}{2} \frac{d}{d t}\left\|\sqrt{\rho}\nabla \cdot \mathcal{Z}^{\alpha} u\right\|_{2}^{2}+\frac{3 \varepsilon^{2}}{4} \frac{d}{d t}\left\|\sqrt{\rho}\nabla \mathcal{Z}^{\alpha} \theta\right\|_{2}^{2} \\
    & +\varepsilon^{2}\int\nabla \mathcal{Z}^{\alpha}\left(\nabla \rho \cdot u\right)\cdot \nabla \mathcal{Z}^{\alpha} \rho dx
    +\varepsilon^{2}\int\nabla \mathcal{Z}^{\alpha}\left(\bar{\rho} \nabla \cdot u\right)\cdot \nabla \mathcal{Z}^{\alpha} \rho dx \\
    &+\varepsilon^2\int\rho\nabla\cdot\mathcal{Z}^\alpha\partial_t u[\mathcal{Z}^\alpha,\nabla\cdot]udx-\varepsilon^2\frac{1}{2}\int\partial_t\rho|\nabla\cdot\mathcal{Z}^\alpha u|^2dx\\
    &-\varepsilon^2\int\rho\mathcal{Z}^\alpha\partial_t u\cdot n\mathcal{Z}^\alpha \nabla\cdot u dS_x+\varepsilon^2\int\nabla\rho\cdot \mathcal{Z}^\alpha\partial_t u\mathcal{Z}^\alpha\nabla\cdot udx\\
    & -\varepsilon^2 \int\left(\rho u \cdot \nabla \mathcal{Z}^\alpha u+\mathcal{Z}^\alpha \nabla(\bar{\rho} \bar{\theta})\right) \nabla \mathcal{Z}^\alpha \nabla \cdot u d x. \\
    & +\varepsilon^3 \int \frac{4}{3} \mu(\theta) \mathcal{Z}^\alpha \nabla \nabla \cdot u \cdot \nabla \mathcal{Z}^\alpha \nabla \cdot u d x+\varepsilon^3 \int \mu(\theta) \mathcal{Z}^\alpha \nabla \times \omega \nabla \mathcal{Z}^\alpha \nabla \cdot u d x\\
      & +\varepsilon^2 \sum_{i=1}^3\int C_i^\alpha \nabla \mathcal{Z}^\alpha \nabla \cdot u d x-\frac{3}{4} \int \partial_t \rho\left|\nabla \mathcal{Z}^\alpha \theta\right|^2 d x+\frac{3}{2} \int \nabla \rho \cdot \mathcal{Z}^\alpha \partial_t \theta \nabla \mathcal{Z}^\alpha \theta d x \\
      & +\frac{3}{2}\varepsilon^2 \int \nabla\left(\rho u \nabla \mathcal{Z}^\alpha \theta\right) \nabla \mathcal{Z}^\alpha \theta d x+\varepsilon^2 \int \nabla\left[(\rho \theta-1) \mathcal{Z}^\alpha \nabla \cdot u\right] \nabla \mathcal{Z}^\alpha \theta d x \\
      & +\varepsilon^2 \int \nabla \mathcal{Z}^\alpha \nabla \cdot u\left[\nabla, \mathcal{Z}^\alpha\right](\rho+\theta) d x-\varepsilon^3 \int \nabla\left(k(\theta) \mathcal{Z}^\alpha \Delta \theta\right)\cdot \nabla \mathcal{Z}^\alpha \theta d x \\
      & -\varepsilon^2 \int \sum_{i=4}^7 \nabla C_i^\alpha \nabla \mathcal{Z}^\alpha \theta d x =0,
      \end{aligned}
\end{equation}
which further implies
\begin{multline}\label{coH1*}
      \varepsilon^{2}\left\|\nabla \mathcal{Z}^{\alpha}[\rho, u, \theta](t)\right\|_{2}^{2}+ \varepsilon^{3} \int_{0}^{t}\left\|\nabla \nabla \cdot \mathcal{Z}^{\alpha} u(\tau)\right\|_{2}^{2} d \tau+  \varepsilon^{3} \int_{0}^{t}\left\|\nabla^{2} \mathcal{Z}^{\alpha} \theta(\tau)\right\|_{2}^{2} d \tau \\
      \lesssim  \Upsilon(0)+ \left(\varepsilon+\lambda_0\right) \Upsilon(t).
\end{multline}
Thereupon \eqref{coH1} immediately follows from \eqref{coomega}, \eqref{coH1*} and Proposition \ref{Pro:div}. The proof of Lemma \ref{Lm:coH1} is thereby completed.
\end{proof}
The next Lemma aims to derive the estimates for $\mathcal{Z}^\alpha \nabla^2 \rho$:
\begin{lemma}\label{Lm:rhocoH2}
  Under the \emph{a priori} assumption \eqref{assum}, the following estimate holds for $\alpha=(\alpha_0,\alpha_1)$ with $\alpha_0\leq 1,|\alpha_1|=1$ or $\alpha_0=0,|\alpha_1|=2$:
  \begin{equation}\label{rhocoH2}
    \varepsilon^4\|\mathcal{Z}^\alpha\nabla^2\rho \|_2^2+\varepsilon^3\int_{0}^{t}\|\mathcal{Z}^\alpha\nabla^2\rho \|_2^2d\tau\lesssim \Upsilon(0)+(\lambda_0+\eps)\Upsilon(t).
  \end{equation} 
\end{lemma}
\begin{proof}Proceeding analogously to the derivation of \eqref{rhoH2}, we obtain
  \begin{equation}\label{rhocoH2*}
    \begin{aligned}
      \frac{2 \mu(1) \varepsilon^{4}}{3} & \left\| \mathcal{Z}^\alpha  \nabla^{2} \rho\right\|_{2}^{2}
      +\frac{4 \mu(1) \varepsilon^{4}}{3}\int_{0}^{t}\int  \mathcal{Z}^\alpha  \nabla^{2}\left(\nabla \rho \cdot u\right)   \mathcal{Z}^\alpha  \nabla^{2} \rho dxd\tau \\
      & +\frac{4 \mu(1) \varepsilon^{4}}{3}\int_{0}^{t}\int  \mathcal{Z}^\alpha  \nabla^{2}\left(\bar{\rho} \nabla \cdot u\right)   \mathcal{Z}^\alpha  \nabla^{2} \rho dxd\tau \\
      & +\varepsilon^{3}\int_{0}^{t}\int  \mathcal{Z}^\alpha  \nabla P_{1}\left\{\rho\left(\partial_{t} u+u \cdot \nabla u\right)\right\}   \mathcal{Z}^\alpha  \nabla^{2} \rho dxd\tau \\
      & +\varepsilon^{3}\int_{0}^{t}\int  \mathcal{Z}^\alpha  \nabla^{2} (\rho\theta)   \mathcal{Z}^\alpha  \nabla^{2} \rho dxd\tau+\varepsilon^{3} \sum_{\left|\beta\right| \geq 1,\beta+\gamma=\alpha} C_{\alpha,\beta} \int_{0}^{t}\int \mathcal{Z}^{\beta} \theta \mathcal{Z}^{\gamma} \nabla^{2} \rho   \mathcal{Z}^\alpha  \nabla^{2} \rho dxd\tau \\
      & -\frac{4 \varepsilon^{4}}{3}\int_{0}^{t}\int  \mathcal{Z}^\alpha  \nabla P_{1}\left\{(\mu(\theta)-\mu(1)) \nabla \nabla \cdot u\right\}   \mathcal{Z}^\alpha  \nabla^{2} \rho dxd\tau \\
      & -\varepsilon^{4}\int_{0}^{t}\int  \mathcal{Z}^\alpha  \nabla P_{1}\left\{\nabla \mu(\theta) \cdot \sigma(u)\right\}   \mathcal{Z}^\alpha  \nabla^{2} \rho dxd\tau\\
      & +\frac{4 \varepsilon^{4}}{3}\int_{0}^{t}\int  \mathcal{Z}^\alpha  \nabla P_{1}\left\{(\mu(\theta)-\mu(1)) \nabla \times \nabla \times u\right\}   \mathcal{Z}^\alpha  \nabla^{2} \rho dxd\tau=0.
      \end{aligned}
  \end{equation}
  For the second term on the left-hand side of \eqref{rhocoH2*}, one needs to be more careful, since we do not have any dissipation of $\mathcal{Z}^\alpha\nabla^3\rho$. First, one has
  \begin{equation}\label{rhocoH21}
    \frac{4 \mu(1) \varepsilon^{4}}{3}\int_{0}^{t}\int  \mathcal{Z}^\alpha  \nabla^{2}\left(\nabla \rho \cdot u\right)   \mathcal{Z}^\alpha  \nabla^{2} \rho dxd\tau
    \lesssim \eps^4\int_{0}^{t}\int \mathcal{Z}^\alpha (\nabla^3 \rho\cdot u)\mathcal{Z}^\alpha\nabla^2\rho dxd\tau+\lambda_0\Upsilon(t).
  \end{equation}
  Note that $u \cdot \nabla = u_1 \p_{y_1} + u_2 \p_{y_2} +  u \cdot N  \p_z$, where $N$ was introduced in \eqref{N}. We have 
\begin{equation}\label{rhocoH22}
      \begin{aligned}
        \mathcal{Z}^\alpha(\nabla^3\rho\cdot u) =& \mathcal{Z}^\alpha \left(\sum_{i=1,2}\p_{y_i} \nabla^2 \rho u_i + u \cdot N \p_z \nabla^2 \rho\right)\\
        = &\mathcal{Z}^\alpha \left(\sum_{i=1,2}\p_{y_i} \nabla^2 \rho\right) u_i + \mathcal{Z}^\alpha \p_{z} \nabla^2 \rho (u\cdot N)+\sum_{|\beta|\geq1,\beta+\gamma=\alpha} \sum_{i=1}^{2}\mathcal{Z}^\beta u_i\mathcal{Z}^\gamma\partial_{y_i}\nabla^2\rho\\
        & +\sum_{|\beta|\geq1,\beta+\gamma=\alpha}\mathcal{Z}^\beta(u\cdot N)\mathcal{Z}^\gamma\partial_{z}\nabla^2\rho, \\
        = & \underbrace{\sum_{i=1,2}\p_{y_i}(\mathcal{Z}^\alpha \nabla^2 \rho ) u_i + \p_z (\mathcal{Z}^\alpha \nabla^2 \rho) u \cdot N +\sum_{i=1,2}[\mathcal{Z}^\alpha, \p_{y_i}] \nabla^2 \rho u_i + [\mathcal{Z}^\alpha,\p_z ]\nabla^2 \rho    (u\cdot N)}_{\mathcal{J}_1}\\
        &+\underbrace{\sum_{|\beta|\geq1,\beta+\gamma=\alpha} \sum_{i=1}^{2}\mathcal{Z}^\beta u_i\mathcal{Z}^\gamma\partial_{y_i}\nabla^2\rho}_{\mathcal{J}_2}
        +\underbrace{\sum_{|\beta|\geq1,\beta+\gamma=\alpha}\mathcal{Z}^\beta(u\cdot N)\mathcal{Z}^\gamma\partial_{z}\nabla^2\rho}_{\mathcal{J}_3}.
      \end{aligned}
\end{equation}
 Then, from the \emph{a priori} assumption \eqref{assum} and integrating by parts, it holds that
\begin{equation}\label{rhocoH23}
        \begin{aligned}
          &\varepsilon^4\int_{0}^{t}\int \mathcal{J}_1 \mathcal{Z}^\alpha\nabla^2\rho dxd\tau\\
          = &\varepsilon^4\int_{0}^{t}\int -\frac{|\nabla^2\mathcal{Z}^\alpha\rho|^2}{2}\nabla\cdot u + \left(u_i + \frac{u \cdot N}{\varphi(z)}  \right) \mathcal{Z}^\alpha \nabla^2 \rho \mathcal{Z}^\alpha \nabla^2 \rho dx d\tau \\
          \lesssim & \lambda_0 \Upsilon(t),
        \end{aligned}
\end{equation}
where we used the Hardy's inequality \eqref{hardy}. For $\mathcal{J}_3$, one notices that
\begin{equation}
  \mathcal{Z}^\beta(u\cdot N)\mathcal{Z}^\gamma\partial_{z}\nabla^2\rho
  =\sum_{\bar{\beta}\leq \beta}C_{\bar{\beta}}(z)  \mathcal{Z}^{\bar{\beta}}\left(\frac{u\cdot N}{\varphi(z)}\right)\cdot \varphi(z)(\partial_z \mathcal{Z}^\gamma\nabla^2\rho+[\mathcal{Z}^\gamma,\partial_z]\nabla^2\rho),
\end{equation}
where $C_{\bar{\beta}}(z)$ is a bounded smooth function of $z$. Therefore, for $\alpha=(\alpha_0,\alpha_1)$ with $\alpha_0=0,|\alpha_1|=2$, Hardy's inequality yields the estimate
  \begin{equation}\label{rhocoH24}
    \begin{aligned}
      \eps^4 \int_{0}^{t} \int \mathcal{J}_3 \mathcal{Z}^\alpha \nabla^2 \rho dx d\tau &=  \varepsilon^4\int_{0}^{t}\int\sum_{|\beta|\geq1,\beta+\gamma=\alpha}\mathcal{Z}^\beta(u\cdot N)\mathcal{Z}^\gamma\partial_{z}\nabla^2\rho\mathcal{Z}^\alpha\nabla^2\rho dxd\tau\\
      &\lesssim  \varepsilon^4 \int_{0}^{t}\sum_{|\beta|\leq|\alpha|}\|\mathcal{Z}^\beta\nabla u\|_\infty\|\nabla^2\rho\|^2_{H^2_{co}} d\tau\\
      & \lesssim \eps^4 \sup_{\tau \in [0,t]} \|\nabla^2 \rho \|_{H^2_{co}} \int_{0}^{t}\|Z^2 \nabla u\|^{\frac{1}{2}}_{H^2} \|Z^2 \nabla u\|^{\frac{1}{2}}_{H^1} \| \nabla^2\rho\|_{H^2_{co}} d\tau\\
      &\lesssim \Upsilon(0)+(\lambda_0+\eps)\Upsilon(t),
  \end{aligned}
  \end{equation}
where we used Agmon's inequality \eqref{agmon}, H{\"o}lder's inequality, as well as the \emph{a prior} assumption \eqref{assum}. The situation with $\alpha_0\leq 1, \alpha_1 =1$ is analogous but simpler, so we omit its detailed treatment here. Similarly, we have 
\begin{equation}\label{rhocoH25}
  \eps^4 \int_{0}^{t} \int \mathcal{J}_2 \mathcal{Z}^\alpha \nabla^2 \rho dx d\tau \lesssim \Upsilon(0)+(\lambda_0+\eps)\Upsilon(t).
\end{equation}

Consequently, combining \eqref{rhocoH21}, \eqref{rhocoH22}, \eqref{rhocoH23}, \eqref{rhocoH24}, and \eqref{rhocoH25} one obtains that
\begin{equation}
  \int_{0}^{t}\int\mathcal{Z}^\alpha\nabla^2(\nabla\rho\cdot u)\mathcal{Z}^\alpha\nabla^2\rho dxd\tau\lesssim \Upsilon(0)+(\lambda_0+\varepsilon)\Upsilon(t).
\end{equation}
Together with \eqref{rhocoH2*}, we further have
  \begin{equation}
    \varepsilon^4\|\mathcal{Z}^\alpha\nabla^2\rho \|^2_2+\varepsilon^3\int_{0}^{t}\|\mathcal{Z}^\alpha\nabla^2\rho \|^2_2 d\tau\lesssim \Upsilon(0)+(\lambda_0+\eps)\Upsilon(t),
  \end{equation}
  for $\alpha=(\alpha_0,\alpha_1)$ with $\alpha_0\leq 1,|\alpha_1|=1$ or $\alpha_0=0,|\alpha_1|=2$.  The proof of Lemma \ref{Lm:rhocoH2} is completed.
\end{proof}

Finally, applying elliptic estimates similar to those used in the proof of Lemma \ref{Lm:rhouthetaH3}, we obtain
\begin{equation}\label{coellip}
\begin{aligned}
      \varepsilon^{4} \sum_{\substack{\alpha_0\leq 1\\|\alpha| \leq 2}}\left\|\partial_{t}^{\alpha_{0}}\nabla^{2}  Z^{\alpha} u\right\|_{2}^{2}
      &\lesssim  \varepsilon^{2} \sum_{\substack{\alpha_0\leq 1\\|\alpha| \leq 2}}\left\| \partial_{t}^{\alpha_{0}+1} Z^{\alpha} u \right\|_2^{2}
      +\eps^2\sum_{\substack{\alpha_0\leq 1\\|\alpha| \leq 2}}\left\|\partial_{t}^{\alpha_{0}}\nabla Z^\alpha(\rho,\theta) \right\|^2_2
      +l.o.t.,\\
      \varepsilon^{4} \sum_{\substack{\alpha_0\leq 1\\|\alpha| \leq 2}}\left\|\nabla^{2} \partial_{t}^{\alpha_{0}} Z^{\alpha} \theta\right\|_{2}^{2}
       &\lesssim \varepsilon^{2} \sum_{\substack{\alpha_0\leq 1\\|\alpha| \leq 2}}\left\| \partial_{t}^{\alpha_{0}+1} Z^{\alpha} \theta\right\|_2^2+
      \eps^2 \sum_{\substack{\alpha_0\leq 1\\|\alpha| \leq 2}}\left\| \nabla \partial_{t}^{\alpha_{0}} Z^{\alpha} u\right\|_{2}^{2} 
      +l.o.t.,\\
      \varepsilon^{4} \sum_{|\alpha| \leq 2}\left\| \nabla^{3} Z^{\alpha} u\right\|_{2}^{2} &\lesssim \varepsilon^{2} \sum_{|\alpha| \leq 2}\left\|\partial_{t}   Z^{\alpha} u\right\|_{H^{1}}^{2}+\varepsilon^{2} \sum_{|\alpha| \leq 2}\left\|\nabla   Z^{\alpha}[\bar{\rho}, \theta]\right\|_{H^{1}}^{2}
      +l.o.t.,\\
      \varepsilon^{4} \sum_{|\alpha| \leq 2}\left\|\nabla^{3}  Z^{\alpha} \theta\right\|_{2}^{2} &\lesssim \varepsilon^{2} \sum_{|\alpha| \leq 2}\left\| (\partial_{t}  Z^{\alpha} \theta, \nabla  Z^{\alpha} u)\right\|_{H^{1}}^{2}+l.o.t.,
\end{aligned}
\end{equation}
where $l.o.t.$ refers to lower order terms that are easy to control.

\begin{proof}
	[Sketch of the proof of Theorem \ref{Thmcns}] The estimate \eqref{apriori} follows from Corollary \ref{Co:classicalder}, Lemma \ref{Lm:coL2}, Lemma \ref{Lm:coH1}, Lemma \ref{Lm:rhocoH2} and \eqref{coellip}. The local existence of \eqref{cns}, \eqref{slip}, \eqref{cnsinitial} is standard and we refer to \cite{cho2006existence} and \cite[Section 4]{wang2016uniform}. The proof is completed.
\end{proof}

\section{Uniform bounds for the Remainder term}\label{Se:uniform}
This section is devoted to deriving uniform estimates for the remainder term through a strategic application of $L^2-L^6-L^\infty$ framework, which provides control over the nonlinear term $\Gamma(R,R)$. This approach was originally developed in the pioneering works of \cite{guo2010decay} and \cite{Esposito-Guo-Kim-Marra-2018}.

By substituting \eqref{chap2} into the Boltzmann equation \eqref{be} with specular reflection boundary condition \eqref{spec}, we derive the remainder equation
\begin{equation}\label{re}
    \partial_tR+v\cdot\nabla_xR+\frac{1}{\varepsilon}\mathcal{L}R=S,
\end{equation}
with the specular reflection boundary condition
\begin{equation}\label{spec1}
  \gamma_-R(t,x,v)=L\gamma_+R(t,x,v), \qquad \textrm{ on } \Sigma_-,
\end{equation}
and the initial data
\begin{equation}\label{reinitial}
  R(0,x,v)=R^{in}(x,v).
\end{equation}
The source term of \eqref{re} is defined as
\begin{equation}\label{S}
  S:=S_1+S_2+S_3+S_4,
\end{equation}
where
\begin{equation}\label{S*}
\begin{aligned}
      &S_1:=\varepsilon \Gamma(R,R),\\
      &S_2:=\Gamma(R,\frac{G}{\sqrt{\mu}})+\Gamma(\frac{G}{\sqrt{\mu}},R) +  {\eps}^{-1}\Gamma (G,G),\\
      &S_3:=\varepsilon^{-1}\Gamma(\frac{M-\mu}{\sqrt{\mu}},R)+\Gamma(R,\frac{M-\mu}{\sqrt{\mu}}),\\
      &S_4:=\varepsilon^{-1}\frac{1}{\sqrt{\mu}} \left(\partial_tG+v\cdot \nabla_x G+H\right).
\end{aligned}
\end{equation}
Recall that the bilinear form $\Gamma$ is defined in \eqref{def_gamma}. The function $H$ represents the higher order dissipation term, given by
\begin{equation}
    H[{\rho}, {u}, {\theta}]=M \xi \cdot \frac{\nabla_{x} \cdot[\mu({\theta}) \sigma({u})]}{{\rho} \sqrt{{\theta}}}
    +M\left(\frac{1}{3}|\xi|^{2}-1\right) \frac{\frac{1}{2} \mu({\theta}) \sigma({u}): \sigma({u})+\nabla_{x} \cdot\left[\kappa({\theta}) \nabla_{x} {\theta}\right]}{{\rho} {\theta}},
\end{equation}
 satisfying
 \begin{equation*}
  \partial_tM+v\cdot\nabla_x M-Q(M,G)-Q(G,M)=\varepsilon H.
 \end{equation*}
 We can directly verify that $S\in \ker^{\perp}(\mathcal{L})$. Indeed, the orthogonality properties  $\Gamma(f,g)\in \ker^{\perp}(\mathcal{L}), G \in \ker^{\perp}(\mathcal{L})$, combined with the fact that the quantities $(\rho,u,\theta)$ satisfy the full compressible Navier-Stokes system, immediately imply this result. The verification follows from a straightforward computation, which we omit here for brevity. 
\begin{lemma}\label{Lm:reL2}
  Let $R$ be the solution of \eqref{re}, \eqref{spec1} and \eqref{reinitial}. Then the following energy estimate holds
  \begin{equation}\label{reL2}
    \sum_{\alpha_0\leq 1}\frac{1}{2}\frac{d}{dt}\left\|\partial^{\alpha_0}_t R\right\|_{L^2_{x,v}}+\sum_{\alpha_0\leq 1}\frac{1}{\varepsilon}\left\|(\mathbf{I-P})\partial_t^{\alpha_0}R\right\|^2_{L^2_\nu}
    \lesssim \sum_{\alpha_0\leq 1}\varepsilon\left\|\partial^{\alpha_0}_t\nu^{-1/2}S\right\|^2_{L^2_{x,v}}.
  \end{equation}
\end{lemma}
\begin{proof}
Multiplying \eqref{re} by $R$ and integrating over $(x,v)$ yields
\begin{equation}\label{1reL2}
  \frac{1}{2}\frac{d}{dt}\| R \|_{L^2_{x,v}}^2+\frac{1}{2}\int_{\Sigma}(v\cdot n)R^2dS_x dv +\frac{1}{\varepsilon}\int_{\Omega\times \mathbb{R}^3}R\mathcal{L}R dxdv=\int_{\Omega\times\mathbb{R}^3}S Rdxdv.
\end{equation}
By the boundary condition \eqref{spec1},  the integrand $(v\cdot n)R^2$ is odd in $v$, and thus the boundary term vanishes
\begin{equation}\label{3reL2}
  \frac{1}{2}\int_{\Sigma}(v\cdot n)R^2dS_x dv=0.
\end{equation}
For the third term on the left-hand side of \eqref{1reL2}, we obtain from \eqref{hyperco} that
\begin{equation}\label{2reL2}
\int_{\Omega\times\mathbb{R}^3}R\mathcal{L} Rdxdv \gtrsim \|(\mathbf{I-P}) R \|_\nu^2,
\end{equation} 
For the source term on the right-hand side of \eqref{1reL2}, since $S\in \ker^\perp(\mathcal{L})$, by Cauchy-Schwarz's inequality,
\begin{equation}\label{4reL2}
  \int_{\Omega\times\mathbb{R}^3}S Rdxdv\leq \frac{\delta_0}{\varepsilon}\|(\mathbf{I-P})R\|_\nu^2+C_{\delta_0}\varepsilon\|\nu^{-\frac{1}{2}}S\|_{L^2_{x,v}}^2.
\end{equation}
 Thereupon, \eqref{reL2} follows from the above estimates for $\alpha_0=0$ by taking $\delta_0$ small. Similar arguments yield the desired result for $\alpha_0=1$, taking into account that $\partial_tR$ also satisfies the specular boundary condition \eqref{spec1}. The proof of Lemma \ref{Lm:reL2} is completed.
\end{proof}

\subsection{$L^2\mbox{-}L^6$ estimates for the Macroscopic part}
This subsection is dedicated to the $L^2\mbox{-}L^6$ estimates for the macroscopic part of the remainder $R$, mainly relies on the Theorem established in \cite{chen2023macroscopic}.  We define
\begin{equation}
  \mathbf{P}R=\left(a(t,x)+b(t,x)\cdot v+c(t,x)\frac{|v|^2-3}{2}\right)\sqrt{\mu}.
\end{equation} 
As required in \cite{chen2023macroscopic}, when the domain \( \Omega \) is axisymmetric (i.e., \( \operatorname{dim}(R_\Omega) \neq 0 \)), we assume 
\begin{equation}\label{initialcod2}
 \int_{\Omega\times \mathbb{R}^3} R^{in}v\cdot f(x)\sqrt{\mu}dxdv= 0, 
\end{equation}
where \( \mathcal{R}_{\Omega} = \{ f : x \in \Omega \mapsto \mathcal{M} x \in \mathbb{R}^3 , f(x) \cdot n(x) = 0 \text{ for any } x \in \partial \Omega \} \), with \( \mathcal{M} \) denoting the antisymmetric \( 3 \times 3 \) real matrices.
\begin{lemma}\label{Lm:maL2L6}
Let $R$ solve \eqref{re} with specular reflection boundary condition \eqref{spec1} and initial condition satisfies \eqref{initialcod2}. Then, we have the following estimates:
\begin{multline}\label{maL2}
      \varepsilon \sum_{\alpha_0\leq 1}\int_{0}^{t}\|\partial^{\alpha_0}_t(a,b,c)\|^2_{L^2_x}d\tau
       \lesssim \varepsilon|G(t)-G(0)| + \sum_{\alpha_0\leq 1}\|\partial^{\alpha_0}_t R^{in}\|^2_{L^2_x}\\
       +\sum_{\alpha_0\leq 1}\frac{1}{\varepsilon}\int_{0}^{t}\|(\mathbf{I-P})\partial^{\alpha_0}_tR\|^2_\nu d\tau
       +\sum_{\alpha_0\leq 1}\varepsilon\int_{0}^{t}\|\nu^{-\frac{1}{2}}\partial^{\alpha_0}_tS\|^2_{L^2_{x,v}}d\tau,
\end{multline}
where $|G(t)|\lesssim \sum_{\alpha_0\leq 1}\|\partial^{\alpha_0}_tR\|_{L^2_{x,v}}^2$, and
\begin{multline}\label{maL6}
     \varepsilon\sup_{\tau\in[0,t]}\|(a,b,c)(\tau)\|^2_{L^6_x}
     \lesssim \frac{1}{\varepsilon}\|(\mathbf{I-P})R(0)\|^2_\nu+\sum_{\alpha_0\leq 1}\frac{1}{\varepsilon}\int_{0}^{t}\|\partial^{\alpha_0}_t (\mathbf{I-P})R(\tau)\|_\nu^2d\tau\\
     +\varepsilon^3\|wR\|^2_{L^\infty_{x,v}}+\varepsilon\|\nu^{-\frac{1}{2}}S\|^2_{L^2_{x,v}}+\varepsilon\|\partial_tR\|^2_{L^2_{x,v}}.
\end{multline}
\end{lemma}
\begin{proof}
This lemma is a direct consequence of \cite{chen2023macroscopic}. Note that the authors in \cite{chen2023macroscopic} further assume that the initial mass and energy to be 0, since the the total energy and mass are conserved from the remainder equation \eqref{re}:
\begin{equation*}
   \int_{\Omega} a(t,x) d x= \frac{d}{dt}\int_{\Omega}c(t,x) dx=0.
\end{equation*}
Define $\bar{a}_0 = \frac{1}{|\Omega|} \int_{\Omega} a(0,x) d x$ and $\bar{c}_0 = \frac{1}{|\Omega|} \int_{\Omega} c(0,x) d x$. Then by \cite[Theorem 1.3]{chen2023macroscopic}, we have
\begin{multline}
    \varepsilon \sum_{\alpha_0\leq 1}\int_{0}^{t}\|\partial^{\alpha_0}_t(a- \bar{a}_0,b,c - \bar{c}_0)\|^2_{L^2_x}d\tau
    \lesssim \varepsilon|G(t)-G(0)|  
    +\sum_{\alpha_0\leq 1}\frac{1}{\varepsilon}\int_{0}^{t}\|(\mathbf{I-P})\partial^{\alpha_0}_tR\|^2_\nu d\tau\\
    +\sum_{\alpha_0\leq 1}\varepsilon\int_{0}^{t}\|\nu^{-\frac{1}{2}}\partial^{\alpha_0}_tS\|^2_{L^2_{x,v}}d\tau.
\end{multline}
Taking $\eps_0$ small, for any $0<\eps < \eps_0$, and then using the fact that 
\begin{equation*}
  \sum_{\alpha_0\leq 1}\|\partial^{\alpha_0}_t(\bar{a}_0,\bar{c}_0) \|_{L^2_x}^2 \lesssim \sum_{\alpha_0\leq 1}\|\partial^{\alpha_0}_t R^{in}\|^2_{L^2_x},
\end{equation*}
we immediately conclude \eqref{maL2}.

Now we focus on proving \eqref{maL6}. It is known from \cite[Theorem 1.2]{chen2023macroscopic} that
\begin{equation}\label{1maL6}
  \varepsilon\|(a,b,c)\|_{L^6_x}^2
  \lesssim \varepsilon\|(\mathbf{I-P})R\|^2_{L^6_{x,v}}+\frac{1}{\varepsilon}\|(\mathbf{I-P})R\|^2_{L^2_{x,v}}
  +\varepsilon\|\nu^{-\frac{1}{2}}(S-\partial_tR)\|_{L^2_{x,v}}^2.
\end{equation}
For the first term on the right-hand side of \eqref{1maL6}, we use H{\"o}lder's inequality to obtain
\begin{equation}\label{2maL6}
  \varepsilon\|(\mathbf{I-P})R\|^2_{L^6_{x,v}}\lesssim \varepsilon^3\|(\mathbf{I-P})R\|_{L^\infty_{x,v}}^2+\varepsilon^{-1}\|(\mathbf{I-P})R\|_{L^2_{x,v}}^2.
\end{equation}
Applying  Cauchy-Schwarz's  inequality yields
\begin{equation}\label{3maL6}
\begin{aligned}
     \frac{1}{\varepsilon}\|(\mathbf{I-P})R\|^2_{L^2_{x,v}}
     &\lesssim \frac{1}{\varepsilon}\|(\mathbf{I-P})R(0)\|^2_{L^2_{x,v}}+\frac{1}{\varepsilon}\int_{0}^{t}\|(\mathbf{I-P})R\|_{L^2_{x,v}}\|(\mathbf{I-P})\partial_tR\|_{L^2_{x,v}}d\tau\\
     &\lesssim \frac{1}{\varepsilon}\|(\mathbf{I-P})R(0)\|^2_{L^2_{x,v}}
     +\frac{1}{\varepsilon}\sum_{\alpha_0\leq 1}\int_{0}^{t}\|(\mathbf{I-P})\partial^{\alpha_0}_tR\|^2_{L^2_{x,v}}d\tau.
\end{aligned}
\end{equation}
Therefore, from plugging \eqref{2maL6} and \eqref{3maL6} into \eqref{1maL6}, we obtain the desired estimate.
\end{proof}
\begin{remark}
  From Theorem 1.3 in \cite{chen2023macroscopic}, we can only initially obtain $\|\partial^{\alpha_0}_tS\|^2_{L^2_{x,v}}$ instead of $\|\nu^{-\frac{1}{2}}\partial^{\alpha_0}_tS\|^2_{L^2_{x,v}}$ on the right-hand side of \eqref{maL2}. However, upon examining the details of their proof, it turns out that it is indeed possible to derive the latter expression. Consequently, equation \eqref{maL2} holds true, see \cite[pp.797-800]{chen2023macroscopic} for details.
\end{remark}
\begin{remark}
  Lemma \ref{Lm:maL2L6} corresponds to Lemma 3.2 in \cite{duan2021compressible}, where the authors establish a similar estimate under the diffuse boundary condition. Their proof relies crucially on the choice of special test functions to derive the macroscopic estimates. To the best of our knowledge, the derivation of such macroscopic estimates for the general Maxwell reflection boundary condition remains an open problem in the field.
\end{remark}

\subsection{$L^\infty$ Estimate}
In this subsection, we study the $L^\infty$ estimates for the solutions of the remainder equation \eqref{re}. As in \cite{guo2010decay} or \cite{guo2021hilbert}, we introduce the following backward characteristics. Given $(t,x,v)$, we define $[X(\tau),V(\tau)]$ such that
\begin{equation*}
  \begin{cases}
    \frac{d}{d\tau }X (\tau)=V (\tau),\quad \frac{d}{dx}V (\tau)=0,\\
    [X(t;t,x,v),V(t;t,x,v)]=[x,v].
  \end{cases}
\end{equation*}
Then $[X(\tau;t,x,v),V(\tau;t,x,v)]=[x-(t-\tau)v,v]=[X (\tau),V (\tau)]$. For $(x,v)\in \Omega\times \mathbb{R}^3$, the backward exist time $t_b(x,v)>0$ is defined as the first moment when the backward characteristics $[X(s;0,x,v);V(s;0,x,v)]$ emerge from $\partial\Omega$:
$$t_b(x,v)=\inf\{t>0:x-tv\notin \Omega\},$$
which indicates that $x-t_b v\in \partial\Omega$. We also define $x_b(x,v)=x-t_bv$. Note that we use $t_b(x,v)$ whenever it is well-defined. We always have $v\cdot n(x_b)\leq 0$. For any point $(t,x,v)\notin\gamma_0\cap\gamma_-$, we define $(t_0,x_0,v_0)=(t,x,v)$ and recursively,
\begin{equation*}
  (t_{k+1},x_{k+1},v_{k+1})=(t_k-t_b(t_k,x_k,v_k),x_b(x_k,v_k),R_{x_{k+1}}v_k),
\end{equation*}
where $R_{x_{k+1}}v_k=v_k-2\left(v_k-2(v_k\cdot n(x_{k+1}))n(x_{k+1})\right)$. The specular backward characteristics are given by
\begin{equation*}
\begin{aligned}
    &X_{cl}(\tau;t,x,v)=\mathbf{1}_{[t_{k+1},t_k)} (\tau)\{x-(t_k-\tau)v_k\},\\
    &V_{cl}(\tau;t,x,v)=\mathbf{1}_{[t_{k+1},t_k)} (\tau)v_k.
\end{aligned}
\end{equation*}
\begin{lemma}\label{Lm:Kest}\cite[Lemma 3]{guo2010decay}
Recall \eqref{K}. We have
  \begin{equation}\label{Kest1}
   |k(v, v')| \leq C\{|v-v'|+|v-v'|^{-1}\} e^{-\frac{|v-v'|^{2}}{8}} e^{-\frac{\left||v|^{2}-|v'|^{2}\right|^{2}}{8|v-v'|^{2}}}.
  \end{equation}
    Let $0\leq \vartheta <\frac{1}{4}$ and $k\geq 0$. Then for $\eta>0$ sufficiently small,
  \begin{equation}\label{Kest2}
    \int_{\mathbb{R}^3}e^{\eta|v-v'|^2}|k(v,v')|\frac{w_k(v) e^{\vartheta|v|^2}}{w_k(v') e^{\vartheta|v'|^2}}dv'\leq \frac{C}{1+|v|},
  \end{equation}
  for some $ C>0$.
\end{lemma}

We denote the weighted non-local operator
\begin{equation}\label{3.2.1}
  K_wg=wK(\frac{g}{w})=\int_{\mathbb{R}_3}k_w(v,v')g(v')dv',
\end{equation}
where
\begin{equation}
  k_w(v,v'):=k(v,v')\frac{w_k(v)}{w_k(v')}.
\end{equation}

\begin{lemma}\label{Lm:reinfty}
  Let $R$ be a solution of \eqref{re}, \eqref{spec1} and \eqref{reinitial}. Then, it holds that
 \begin{equation} \label{reinfty}
       \begin{aligned}
        \varepsilon\sup_{\tau\in[0,t]}\|R_w (\tau)\|_{L^\infty_{x,v}}
        \lesssim &\varepsilon\|R_w(0)\|_{L^\infty_{x,v}}+\sup_{\tau \in [0,t]}\varepsilon^{\frac{1}{2}}\|\mathbf{P} R  (\tau)\|_{L^6_{x,v}}+\varepsilon^{-\frac{1}{2}}\sup_{\tau\in[0,t]}\|(\mathbf{I-P})R (\tau) \|_{L^2_{x,v}}\\
        &+\varepsilon^2\sup_{\tau\in[0,t]}\|w_{k-1}S (\tau)\|_{L^\infty_{x,v}}.
       \end{aligned}
 \end{equation}
 \end{lemma}
\begin{proof}
  Based on the remainder equation \eqref{re} and the definition \eqref{wk}, we deduce that $R_w $ satisfies
  \begin{equation}\label{reweight}
    \begin{cases}
          \partial_t R_w +v\cdot\nabla_xR_w +\frac{\nu}{\varepsilon}R_w -\frac{1}{\varepsilon}K_wR_w =w_kS,\\
          \gamma_-R_w (t,x,v)=L\gamma_+R_w(t,x,v) \quad (x,v)\in \Sigma_-.
    \end{cases}
    \end{equation}
  By integrating $\eqref{reweight}_1$ along the backward characteristics, we obtain
\begin{equation}\label{reinfty1}
   \begin{aligned}
    R_w(t, x, v)=& e^{-\frac{\nu t}{\varepsilon}} R_w\left(0, x_k-t_k v_k, v_k\right)\\
    &+\sum_{j=0}^k \int_{t_{j+1}}^{t_j} e^{\frac{\nu(\tau -t)}{\varepsilon}}\left(\frac{1}{\varepsilon} K_w R_w+w_k S\right)\left(\tau, x_j-\left(t_j-\tau \right)v_j, v_j\right) d\tau,
   \end{aligned}
\end{equation}
which further implies

\begin{equation}
  \begin{aligned}
         \left|R_w(t,x, v)\right|&\leq \underbrace{e^{-\frac{\nu t}{\varepsilon}}\left|R_w\left(0, x_k-t_k v_k, v_k\right)\right|}_{I_1}\\
         &\quad+\underbrace{\sum_{j=0}^k \int_{t_{j+1}}^{t_j} e^{\frac{\nu(\tau -t)}{\varepsilon}}\left| w_k S\left(\tau, x_j-(t_j-\tau )v_j, v_j\right)\right| d\tau }_{I_2} \\
         &\quad+\underbrace{\sum_{j=0}^k \int_{t_{i+1}}^{t_j} e^{\frac{\nu(\tau -t)}{\varepsilon}} \frac{1}{\varepsilon} \left| K_w R_w(\tau, x_j-(t_j-\tau )v_j, v_j) d\tau \right|}_{I_3}.
  \end{aligned}
\end{equation}
We straightforward obtain
$$
I_1 \lesssim e^{-\frac{\nu t}{\varepsilon}}\left\|R_w(0)\right\|_{L^\infty_{x,v}}, \quad I_2 \lesssim \varepsilon \sup_{\tau \in[0, t]}\left\|w_{k-1} S(\tau)\right\|_{L^\infty_{x,v}}.
$$
Regarding the estimate for $I_3$, we denote
$\left(t_0^{\prime}, x_0^{\prime}, v_0^{\prime}\right)=\left(\tau, x_j-\left(t_j-\tau \right) v_j, v^{\prime}\right)$ and define a new back-time cycle as
$$
(t_{j+1}^{\prime} x_{j+1}^{\prime}, v_{j+1}^{\prime})=(t_j^{\prime}-t_b(x_j^{\prime}, v_j^{\prime}), x_b(x_j^{\prime}, v_j^{\prime}), R_{x_{j+1}^\prime}v_{j}^{\prime}),
$$
where $R_{x_{j+1}^\prime}v_{j}^{\prime}=v_{j}^{\prime}-2\left(v_{j}^{\prime}\cdot n(x_{j+1})\right)n(x_{j+1})$. We then use \eqref{reinfty1} again to obtain
\begin{equation}
  \begin{aligned}
   I_3&\leq \underbrace{\sum_{j=0}^{k}\int_{t_{j+1}}^{t_j}\int_{\mathbb{R}^3}\frac{1}{\varepsilon} e^{-\frac{\nu t}{\varepsilon}}\left|k_w(v_j,v')R_w\left(0, x'_{k'}-t'_{k'} v'_{k'}, v'_{k'}\right)\right|}_{I_{3,1}}d\tau dv'\\
  &\quad+\underbrace{\sum_{j=0}^k\sum_{j'=0}^{k'} \int_{t_{j+1}}^{t_j} \int_{t'_{j'+1}}^{t'_{j'}}\int_{\mathbb{R}^3} \frac{1}{\varepsilon}e^{\frac{\nu(\tau '-t)}{\varepsilon}}\left| k_w(v_j,v')w_k S\left(\tau, x'_{j'}-(t'_{j'} - \tau') v'_{j'}, v'_{j'}\right)\right| d\tau d\tau 'dv'}_{I_{3,2}} \\
  &\quad+\underbrace{\sum_{j=0}^k\sum_{j'=0}^{k'} \int_{t_{j+1}}^{t_j} \int_{t'_{j'+1}}^{t'_{j'}}\int_{\mathbb{R}^3} \frac{1}{\varepsilon^2}e^{\frac{\nu(\tau '-t)}{\varepsilon}}\left| k_w(v_j,v')K_wR_w\left(\tau, x'_{j'}-(t'_{j'} - \tau')v'_{j'}, v'_{j'}\right)\right| d\tau d\tau 'dv'}_{I_{3,3}}.
  \end{aligned}
\end{equation}
It is straightforward to see that
$$
I_{3, 1} \lesssim e^{-\frac{\nu t}{\varepsilon}}\left\|R_w(0)\right\|_{L^{\infty}_{x,v}}, \quad I_{3,2} \lesssim\varepsilon \sup_{\tau  \in[0,t]}\left\|w_{k-1} S(\tau)\right\|_{L^{\infty}_{x,v}}.
$$
Next, we compute $I_{3,3}$. We divide our computations into the following cases:

  Case 1. $|v_j|\geq N_0$  with $N_0$ suitably large. It follows from lemma \ref{Lm:Kest} that
  \begin{equation}\label{3.2.20}
    \int_{\mathbb{R}^3}|k_w(v_j,v')|dv'\lesssim \frac{1}{N_0},
  \end{equation}
  which yields that
  \begin{equation}\label{3.2.22}
    |I_{3,3}|\lesssim \frac{1}{N_0}\sup_{\tau \in[0,s] }\|R_w  (\tau)\|_{L^{\infty}_{x,v}}.
  \end{equation}

  Case 2. $|v_j|\leq N_0$ and $|v'|\geq 2N_0$ or $|v'|\leq 2N_0$ and $|v''|\geq 3N_0$. Notice that $v_j'=v'$, and thus we have either $|v_j-v'|\geq N_0$ or $|v_j'-v''|\geq N_0$ and either of the following holds:
  \begin{equation}\label{3.2.23}
   | k_w(v_1,v')|\leq C e^{-\frac{\eta N_0^2}{8}}| k_w(v_1,v')|e^{\frac{\eta|v_1-v'|^2}{8}},
  \end{equation}
  or
  \begin{equation}\label{3.2.24}
    |k_w(v_1',v'')|\leq C e^{-\frac{\eta N_0^2}{8}} |k_w(v_1',v'')|e^{\frac{\eta|v_1'-v''|^2}{8}},
  \end{equation}
  for some small $\eta$. Hence we have
  \begin{equation}\label{3.2.25}
    |I_{3,3}|\leq C_\eta e^{-\frac{\eta}{8}N_0^2}\sup_{\tau \in[0,s] }\|R_w (\tau)\|_{L^{\infty}_{x,v}}.
  \end{equation}

  Case 3a. $|v_1|\leq N_0,|v'|\leq 2N_0,|v''|\leq 3N_0$, and $\tau '\geq t_{j'}'-\varepsilon\kappa_*$,

  \begin{equation}\label{3.2.26}
    |I_{3,3}|\leq C_{N_0}\kappa_*\sup_{\tau \in[0,s] }\|R_w  (\tau)\|_{L^{\infty}_{x,v}}.
  \end{equation}

  Case 3b. $|v_1|\leq N_0,|v'|\leq 2N_0,|v''|\leq 3N_0$, and $\tau '\leq t_{j'}'-\varepsilon\kappa_*$. This is the last remaining case. We denote $D=\{|v'|\leq 2N_0,|v''| \leq 3N_0\}$. We now can bound $I_{3,3}$ by

\begin{multline}\label{I33}
      I_{3,3}\leq  \sum_{j=0}^k\sum_{j'=0}^{k'} \int_{t_{j+1}}^{t_j} \int_{t'_{j'+1}}^{t'_{j'}-\kappa_*\varepsilon}\iint_{D} \frac{1}{\varepsilon^2}e^{\frac{\nu(\tau '-t)}{\varepsilon}}\left| k_w(v_j,v')k_w(v'_{j'},v'')\right|\\
    \times \left| R_w\left(\tau, x'_{j'}-(t'_{j'} - \tau') v'_{j'}, v'_{j'}\right)\right| d\tau d\tau 'dv'dv''.
\end{multline}
  Based on the Lemma \ref{Lm:Kest}, $k_w(v_1,v'$) has a possible integrable singularity of $\frac{1}{|v_1-v'|}$. We choose a number $m(N_0)$ to define
  \begin{equation*}
    k_{m,w}(p,v')=\mathbf{1}_{|p-v'|\geq \frac{1}{m},|v'|\leq m}k_w(p,v'),
  \end{equation*}
with
  \begin{equation}\label{3.2.28}
    \sup_{|p|\leq 3N_0}\int_{|v'|\leq 3N_0}|k_{m,w}(p,v')-k_w(p,v')|dv'\leq \frac{1}{N_0}.
  \end{equation}
  We split

\begin{multline*}
      k_w(v_1,v')k_w(v_1',v'')  =\{k_w(v_1,v')-k_{m,w}(v_1,v')\}k_{w}(v_1',v'')\\
      +\{k_w(v_1',v'')-k_{m,w}(v_1',v'')\}k_{m,w}(v_1,v')+k_{m,w}(v_1,v')k_{m,w}(v',v'').
\end{multline*}
  The first two differences leads to a small contribution to $I_{3,3}$:
  \begin{equation}\label{3.2.29}
    \frac{C}{N_0}\sup_{\tau \in [0,,t]}\|R_w \|_{L^{\infty}_{x,v}}.
  \end{equation}
  For the main contribution of $k_{w,m}(v_1,v')k_{w,m}(v_1',v'')$, using \eqref{I33} and H{\"o}lder's inequality, we obtain
  $$
\begin{aligned}
I_{3,3}\leq & \sum_{j=0}^k \int_{t_{j+1}}^{t_j} e^{\frac{\nu(\tau-t)}{\varepsilon}} d\tau \left( \iint_D\sum_{j'=0}^{k'} \int_{t'_{j'+1}}^{t'_{j'}-k_{*} \varepsilon} e^{\frac{\nu\left(\tau^{\prime}-\tau \right)}{\varepsilon}} |k_{w,m}\left(v_j, v^{\prime}\right) k_{w,m}\left(v_{j^{\prime}}, v^{\prime\prime}\right)|^2 d \tau d\tau^\prime d v^{\prime} d v^{\prime\prime}\right)^{1 / 2} \\
& \times\left(\iint_D\sum_{j' = 0}^{k'} \int_{t_{j+1}^{\prime}}^{t_j^{\prime}-k_{*} \varepsilon} e^{\frac{\nu\left(\tau^{\prime}-\tau \right)}{\varepsilon}} \left|\mathbf{(I-P)}R_w\left(\tau^{\prime}, x'_{j'}-(t'_{j'} - \tau') v'_{j'}, v^{\prime \prime}\right)\right|^2 d\tau d\tau'  d v^{\prime} d v^{\prime \prime}\right)^{1 / 2}\\
&+\sum_{j=0}^k \int_{t_{j+1}}^{t_j} e^{\frac{\nu(\tau-t)}{\varepsilon}} d\tau \left( \iint_D\sum_{j=0}^k \int_{t_{j+1}}^{t_j-k_{*} \varepsilon} e^{\frac{\nu\left(\tau^{\prime}-\tau\right)}{\varepsilon}} |k_{w,m}\left(v_j, v^{\prime}\right) k_{w,m}\left(v_{j^{\prime}}, v^{\prime\prime}\right)|^{\frac{6}{5}} d\tau  d\tau ' d v^{\prime} d v^{\prime\prime}\right)^{5/6} \\
& \times\left(\iint_D\sum_{j=0}^{k} \int_{t_{j+1}^{\prime}}^{t_j^{\prime}-k_{*} \varepsilon} e^{\frac{\nu\left(\tau^{\prime}-\tau\right)}{\varepsilon}} \left|\mathbf{P}R_w\left(\tau^{\prime}, x'_{j'}-(t'_{j'} - \tau') v'_{j'}, v^{\prime \prime}\right)\right|^6 d\tau ' d v^{\prime} d v^{\prime \prime}\right)^{1 / 6}.
\end{aligned}
$$
We make a change of variable: $ v_j^{\prime} \rightarrow y=x'_{j'}-(t'_{j'} - \tau') v'_{j'}$. Notice that $x_j^{\prime}$ is independent of $v_j^{\prime}$. We see that
$$
\left|\operatorname{det}\left(\frac{\partial y}{\partial v^\prime_j}\right)\right| \geqslant\left(\varepsilon k_*\right)^3>0, \quad \text {for} \quad \tau^{\prime} \in [t_{j'+1}^\prime, t_{j'}^{\prime}-k_* \varepsilon).
$$
Consequently,
\begin{multline}
 I_{3,3} \leqslant \frac{C_{N_0}}{\left(\varepsilon k_*\right)^{1/2}} \sup_{\tau\in[0,t]} \|\mathbf{P}R(\tau)\|_{L^6_{x,v}}+\frac{C_{N_0}}{\left(\varepsilon k_*\right)^{3/2}} \sup_{\tau\in[0,t]} \|\mathbf{(I-P)}R(\tau)\|_{L^2_{x,v}}\\
  +C_{N_0} k_{*} \sup_{\tau\in[0, t]}\left\|\mathbf{(I-P)}R_w(\tau)\right\|_{L^{\infty}_{x,v}}.
\end{multline}

  Collecting all the above estimates, one has
  \begin{equation}\label{3.2.31}
  \begin{aligned}
     \sup_{\tau\in[0, t]} \|R_w  (\tau)\|_{L^{\infty}_{x,v}}
    \lesssim & \eps \|R_w (0)\|_{L^{\infty}_{x,v}}+(\frac{1}{N_0}+\kappa_*)\sup_{\tau\in[0, t]}\|R_w  (\tau)\|_{L^{\infty}_{x,v}}\\
    & +\frac{C_{N_0,k_*}}{\varepsilon^\frac{3}{2}}\sup_{\tau\in[0, t]}\|\mathbf{(I-P)}R (\tau)\|_{L^2_{x,v}}
      +\frac{C_{N_0,k_*}}{\varepsilon^\frac{1}{2}}\sup_{\tau\in[0, t]}\|\mathbf{P}R (\tau)\|_{L^6_{x,v}},
  \end{aligned}
  \end{equation}
  which concludes the Lemma \ref{Lm:reinfty} by choosing sufficiently small $ m > 0, \kappa_*> 0$ and
  large $N_0 >0$.
\end{proof}

\subsection{Remainder Estimate}
In this subsection, with Lemma \ref{Lm:reL2}, Lemma \ref{Lm:maL2L6} and Lemma \ref{Lm:reinfty}, we employ $L^2\mbox{-}L^6\mbox{-}L^\infty$ framework to establish the uniform estimates for the remainder term $R$ in \eqref{re}. We begin by stating the following foundational lemma:
\begin{lemma}\label{Lm:gamma}\cite{guo2010decay} It holds that
\begin{equation}\label{gamma}
    \|w_{k-1}\Gamma(f,g)\|_{L^\infty_{x,v}}\lesssim \|w_kf\|_{L^\infty_{x,v}}\|w_kg\|_{L^\infty_{x,v}},
\end{equation}
and
\begin{equation}\label{gamma2}
  \|\nu^{-\frac{1}{2}} \Gamma(f,g) \|_{L^2_{v}} \lesssim  \|f\|_{L^2_{v}}\|\nu^{\frac{1}{2}}g\|_{L^2_{v}} + \| g\|_{L^2_{v}}\|\nu^{\frac{1}{2}}f\|_{L^2_{v}}.
\end{equation}
\end{lemma}
\begin{proof}
  The first equality is established in \cite[Lemma 5]{guo2010decay}, while the second follows from \cite[Lemma 2.3]{Guo-2002-CPAM}.
\end{proof}
\begin{lemma}\label{Lm:reuniform} 
  Recall \eqref{Rnorm} for the definition of the norm $\|\cdot\|_{\mathbf{X}}$. Assume that the initial data satisfies
\begin{equation}
  \sum_{\alpha_0\leq 1}\left\|\partial^{\alpha_0}_tR^{in}\right\|^2_{L^2_{x,v}}+\frac{1}{\eps}\left\|(\mathbf{I-P})R^{in}\right\|^2_{L^2_{x,v}}+\varepsilon\left\|w_k R^{in} \right\|^2_{L^\infty_{x,v}}\lesssim\delta,
\end{equation}
for $\delta$ suitably small, as well as the relation \eqref{initialcod2} if the domain $\Omega$ is axisymmetric. Then there exists a unique $R$ defined on $(t,x,v)\in[0,T]\times\Omega\times \mathbb{R}^3$, where $T$ is given in Theorem \ref{Thmcns}, that
   solves \eqref{re}, \eqref{spec1} and \eqref{reinitial}. Furthermore, it holds that
  \begin{equation}\label{reuniform}
    \|R\|^2_{\mathbf{X}}\lesssim\delta.
  \end{equation}
\end{lemma}
\begin{proof}
  In view of Lemma \ref{Lm:reL2}, Lemma \ref{Lm:maL2L6} and Lemma \ref{Lm:reinfty}, we have
  \begin{equation}\label{reuniform1}
\begin{aligned}
  \|R\|_{\mathbf{X}}^2
  \lesssim &\varepsilon \sum_{\alpha_0\leq 1}\int_{0}^{t}\left\|\nu^{-\frac{1}{2}}\partial^{\alpha_0}_tS\right\|^2_{L^2_{x,v}}d\tau
  +\varepsilon^4\sup_{\tau\in[0,t]}\left\|w_{k-1}S\right\|^2_{L^\infty_{x,v}}
  +\frac{1}{\varepsilon}\left\|(\mathbf{I-P})R(0)\right\|_\nu^2\\
  &+\varepsilon\left\|\nu^{-\frac{1}{2}}S(0)\right\|^2_{L^2_{x,v}}+\sum_{\alpha_0\leq 1}\left\|\partial^{\alpha_0}_tR(0)\right\|^2_{L^2_{x,v}}+\varepsilon^2\left\|w_kR(0)\right\|^2_{L^\infty_{x,v}}.
\end{aligned}
  \end{equation}
Recall that $S$ is defined in \eqref{S}. By virtue of Lemma \ref{Lm:gamma}, it holds that for $\alpha_0\leq 1$,
  \begin{equation}
\begin{aligned}
     \varepsilon\int_{0}^{t}\left\|\nu^{-\frac{1}{2}}\partial^{\alpha_0}_tS_1\right\|^2_{L^2_{x,v}}
     &=
     \varepsilon^3 \int_{0}^{t}\left\|\nu^{-\frac{1}{2}}\Gamma(\partial^{\alpha_0}_tR,R)\right\|_{L^2_{x,v}}^2d\tau\\
    & \lesssim \varepsilon^3\int_{0}^{t}\|w_kR\|^2_{L^\infty_{x,v}}\| {\nu}^{\frac{1}{2}}\partial^{\alpha_0}_tR\|^2_{L^2_{x,v}}d\tau\\
    &\lesssim \varepsilon^2\|w_kR\|^2_{L^\infty_{x,v}}\left(\varepsilon\int_{0}^{t}\left\|\partial_t^{\alpha_0}(a,b,c)\right\|_2^2d\tau+\int_{0}^{t}\left\|(\mathbf{I-P})R\right\|_\nu^2d\tau\right)\\
    &\lesssim \|R\|^4_{\mathbf{X}} + \delta \|R \|^2_{\mathbf{X}}.
\end{aligned}
  \end{equation}
  Recall \eqref{G} for the definition of $G$. Then we have
  \begin{equation}
    \begin{aligned}
         \varepsilon\int_{0}^{t}\left\|\nu^{-\frac{1}{2}}\partial^{\alpha_0}_tS_2\right\|^2_{L^2_{x,v}}
         &\lesssim\sum_{\beta_0+\gamma_0\leq \alpha_0}
 \varepsilon\int_{0}^{t}\left\|\nu^{-\frac{1}{2}}\Gamma(\partial^{\beta_0}_tR,\frac{\partial^{\gamma_0}_tG}{\sqrt{\mu}})\right\|_{L^2_{x,v}}^2d\tau\\
        & \lesssim \sum_{\beta_0+\gamma_0\leq \alpha_0} \varepsilon\int_{0}^{t}\left\|\nu^{\frac{1}{2}}\partial^{\beta_0}_t R\right\|^2_{L^2_{x,v}}\left\|\frac{\partial^{\gamma_0}_tG}{\sqrt{\mu}} \right\|^2_{L^2_{x,v}}d\tau\\
        &\lesssim \lambda^2_0\varepsilon^3\|R\|^2_{\mathbf{X}},
    \end{aligned}
      \end{equation}
where we used $\left\|\frac{\partial^{\gamma_0}_tG}{\sqrt{\mu}} \right\|^2_{L^2_{x,v}} \lesssim \|\partial^{\gamma_0}_t \nabla(u,\theta)\|^2_{L^2_{x,v}} \lesssim  \lambda_0^2 \eps^2$ by Theorem \ref{Thmcns}. Similarly, we have the following estimates for $S_3$ and $S_4$
\begin{equation}
  \varepsilon\int_{0}^{t}\left\|\nu^{-\frac{1}{2}}\partial^{\alpha_0}_t(S_3+S_4)\right\|^2_{L^2_{x,v}}\lesssim \lambda_0^2\|R\|_{\mathbf{X}}^2.
\end{equation}
In total, we have
\begin{equation}\label{SL2}
  \sum_{\alpha_0\leq 1}\varepsilon\int_{0}^{t}\left\|\nu^{-\frac{1}{2}}\partial^{\alpha_0}_t S\right\|^2_{L^2_{x,v}}\lesssim \|R\|^4_{\mathbf{X}}+\lambda_0(\lambda_0+\varepsilon)\|R\|_{\mathbf{X}}^2.
\end{equation}

As for the $L^\infty_{x,v}$ estimates of $S$, by Lemma \ref{Lm:gamma}, one has 
\begin{equation}\label{S1infty}
  \varepsilon^4\|w_{k-1}S_1\|_{L^\infty_{x,v}}^2
  = \varepsilon^6\left\|w_{k-1}\Gamma(R,R)\right\|_{L^\infty_{x,v}}^2\lesssim\varepsilon^6\|w_kR\|_{L^\infty_{x,v}}^4\lesssim\varepsilon^2\|R\|_{\mathbf{X}}^4,
\end{equation}
and
\begin{equation}\label{Siinfty}
  \sum_{i=2}^{4}\varepsilon^4\left\|w_{k-1}S_i\right\|_{L^\infty_{x,v}}^2
\lesssim\lambda_0^2\varepsilon^2\|R\|_{\mathbf{X}}^2+\lambda_0^2\varepsilon.
\end{equation}
Collecting \eqref{reuniform1}, \eqref{SL2}, \eqref{S1infty} and \eqref{Siinfty}, we conclude
\begin{equation}
  \|R\|_{\mathbf{X}}^2\lesssim\|R\|^4_{\mathbf{X}}+\lambda_0^2\varepsilon+\delta.
\end{equation}
Therefore, the existence of \eqref{re}, \eqref{spec1}, and \eqref{reinitial}, as well as the estimate \eqref{reuniform} follow from a standard iteration argument, see \cite[Section 6]{duan2021compressible} for instance.
\end{proof}

\begin{proof}
  [Completing the proof of Theorem \ref{Tm:main}]: The estimate \eqref{main} directly follows from Lemma \ref{reuniform}. Based on the expansion
  \begin{equation}
    F=M+\eps G+\eps^2 \sqrt{\mu}R,
  \end{equation}
  the uniqueness and positivity of $F$ follows from the same arguments as \cite[Theorem 3]{guo2010decay}. We now focus on the proof of \eqref{main1}. From \eqref{reuniform}, it holds that
  \begin{equation*}
    \|R\|^2_{L^2_{x,v}}\lesssim \delta,\quad \varepsilon^2\|R\|^2_{L_{x,v}^\infty}\lesssim\delta.
  \end{equation*}
Recalling $G$ as defined in \eqref{G} and utilizing Theorem \ref{Thmcns}, we establish
  \begin{equation*}
    \left\|\frac{F-M}{\sqrt{\mu}}\right\|_{L^2_{x,v}}\leq \varepsilon\left\|\frac{G}{\sqrt{\mu}}\right\|_{L^2_{x,v}}+\varepsilon^2\left\|R\right\|_{L^2_{x,v}}\lesssim \varepsilon^2,
  \end{equation*}
  and
  \begin{equation*}
    \left\|w_k\frac{F-M}{\sqrt{\mu}}\right\|_{L^\infty_{x,v}}\leq \varepsilon\left\|\frac{G}{\sqrt{\mu}}\right\|_{L^\infty_{x,v}}+\varepsilon^2\left\|R\right\|_{L^\infty_{x,v}}\lesssim \varepsilon.
  \end{equation*}
  Thus, \eqref{main1} follows. The proof is complete.
\end{proof}

\appendix

\section{Some useful propositions}\label{Se_useful_estimate}
This section establishes key notations and presents essential lemmas for our analysis. The proof of uniform bound for the compressible Navier-Stokes system relies fundamentally on the use of conormal derivatives, which we define below along with relevant propositions. The conormal derivative coincides with the standard normal derivative in the interior of the domain, while at the boundary, it incorporates only the tangential components of the derivative. This approach was first introduced in \cite{masmoudi2012uniform} for studying the vanishing viscosity limit of incompressible Navier-Stokes system. Subsequent applications to compressible flows can be found in \cites{duan2021compressible,wang2016uniform,xiao2007vanishing}.
\begin{proposition}\label{Pro:div}Let $m\in \mathbb{N}_+$. Let $u\in H^m$ be a vector-valued function. Then we have
  \begin{equation}\label{div}
    \|u\|_{H^m} \lesssim \|\nabla_x\times u\|_{H^{m-1}}+\|\nabla_x\cdot u\|_{H^{m-1}}+|u\cdot n|_{H^{m-\frac{1}{2}}} ,
  \end{equation}
  and
  \begin{equation}\label{curl}
    \|u\|_{H^m} \lesssim  \|\nabla_x\times u\|_{H^{m-1}}+\|\nabla_x\cdot u\|_{H^{m-1}}+|n\times u|_{H^{m-\frac{1}{2}}},
  \end{equation}
  where $|\cdot|_{H^{m-\frac{1}{2}}}$ is the norm  in $H^{m-\frac{1}{2}}(\partial\Omega)$.
\end{proposition}
\begin{proof}
  See \cites{temam1977navier,xiao2007vanishing}.
\end{proof}

\begin{proposition}\label{Pro:infty}\cite[Proposition 2.2]{masmoudi2012uniform}
  Let $m_1\geq 0$ and $m_2\geq 0$ be integers, $f\in H^{m_1}_{co}(\Omega)$ and $\nabla_x f\in H^{m_2}_{co}.$ Then,
  \begin{equation}\label{infty}
    \|f\|^2_{L^\infty_x} \lesssim \|\nabla_x f\|_{H^{m_2}_{co}}\|f\|_{H^{m_1}_{co}},
  \end{equation}
provided $m_1+m_2\geq 3$, and
\begin{equation}\label{boundaryest}
  |f|^2_{H^s(\partial\Omega)} \lesssim \|\nabla_x f \|_{H^{m_2}_{co}}\|f\|_{H^{m_1}_{co}},
\end{equation}
provided $m_1+m_2\geq 2s\geq 0.$
\end{proposition}
\begin{proposition}\label{Pro:Agmon}(Agmon's inequality) If $f\in H^2(\Omega)$, then it holds that
\begin{equation}\label{agmon}
    \|f\|_{L_x^\infty}^2 \lesssim \| f\|_{H^1}\|  f\|_{H^2}.
\end{equation}
\end{proposition}

Furthermore, we give a Helmholtz decomposition result, which is essentially used in the study of the compressible Navier-Stokes system.
\begin{proposition}\label{Pro:Helm}
  Let $\Omega\in\mathbb{R}^3$ be open bounded domain with a smooth boundary. For every vector field  $f \in L^{2}\left(\Omega ; \mathbb{R}^{3}\right)$, there exist $ \varphi: \Omega \rightarrow \mathbb{R}$  and  $\psi: \Omega \rightarrow \mathbb{R}^{3}$  such that
\begin{equation*}
    f=\nabla_{x} \mathfrak{f}_1+\nabla_{x} \times \mathfrak{f}_2,  \mathfrak{f}_1 \in H^{1}(\Omega ; \mathbb{R}), \nabla_{x} \cdot\mathfrak{f}_2=0,\left.n \cdot\left(\nabla_{x} \times \mathfrak{f}_2\right)\right|_{\partial \Omega}=0.
\end{equation*}
Moreover, we have the following estimate
\begin{equation*}
    \|\mathfrak{f}_1\|_{H^{1}}+\|\mathfrak{f}_2\|_{H^{1}} \lesssim \|f\|_{L^2_x}.
\end{equation*}
\end{proposition}
The above two propositions are quite fundamental. For further details and proofs, we refer the reader to \cite{James-Josel-Sadowski-}.

\section{Elliptic estimates}\label{App.A}
In this appendix, we collect preliminary results from Agmon-Douglis-Nirenberg \cite{agmon1964estimates}.
\begin{lemma}\label{Lm:ADN}
   Suppose that $u,P $ are solutions of the generalized Stokes problem:
   \begin{equation}\label{sotkes}
\begin{cases}
        -\Delta u+\nabla_x P=f \quad &\text{in} \quad\Omega,\\
        \nabla_x\cdot u=g \quad &\text{in}\quad \Omega,\\
        u\cdot n=0 \quad &\text{on}\quad \partial\Omega,\\
        [S(u)\cdot n]^{\tan}=0\quad &\text{on}\quad \partial\Omega.
\end{cases}
  \end{equation}
   then we have
   \begin{equation}\label{ADN}
    \|u\|_{W^{2,p}}\lesssim \|f\|_2^2+\|g\|_{H^1}^2.
   \end{equation}
\end{lemma}

\begin{proof}
  This lemma is results from Agmon-Douglis-Nirenberg \cite{agmon1964estimates}.

   Let $u_4=P, u=(u_1,u_2,u_3,u_4),f=(f_1,\dots,f_n,g).$ Then we express the elliptic system \eqref{sotkes} in matrix form as:
\begin{equation}
    \sum_{j=1}^{4}l_{ij}(\partial)u_j=f_j,
\end{equation}
where  $l_{i j}(\xi)$, $\xi=\left(\xi_{1}, \ldots, \xi_{4}\right) \in \mathbb{R}^{3}$, is the matrix
$$l(\xi)=\begin{bmatrix}
 |\xi|^{2} &0 & 0&-\xi_1 \\
0 &  |\xi|^{2} & 0 &-\xi_2 \\
0& 0 &  |\xi|^{2}&-\xi_3 \\
 \xi_1  & \xi_2  & \xi_3  & 0
\end{bmatrix}.$$
We take (see  \cite[pp.38]{agmon1964estimates}),
 $ s_{i}=0, t_{i}=2,1 \leqslant i \leqslant n, s_{n+1}=-1,  t_{n+1}=1 $.
 As requested, degree $l_{i j}(\xi) \leqslant s_{i}+t_{j}, 1 \leqslant i, j \leqslant n+1$, and we have  $l_{i j}^{\prime}(\xi)= l_{i j}(\xi).$

We easily compute
$$ \operatorname{det}(l(\xi))=|\xi|^{6},$$
which ensures the ellipticity of the system. Let $n$ denote the normal vector and $\Lambda \neq 0 $ any tangent vector. Then $D(\Lambda+\tau n)$  has 3 roots in  $\tau$:
\begin{equation*}
   \tau_{1}=\tau_{2}=\tau_{3}=i|\Lambda|.
\end{equation*}
We define
\begin{equation}\label{M+}
   M^{+}(\Lambda, \tau)=(\tau-i|\Lambda|)^{3}.
\end{equation}
The adjoint matrix of  $l(\xi)$  is
\begin{equation}\label{Ladjoint}
   L(\xi)=\left(L_{i j}(\xi)\right)=\begin{bmatrix}
  \left(\xi_{2}^{2}+\xi_{3}^{2}\right)|\xi|^{2} & -\xi_{1} \xi_{2}|\xi|^{2} & -\xi_{1} \xi_{3}|\xi|^{2}&\xi_1|\xi|^4 \\
  -\xi_{1} \xi_{2}|\xi|^{2} &  \left(\xi_{1}^{2}+\xi_{3}^{2}\right)|\xi|^{2} & -\xi_{2} \xi_{3}|\xi|^{2} &\xi_2|\xi|^4 \\
  -\xi_{1} \xi_{3}|\xi|^{2} & -\xi_{2} \xi_{3}|\xi|^{2} &  \left(\xi_{1}^{2}+\xi_{2}^{2}\right)|\xi|^{2}&\xi_3|\xi|^4 \\
  -\xi_1|\xi|^4 &-\xi_2|\xi|^4 &-\xi_3|\xi|^4 & |\xi|^6
  \end{bmatrix}.
\end{equation}

Next we express boundary operator in \eqref{sotkes} in matrix form:
\begin{equation*}
   B(\xi)=\begin{bmatrix}
  \xi \cdot n+\xi_{1} n_{1}-2 n_{1}^{2} \xi \cdot n & n_{2} \xi_{1}-2 n_{1} n_{2} \xi \cdot n & n_{3} \xi_{1}-2 n_{1} n_{3} \xi \cdot n &0 \\
  n_{1} \xi_{2}-2 n_{2} n_{1} \xi \cdot n & \xi \cdot n+\xi_{2} n_{2}-2 n_{2}^{2} \xi \cdot n & n_{3} \xi_{2}-2 n_{2} n_{3} \xi \cdot n &0 \\
  n_{1} \xi_{3}-2 n_{3} n_{1} \xi \cdot n & n_{2} \xi_{3}-2 n_{3} n_{2} \xi \cdot n & \xi \cdot n+\xi_{3} n_{3}-2 n_{3}^{2} \xi \cdot n &0\\
  n_{1} & n_{2} & n_{3}&0
  \end{bmatrix}.
\end{equation*}
We set $r_{h}=-1$ for $h=1, \ldots, 3$ and $r_h=-2$ for $h=4$. Then, as requested, degree $B_{h j}\leq r_h+t_j $ and we have  $B_{h j}^{\prime}=B_{h j}.$
Since  $\Lambda \cdot n=0$, we have $ (\Lambda+\tau n) \cdot n=\tau$, and $(\Lambda+\tau n)_{i} n_{j}=\Lambda_{i} n_{j}+\tau n_{i} n_{j} $, hence
$$
B(\Lambda+\tau n)=\begin{bmatrix}
  \tau\left(1-n_{1}^{2}\right)+\Lambda_{1} n_{1} & \Lambda_{1} n_{2}-\tau n_{1} n_{2} & \Lambda_{1} n_{3}-\tau n_{1} n_{3} \\
  \Lambda_{2} n_{1}-\tau n_{1} n_{2} & \tau\left(1-n_{2}^{2}\right)+\Lambda_{2} n_{2} & \Lambda_{2} n_{3}-\tau n_{2} n_{3} \\
  \Lambda_{3} n_{1}-\tau n_{1} n_{3} & \Lambda_{3} n_{2}-\tau n_{2} n_{3} & \tau\left(1-n_{3}^{2}\right)+\Lambda_{3} n_{3} \\
  n_{1} & n_{2} & n_{3}
  \end{bmatrix}.$$

By direct computation, the determinant of the submatrix that consists of the first three rows in  $B(\Lambda+\tau n)$  is $0$, indicating that the rank of the first three rows is $2$.

Without loss of generality, assuming the first two rows are linearly independent, we claim $ n_{3} \neq 0 $. Suppose $ n_{3}=0 $, in such a case, we would have $ n_{1}^{2}+n_{2}^{2}=1, \Lambda_{1} n_{1}+\Lambda_{2} n_{2}=0$. The first two rows would then transform to
\begin{equation}
  \left[\begin{array}{ccc}
  \tau n_{2}^{2}-\Lambda_{2} n_{2} & \Lambda_{1} n_{2}-\tau n_{1} n_{2} & 0 \\
  \Lambda_{2} n_{1}-\tau n_{1} n_{2} & \tau n_{1}^{2}-\Lambda_{1} n_{1} & 0
  \end{array}\right].
\end{equation}
These row vectors  $n_{2}\left(\tau n_{2}-\Lambda_{2}, \Lambda_{1}-\tau n_{1}\right)$ and $n_{1}\left(\Lambda_{2}-\tau n_{2}, \tau n_{1}-\Lambda_{1}\right)$,  are linearly dependent, leading to a contradiction. Therefore, $n_{3} \neq 0 $. We rewrite the matrix  $\mathfrak{B}$  as
\begin{equation}\label{B}
   \mathfrak{B}(\Lambda+\tau n)=\begin{bmatrix}
  \tau\left(1-n_{1}^{2}\right)+\Lambda_{1} n_{1} & \Lambda_{1} n_{2}-\tau n_{1} n_{2} & \Lambda_{1} n_{3}-\tau n_{1} n_{3} & 0\\
  \Lambda_{2} n_{1}-\tau n_{1} n_{2} & \tau\left(1-n_{2}^{2}\right)+\Lambda_{2} n_{2} & \Lambda_{2} n_{3}-\tau n_{2} n_{3}& 0 \\
  n_{1} & n_{2} & n_{3}& 0
  \end{bmatrix}.
\end{equation}
Next, to verify the Complementing Condition, we compute the matrix multiplication with modulo  $M^{+}(\Lambda, \tau)$  defined in \eqref{M+}:
\begin{equation}\label{BL}
   \mathfrak{B}(\Lambda+\tau n) L(\Lambda+\tau n) \quad \bmod\,M^{+}(\Lambda, \tau).
\end{equation}

First, we compute $ L(\Lambda+\tau n)\quad \bmod\,M^{+}(\Lambda, \tau)$. Note that the elements in \eqref{Ladjoint} involve $|\xi|^{4}$, $\xi_j|\xi|^4$ and $\xi_{j}\xi_{k}|\xi|^{2}$. By direct computation, we obtain:
\begin{equation}\label{a}
  |\Lambda+\tau n|^{4} \quad \bmod\, (\tau-i|\Lambda|)^{3}=-4|\Lambda|^{2}(\tau-i|\Lambda|)^{2},
\end{equation}
\begin{equation}\label{c}
  |\Lambda+\tau n|^4(\Lambda_j+\tau n_j) \quad \bmod\,  (\tau-i|\Lambda|)^{3}=-4|\Lambda|^{2}\left(\Lambda_j+i|\Lambda|n_j\right):=\mathcal{A}_j,
\end{equation}
and
\begin{equation}\label{b}
  \begin{aligned}
 &\left(\Lambda_{j}+\tau n_{j}\right)\left(\Lambda_{k}+\tau n_{k}\right)|\Lambda+\tau n|^{2} \quad \bmod\, (\tau-i|\Lambda|)^{3} \\
 &\quad=(\tau-i|\Lambda|)\left[\left(\Lambda_{j} \Lambda_{k}+3 i|\Lambda| \Lambda_{j} n_{k}+3 i|\Lambda| \Lambda_{k} n_{j}-5 n_{j} n_{k}|\Lambda|^{2}\right) \tau\right. \\
 &\qquad\left.+i \Lambda_{j} \Lambda_{k}|\Lambda|+n_{j} \Lambda_{k}|\Lambda|^{2}+\Lambda_{j} n_{k}|\Lambda|^{2}+3 i n_{j} n_{k}|\Lambda|^{3}\right]
 :=(\tau-i|\Lambda|)\mathcal{B}_{jk}.
\end{aligned}
\end{equation}
Then we have
\begin{equation}\label{d}
  \begin{array}{l}
  L(\Lambda+\tau n) \quad \bmod\, M^{+}(\Lambda, \tau) \\
  =-(\tau-i|\Lambda|)
\begin{bmatrix}
    \mathcal{B}_{22}+ \mathcal{B}_{33} & -\mathcal{B}_{12}&-\mathcal{B}_{13}&\mathcal{A}_1 \\
    -\mathcal{B}_{12} & \mathcal{B}_{11}+ \mathcal{B}_{33} & -\mathcal{B}_{23}&\mathcal{A}_2 \\
    -\mathcal{B}_{13} & -\mathcal{B}_{23} & \mathcal{B}_{11}+ \mathcal{B}_{22}&\mathcal{A}_3\\
    -\mathcal{A}_1&-\mathcal{A}_2&-\mathcal{A}_3&0
\end{bmatrix}.
  \end{array}
\end{equation}
To compute the matrix multiplication \eqref{BL}, due to the extra  $\tau $ factor in \eqref{B}, we apply \eqref{a} and \eqref{b} to further compute

\begin{equation}
  \begin{array}{l}
  \tau|\Lambda+\tau n|^{4} \quad \bmod\, (\tau-i|\Lambda|)^{3} \\
  =-4|\Lambda|^{2}(\tau-i|\Lambda|)^{2} \tau \quad \bmod\, (\tau-i|\Lambda|)^{3}=-4 i|\Lambda|^{3}(\tau-i|\Lambda|)^{2}.
  \end{array}
\end{equation}

\begin{equation}
  \begin{array}{l}
  \tau\left(\Lambda_{j}+\tau n_{j}\right)\left(\Lambda_{k}+\tau n_{k}\right)|\Lambda+\tau n|^{2}\quad \bmod\, \left(\tau-i|\Lambda|^{3}\right) \\
  =(\tau-i|\Lambda|) \tau \mathcal{B}_{j k} \quad \bmod\,\left(\tau-i|\Lambda|^{3}\right) \\
  =(\tau-i|\Lambda|)\left[\left(3 i|\Lambda| \Lambda_{j} \Lambda_{k}-5 n_{j} \Lambda_{k}|\Lambda|^{2}-5 \Lambda_{j} n_{k}|\Lambda|^{2}-7 i n_{j} n_{k}|\Lambda|^{3}\right) \tau\right. \\
  \left.+\Lambda_{j} \Lambda_{k}|\Lambda|^{2}+3 i|\Lambda|^{3} n_{j} \Lambda_{k}+3 i|\Lambda|^{3} \Lambda_{j} n_{k}-5 n_{j} n_{k}|\Lambda|^{4}\right].
  \end{array}
\end{equation}
We denote
\begin{equation}
  \begin{aligned}
  \mathcal{C}_{j k}:= & \left(3 i|\Lambda| \Lambda_{j} \Lambda_{k}-5 n_{j} \Lambda_{k}|\Lambda|^{2}-5 \Lambda_{j} n_{k}|\Lambda|^{2}-7 i n_{j} n_{k}|\Lambda|^{3}\right) \tau \\
  & +\Lambda_{j} \Lambda_{k}|\Lambda|^{2}+3 i|\Lambda|^{3} n_{j} \Lambda_{k}+3 i|\Lambda|^{3} \Lambda_{j} n_{k}-5 n_{j} n_{k}|\Lambda|^{4}.
  \end{aligned}
\end{equation}
Following similar calculations as in \cite{chen2023macroscopic}, we finally conclude that
\begin{equation*}
  \mathfrak{B}(\Lambda+\tau n)L(\Lambda+\tau n) \quad\text{ mod }(\tau-i|\Lambda|)^3=-(\tau-i|\Lambda|)M(\Lambda+\tau n),
\end{equation*}
with
\begin{equation}
  \det\left(M(\Lambda+\tau n)\right)=32|\Lambda|^8n_3(\tau-i|\Lambda|)^3\neq 0,
\end{equation}
where we have used the fact that $|\Lambda|\neq 0$ and $n_3\neq 0$.
Thus, all three rows of \eqref{BL} are linearly independent and the Complementing Condition holds. We then apply \cite[Theorem 2.5,pp.78]{agmon1964estimates} to conclude the proof of the lemma.
\end{proof}

\section*{Acknowledgments}
The author N. Jiang is supported by the grants from the National Natural Science Foundation of China under contract Nos. 11971360, 12371224 and 12221001. We appreciate any comments and suggestions from the anonymous reviewers.\\

{\bf Data availability} The authors declare that the data supporting the findings of this study are available within the paper.\\

{\bf Conflict of interest} The authors have no Conflict of interest.

\bibliography{reff} 
\bibliographystyle{abbrv}
\end{sloppypar}
\end{document}